\documentclass[final,1p,times]{elsarticle}

\usepackage{latexsym, amsmath, amssymb}
\usepackage[ruled,lined]{algorithm2e}
\usepackage{array}

\usepackage[colorlinks,linkcolor=blue, anchorcolor=blue, citecolor=blue]{hyperref}
\usepackage{color,tabularx}
\usepackage{setspace}
\usepackage{booktabs}
\usepackage{subfigure}
\usepackage{multirow}
\usepackage{lineno}
\usepackage{stmaryrd}
\usepackage[T1]{fontenc}
\biboptions{numbers,sort&compress}
\usepackage{graphicx} 

\allowdisplaybreaks

\newcommand{\p}{\partial}

\numberwithin{equation}{section}
\newtheorem{thm}{Theorem}[section]
\newtheorem{lem}[thm]{Lemma}

\newtheorem{rem}[thm]{Remark}
\newtheorem{cor}[thm]{Corollary}
\newtheorem{example}{Example}[section]

\journal{}
\begin{document}

\begin{frontmatter}

\title{ An efficient  two-grid fourth-order compact difference scheme with variable-step BDF2 method for the semilinear parabolic equation}

\author[1]{Bingyin Zhang} \ead{zhangbingyin@stu.ouc.edu.cn}
\author[1]{Hongfei Fu \corref{cor1}} \ead{fhf@ouc.edu.cn}
\address[1]{School of Mathematical Sciences, Ocean University of China, Qingdao, Shandong 266100, China}
\cortext[cor1]{Corresponding author.}

\begin{abstract}
Due to the lack of  corresponding analysis on appropriate mapping operator between two grids, high-order two-grid difference algorithms are rarely studied. In this paper, we firstly discuss the boundedness of a local bi-cubic Lagrange interpolation operator. And then, taking the semilinear parabolic equation as an example, we first construct a  variable-step high-order nonlinear difference algorithm using compact difference technique in space and  the second-order backward differentiation formula (BDF2) with variable temporal stepsize in time. With the help of discrete orthogonal convolution (DOC) kernels and a cut-off numerical technique, the unique solvability and corresponding error estimates of the high-order nonlinear difference scheme are established under assumptions that the temporal stepsize ratio satisfies $ r_{k} := \tau_{k}/\tau_{k-1} < 4.8645 $ and  the maximum temporal stepsize satisfies $ \tau = o( h^{\frac{1}{2}} ) $.
Then, an efficient two-grid high-order difference algorithm is developed by combining a small-scale variable-step high-order nonlinear difference algorithm on the coarse grid and a large-scale variable-step high-order linearized difference algorithm on the fine grid, in which the constructed piecewise bi-cubic Lagrange interpolation mapping operator is adopted to project the coarse-grid solution to the fine grid. Under the same temporal stepsize ratio restriction $ r_{k} < 4.8645 $ and  a weaker maximum temporal stepsize condition $ \tau = o(H^{\frac{1}{2}} ) $, optimal fourth-order in space and second-order in time error estimates of the two-grid difference scheme is established if the coarse-fine grid stepsizes satisfy $ H=\mathcal{O}(h^{4/7})$.
Finally, several numerical experiments are carried out to demonstrate the effectiveness and efficiency of the proposed scheme.
\end{abstract}


\end{frontmatter}

\section{Introduction}\label{sec:Intro}
It is well known that the analytical solutions of nonlinear parabolic PDEs arising from a variety of physical and engineering applications are not available in most cases. Thus, numerous efforts have been devoted to the development of efficient numerical schemes, see \cite{1999_NM_Akrivis,2017_JSC_Tang,1995_SINUM_Xu,SINUM_Dawson_1998,2002_CMAME_Novo,2020_SINUM_Li}. 
Generally speaking, fully-implicit numerical schemes are usually proved to be unconditionally stable. Unfortunately, at each time step, one has to solve a system of nonlinear equations \cite{2017_JSC_Tang,2013_IJNAM_Li}, in which an extra iterative process must be imposed, and this in turn may cause severe computational costs. Instead, a very popular and widely-used approach is the so-called implicit-explicit scheme, which treats the linear term implicitly and the nonlinear term explicitly. However, if the corresponding globally continuous condition of the nonlinear term (e.g., $ \vert f''(w) \vert \leq K_{f}$ for $ w \in \mathbb{R}$) cannot be imposed or the boundedness of numerical solution in $L^{\infty}$ norm cannot be obtained, this method usually suffers from a very restrictive temporal stepsize condition caused by use of inverse inequality for convergence, e.g., $ \tau = \mathcal{O}( h^{d/2p} ) $, where $h$ is the spatial mesh size, $d$ is the space dimension and $ p $ is the accuracy of the time discretization. Therefore, such restrictions may lead to the use of a small temporal stepsize, and thus much computational time may be consumed. Here we refer to \cite{1999_NM_Akrivis,2005_BIT_Akrivis,2014_SISC_Sun,2008_SINUM_Wang} for an incomplete list of references.

Another efficient and powerful strategy is the two-grid method which is proposed by  Xu et al. \cite{1995_SINUM_Xu,SINUM_Xu_1996}. The basic idea of this kind of method is to reduce the solution of the large-scale nonlinear problem on the fine grid to a small-scale nonlinear problem on the coarse grid and a large-scale linear problem on the fine grid. Hence, basically it  includes two solution steps:  First, one solve a nonlinear system on the coarse grid to obtain a rough approximation, and then solve a linearized system resulting from the rough solution to derive a corrected solution on the fine grid. Up to now, this technique has been widely applied to numerically solve many types of nonlinear PDEs, e.g., \cite{SINUM_Dawson_1998,2019_ANM_Hou,2021_JSC_Chen} for parabolic equations, \cite{2015_SINUM_Rui,2018_JSC_Liu} for Darcy-Forchheimer equations, \cite{2018_JSC_Novo} for Navier-Stokes equations, \cite{JSC_Li_2017,2022_JCM_Li,2017_JCP_Liu} for time-fractional equations and \cite{2019_ACM_Chen,2022_ANM_Fu} for other nonlinear equations. Compared to traditional implicit-explicit scheme, the main advantages of the two-grid scheme are twofold: (i) if the boundedness of numerical solution in $L^{\infty}$ norm is not obtained, the temporal stepsize restriction of the method under the local Lipschitz continuous condition on nonlinear term usually is only related to $ H $ instead of $h$, which is much more weaker, see \cite{SINUM_Dawson_1998,2022_ANM_Fu}; 
(ii) the two-grid method which treats the  nonlinearity on the coarse grid and solving the  linear system on the fine grid \cite{1995_SINUM_Xu}, is much more stable and accurate than the implicit-explicit one when solving the nonlinear PDEs whose solutions change rapidly with respect to time.
In this paper, we will investigate the numerical stability of these two methods for general semilinear parabolic PDEs by carrying out representative numerical examples.

Unlike the finite element method, which generates pointwise solution in space and thus is easy to develop two-grid finite element method \cite{SINUM_Xu_1996,2003_IJNME_Chen,2018_JSC_Novo}, the solution yielded by finite difference method is only on discrete grids, and therefore, an appropriate accuracy-preserving mapping operator from the coarse-grid function space to the fine-grid function space is required to construct and analyze the two-grid difference method, e.g., piecewise linear/bilinear interpolation for  second-order two-grid difference schemes \cite{SINUM_Dawson_1998,2015_SINUM_Rui,JSC_Li_2017,2022_ANM_Fu}. However, due to  the lack of corresponding analysis on appropriate high-order mapping operator, the work about high-order two-grid difference method is meager and, in fact, numerical analysis is also lack. 
This motivates us to develop high-order two-grid difference scheme for semilinear parabolic PDEs with general boundary condition, e.g., Dirichlet or periodic boundary condition, by introducing and analyzing appropriate mapping operator. 

In this paper, to illustrate the application of the proposed high-order two-grid difference method, we focus on the following semilinear parabolic equation
\begin{equation}\label{Model:Bihar}
	\begin{aligned}
		u_{t}(x,y,t) - c \Delta u(x,y,t) = f(u(x,y,t)) + g(x,y,t), \quad (x,y) \in \Omega, \  t \in (0,T],
	\end{aligned}
\end{equation}
where $\Omega = (0,L_x)\times (0,L_y)$,
subject to the initial condition
\begin{equation}\label{Initial:Bihar}
	\begin{aligned}
		u(x,y,0) = u_{0}(x,y), \quad (x,y) \in \Omega,
	\end{aligned}
\end{equation}
and periodic boundary condition or Dirichlet boundary condition
\begin{equation}\label{Bound:Bihar}
	\begin{aligned}
		u(x,y,t) = \psi (x,y,t), \quad (x,y) \in \p \Omega, \  t\in [0,T],
	\end{aligned}
\end{equation}
where $ c > 0 $ and $ \Delta $ is the Laplacian operator $ \p^{2}_{x} + \p^{2}_{y} $, $ f(u) $ is the nonlinear term, $ u_{0} $, $ \psi $ and $ g $ are given smooth functions.

For many time-dependent PDEs,  e.g., Allen-Cahn equations \cite{SINUM_Liao_2020}, whose solutions admiting multiple time scales, adaptive temporal stepsize  strategies \cite{SISC_Qiao_2011,MOC_Li_2019} are heuristic and available methods to improve accuracy or efficiency. 
Due to its strong stability, variable-step BDF2 method is practically valuable for stiff or differential-algebraic problems \cite{1974_SINUM_Gear,1997_SISC_Shampine}. However, compared to those one-step methods, such as the backward Euler and
Crank-Nicolson schemes, the numerical analysis of nonuniform BDF2 method would be challenging. In particular, for a linear parabolic problem, \cite{1998_BIT_Becker} proved that, if $ 0 < r_{k}:= \tau_{k}/\tau_{k-1}  \leq 1.868 $ with $ \tau_{k} := t_{k} - t_{k-1} $ the $k$th temporal stepsize, the variable-step BDF2 scheme is zero-stable and second-order convergence containing a prefactor $ \exp( C \Gamma_{n} ) $, where $ \Gamma_{n} := \sum^{n-2}_{k} [ r_{k} - r_{k+2} ]_{+} $ with $ [x]_{+} $ the positive part of $ x $. Recently, by introducing a generalized discrete Gr\"{o}nwall inequality, Chen et al. \cite{2019_SINUM_Chen} circumvented such a prefactor in error analysis under a little stronger step-ratio restriction $ 0 < r_{k} \leq 1.53 $. In \cite{2019_SINUM_Wang}, the authors developed an implicit-explicit BDF2 method with variable stepsize for the parabolic partial integro-differential equations and proved its stability and convergence with $ 0 < r_{k} \leq 1.91 $. The authors in \cite{SINUM_Liao_2020} considered the fully-implicit BDF2 scheme for the Allen--Cahn equation and established the maximum-norm stability under $ r_{k} < 1 + \sqrt{2} $ by developing a novel kernel recombination and complementary technique. To analyze the variable-step BDF2 scheme for linear reaction-diffusion equations, Liao and Zhang \cite{2021_MOC_Liao} introduced a new concept, namely, discrete orthogonal convolution (DOC) kernels, and they improved the unconditional stability in the $ L^{2} $ norm to $ r_{k} \leq 3.561 $. Subsequently, with the help of DOC kernels and corresponding convolution inequalities, there is a great progress on the stability and error estimates of variable-step BDF2 method for  nonlinear PDEs under  $ r_{k} < 3.561 $ \cite{2021_SCM_Liao,2022_IMA_Liao,2022_CNSNS_Zhao} and the further improved step-ratio restriction $ r_{k} < 4.8645 $ \cite{JSC_Zhang_2022,JSC_Liao_2022,2023_MCS_Li}, respectively.

Among all the variable-step BDF2 methods for nonlinear PDEs in the literature mentioned above, they treat the nonlinear terms fully or partially implicit, in which a nonlinear iteration must be implemented at each time step. Very recently, Zhao et al. \cite{2023_CMS_Zhang} presented a linearized variable-step BDF2 scheme for solving nonlinear parabolic equation, and they proved the unconditional error estimate under $ r_{k} < 4.8645 $ and the maximum temporal stepsize $ \tau \leq C \frac{1}{\sqrt{N}} $ by adopting the error splitting approach. After then, Li et al. extended this method to solve a nonlinear Ginzburg-Landau equation \cite{2023_CNSNS_Li} and coupled Ginzburg-Landau equations \cite{2023_ANM_Li} under the same conditions. In \cite{2023_JSC_Chen}, a positivity-preserving and energy stable BDF2 scheme with variable stepsize was developed for the Cahn--Hilliard equation with nonlinear logarithmic potential, and convergence analysis in $ L^{2} $ norm was established under $ \tau \leq Ch $. For solving gradient flow problems, Hou and Qiao \cite{2023_JSC_Qiao} proposed an unconditionally energy stable implicit-explicit BDF2 scheme with variable temporal stepsize using the SAV method, and derived the error estimates under the mild restriction on the adjacent temporal stepsize ratio $ r_{k} < 4.8645 $. Our goal  is to construct and analyze an efficient  high-order two-grid difference algorithm with nonuniform BDF2 method for the nonlinear parabolic equation \eqref{Model:Bihar}--\eqref{Bound:Bihar}. Compared to the existing literature, our contributions are mainly  threefold:
\begin{itemize}
	\item An efficient  variable-step two-grid fourth-order compact difference method is proposed,  by using compact difference scheme, two-grid method,  variable-step BDF2 formula as well as a developed piecewise bi-cubic Lagrange interpolation operator.
	
	\item Under the local continuous condition imposed on the nonlinear term, see \eqref{Regu:Bihar:e1}, by adopting the DOC kernels and a cut-off technique, we rigorously prove the unique solvability and error estimate for the nonlinear compact difference scheme, under the temporal stepsize ratio restriction $ r_{k} < 4.8645 $ and the maximum temporal stepsize condition $ \tau = o(h^{\frac{1}{2}}) $. Furthermore,  by discussing the boundedness of the proposed piecewise bi-cubic Lagrange interpolation operator with  periodic and Dirichlet boundary conditions, optimal error estimate of the two-grid compact difference scheme is established, under a weaker maximum temporal stepsize condition $ \tau = o(H^{\frac{1}{2}}) $ and  coarse-fine-grid condition $ H=\mathcal{O}(h^{4/7})$. 
	
	\item Several numerical experiments are presented to illustrate the effectiveness and efficiency of the proposed variable-step (adaptive) two-grid compact difference method, and comparisons of computational efficiency and stability with other schemes are also tested.
\end{itemize}

The remainder of this paper is organized as follows. In Section \ref{sec:Map_lag}, we introduce and analyze the high-order mapping operator between the coarse-grid function space and fine-grid function space. In Section \ref{sec:Non}, we first propose a nonlinear compact difference scheme with variable-step BDF2 method for  the semilinear parabolic equation subject to  Dirichlet boundary condition, and then rigorously prove the unique solvability and convergence, based on which we construct an efficient variable-step two-grid compact difference scheme in Section \ref{sec:TG}, and optimal-order error analysis under a weaker  maximum temporal stepsize condition and  a coarse-fine-grid condition is derived.  Moreover, in Section \ref{sec:PBC}, the developed methods and techniques are extended to the context of periodic boundary condition. Several numerical experiments are presented to demonstrate the accuracy and efficiency of the proposed method in Section \ref{sec:exam}. Finally, some concluding remarks are drawn in the last section.

\section{High-order mapping operator between two grids}\label{sec:Map_lag}
In this section, we first propose and analyze a high-order mapping operator  between two grids based on Lagrange interpolation, which plays a significant role in the construction and numerical analysis of high-order two-grid difference method in the subsequent sections.

\subsection{Some notations}
Given two positive integers $ N^{H}_{x} $ and $ N^{H}_{y} $, we define a uniform coarse grid $(x^{H}_{i},   y^{H}_{j}) := (iH_{x},   jH_{y}) $  for   $0 \leq i \leq N^{H}_{x}$ and $0 \leq j \leq N^{H}_{y} $, with corresponding coarse mesh sizes $ H_{x} := L_{x}/N^{H}_{x} $ and $ H_{y} := L_{y}/N^{H}_{y} $. Moreover, for fixed positive integers $ M_{x} ,  M_{y} \geq 2$, denote $ N^{h}_{x} := M_x N^{H}_{x} $, $ N^{h}_{y} := M_y N^{H}_{y} $ and define a uniform fine grid $ ( x^{h}_{i},   y^{h}_{j} ) := (ih_{x},  jh_{y} )$  for   $0 \leq i \leq N^{h}_{x}$ and $0 \leq j \leq N^{h}_{y} $, with corresponding fine mesh sizes $ h_{x} := H_{x}/M_{x} $ and $ h_{y} := H_{y}/M_{y} $. Denote $ H := \max\{ H_{x}, H_{y} \} $ and $ h := \max\{ h_{x}, h_{y} \}$. 

Let  $ \bar{\omega}_{\kappa} := \left\{ ( i, j ) \mid 0 \leq i \leq N^{\kappa}_{x}, \ 0 \leq j \leq N^{\kappa}_{y} \right\} $, 
$ \omega_{\kappa} := \bar{\omega}_{\kappa} \cap \Omega $ and $ \partial \omega_{\kappa} := \bar{\omega}_{\kappa} \cap \partial \Omega $ denote the sets of spatial grids, where $ \kappa = H$ or $ h $.
 Accordingly, we define the following discrete spaces of grid functions
$$
\mathcal{V}_{\kappa} = \left\{v = \{ v_{i,j}\}  \mid (i,j) \in \bar{\omega}_{\kappa} \right\} \quad  {\rm and } \quad \mathcal{V}^{0}_{\kappa} = \left\{v \mid v \in \mathcal{V}_{\kappa} \  \  \text {and} \  \  v_{i,j} = 0 \  {\rm if} \  (i, j) \in \partial \omega_{\kappa} \right\}.
$$
For any $w, q \in \mathcal{V}_{\kappa}$, we introduce the following notations
\begin{equation*}
	\begin{array}{l}
		\displaystyle d_{\kappa,x} w_{i+\frac{1}{2},j} := \frac{1}{\kappa_{x}} \left(w_{i+1,j}-w_{i,j}\right),
		\quad \displaystyle  d^{2}_{\kappa,{x}} w_{i,j} := \frac{1}{\kappa_{x}} \left(\left[d_{\kappa,x} w\right]_{i+\frac{1}{2},j} - \left[d_{\kappa,x} w\right]_{i-\frac{1}{2},j}\right),
	\end{array}
\end{equation*}
\begin{equation*}
	\begin{array}{ll}
		\mathcal{A}_{\kappa,x} w_{i,j} :=
			\left\{
			\begin{split}
				& w_{i,j} + \frac{\kappa_{x}^2}{12} d^{2}_{\kappa,x} w_{i,j} = \frac{1}{12} ( w_{i-1,j} + 10w_{i,j} + w_{i+1,j} ),  &\quad 1\leq i \leq N^{\kappa}_{x}-1,\\
				& w_{i,j} , & \quad   i=0, N^{\kappa}_{x}.	
			\end{split}
		\right.
	\end{array}
\end{equation*}
Similarly,  the notations $ d_{\kappa,y} w_{i,j+\frac{1}{2}} $, $ d^{2}_{\kappa,y} w_{i,j} $ and $ \mathcal{A}_{\kappa,y} w_{i,j} $ can be defined. Furthermore, we denote $ \Delta_{\kappa} := d^{2}_{\kappa,x} + d^{2}_{\kappa,y} $, $ \Lambda_{\kappa} := \mathcal{A}_{\kappa,x} d^{2}_{\kappa,y} + \mathcal{A}_{\kappa,y} d^{2}_{\kappa,x} $ and $ \mathcal{A}_{\kappa} := \mathcal{A}_{\kappa,x} \mathcal{A}_{\kappa,y} $.

Besides, we also introduce  the discrete inner products 
\[
	(w,q)_{\kappa} = \kappa_{x} \kappa_{y} \sum^{N^{\kappa}_{x}-1}_{i=1} \sum^{N^{\kappa}_{y}-1}_{j=1} w_{i,j} q_{i,j}, 
\]
\[	
	(w, q)_{\kappa,x} = \kappa_{x} \kappa_{y} \sum^{N^{\kappa}_{x}-1}_{i=0} \sum^{N^{\kappa}_{y}-1}_{j=1} w_{i+\frac{1}{2},j} q_{i+\frac{1}{2},j}, ~~
	 (w, q)_{\kappa,y} = \kappa_{x} \kappa_{y} \sum^{N^{\kappa}_{x}-1}_{i=1} \sum^{N^{\kappa}_{y}-1}_{j=0}  w_{i,j+\frac{1}{2}} q_{i,j+\frac{1}{2}},
\]
\[
		<w,q>_{\kappa}  =   \kappa_{x} \kappa_{y} \left[ \frac{1}{4}\sum_{i=\{0,N^{x}_{\kappa}\}} \sum_{j=\{0,N^{y}_{\kappa}\}}  
	 + \frac{1}{2}    \sum^{N^{x}_{\kappa}-1}_{i=1}\sum_{j=\{0,N^{y}_{\kappa}\}}   
	  +  \frac{1}{2}\sum_{i=\{0,N^{x}_{\kappa}\}} \sum^{N^{y}_{\kappa}-1}_{j=1}  
		+ \sum^{N^{x}_{\kappa}-1}_{i=1} \sum^{N^{y}_{\kappa}-1}_{j=1} \right] w_{i,j} q_{i,j},
\]
and corresponding discrete $L^{2}$ and $L^{\infty}$  norms
$$
\| w \|_{\kappa} = \sqrt{(w, w)_{\kappa}}, \quad \interleave w \interleave_{\kappa} = \sqrt{<w, w>_{\kappa}}, \quad \|w\|_{\mathcal{A},\kappa} = \sqrt{ ( \mathcal{A}_{\kappa} w, w )_{\kappa} }, \quad \|w\|_{\kappa,\infty} = \max_{(i,j) \in \bar{\omega}_{\kappa} } |w_{i,j}|.
$$
It is easy to check that $
\| \mathcal{A}_{\kappa} w \|_{\kappa} \leq \interleave w \interleave_{\kappa} 
$ for any $ w \in \mathcal{V}_{\kappa} $. Moreover, two well-known and useful lemmas are listed below.

\begin{lem}[\cite{Compact_S12}]\label{lem:cpmpact1} For any $ w \in \mathcal{V}^{0}_{\kappa} $, we have
	$
	\frac{1}{3} \| w \|_{\kappa}^{2} \leq \| w \|_{\mathcal{A},\kappa}^{2} \leq \| w \|_{\kappa}^{2}.
	$
\end{lem}

\begin{lem}[\cite{Compact_S12}]\label{lem:inverse} For any $ w \in \mathcal{V}_{\kappa} $, there exists a positive constant $ C_{0} $, independent of $ \kappa $, such that
	$$
	\| w \|_{\kappa,\infty} \leq  C_{0}\ \kappa^{-1}\| w \|_{\kappa}. 
	$$
\end{lem}

\subsection{Piecewise bi-cubic Lagrange interpolation}
An important tool used in the construction of high-order two-grid method is the local high-order Lagrange interpolation from coarse-grid space to fine-grid space. We shall present and discuss its properties in this subsection.

We first define the one-dimensional piecewise cubic Lagrange interpolation  along $x$-direction. For  each $ x \in ( x^{H}_{i}, x^{H}_{i+1}) $ with $ 0 \leq i \leq N^{H}_{x}-1 $, we use $ \big\{ \phi^{x}_{i,s}(x) \big\}^{3}_{s=0} $ to represent  the  cubic Lagrange interpolation basis functions. For $ 1 \leq i \leq N^{H}_{x}-2 $, $ \phi^{x}_{i,s}(x) $ is defined as
\begin{equation}\label{Def:Basis_2}
	\begin{array}{l}
		\phi^{x}_{i,s} (x) :=
		\left\{
		\begin{split}
			& -\frac{ (x - x^{H}_{i}) (x - x^{H}_{i+1}) (x - x^{H}_{i+2}) }{6H^{3}_{x}}, \qquad s = 0,\\
			& \frac{ (x - x^{H}_{i-1}) (x - x^{H}_{i+1}) (x - x^{H}_{i+2}) }{2H^{3}_{x}}, \qquad \  \  \ s = 1, \\
			& - \frac{ (x - x^{H}_{i-1}) (x - x^{H}_{i}) (x - x^{H}_{i+2}) }{2H^{3}_{x}}, \qquad \ s = 2, \\
			& \frac{ (x - x^{H}_{i-1}) (x - x^{H}_{i}) (x - x^{H}_{i+1}) }{6H^{3}_{x}}, \qquad \quad  \ s = 3.	
		\end{split}
		\right.
	\end{array}
\end{equation}
For $i=0$, i.e., $ x \in (x^{H}_{0}, x^{H}_{1}) $, we define $ \phi^{x}_{0,s}(x) := \phi^{x}_{1,s}(x) $; and for $i=N^{H}_{x}-1$, i.e., $ x \in ( x^{H}_{N^{H}_{x}-1}, x^{H}_{N^{H}_{x}} ) $, we define $ \phi^{x}_{N^{H}_{x}-1,s}(x) := \phi^{x}_{N^{H}_{x}-2,s}(x)$. 
Then,  for any continuous function $w(x)$,  the piecewise cubic Lagrange interpolation operator $\Pi_{H,x} $ along  $x$-direction is defined as
\begin{equation}\label{Def:cubic_L}
			\Pi_{H,x} w(x) :=
			\left\{
			\begin{aligned}
					& \sum^{3}_{s=0} w_{s}\ \phi^{x}_{0,s}(x),     & \quad x \in ( x^{H}_{0},x^{H}_{1} ),  \  i=0,\\
					& \sum^{3}_{s=0} w_{i-1+s}\ \phi^{x}_{i,s}(x),  & \quad x \in ( x^{H}_{i},x^{H}_{i+1} ), \  1 \leq i \leq N^{H}_{x} - 2,\\
					& \sum^{3}_{s=0} w_{N^{H}_{x}-3+s}\ \phi^{x}_{N^{H}_{x}-1,s}(x),  &\quad   x \in ( x^{H}_{N^{H}_{x}-1},x^{H}_{N^{H}_{x}}), \  i=N^{H}_{x}-1.
				\end{aligned}
			\right.
\end{equation}
where $w_{i}=w(x^{H}_{i})$ for $ 0 \leq i\leq N^{H}_{x} $.

Similarly, we can define the  cubic Lagrange interpolation basis functions $ \big\{\phi^{y}_{j,s}(y)\big\}^{3}_{s=0} $  and corresponding piecewise cubic Lagrange interpolation operator  $ \Pi_{H,y} $ along  $y$-direction.  Therefore, the piecewise bi-cubic Lagrange interpolation operator $\Pi_{H}$ can be defined as the tensor product of the one-dimensional piecewise cubic Lagrange interpolation operators  in two directions, that is, $\Pi_{H}:=\Pi_{H,y}\Pi_{H,x}$. When no confusion caused, below we denote $(\Pi_{H}w)_{i,j}=\Pi_{H}w_{i,j}$ for $(i,j) \in \bar{\omega}_{h}$ and $w\in\mathcal{V}_H$.

\begin{lem}[\cite{Cubic}]\label{lem:cubic_error} Assume that $ w \in W^{4,\infty} (\Omega) $, there exists positive constants $ C_{1} $ and $ C_{2} $, independent of $ H_{x}$ and $H_{y} $, such that
$$
\|w - \Pi_{H}w\|_{h,\infty} \leq C_1 H^{4} \quad {\rm and} \quad \|w - \Pi_{H}w\|_{h} \leq C_2 H^{4}. 
$$
\end{lem}

Next, the bounds for $ \big\{ \phi^{x}_{i,s} (x) \big\}^{3}_{s=0} $ are given in the following lemma in order to support the proof of boundedness conclusions for the operator $\Pi_{H}$.
\begin{lem}\label{lem:basis_bounds} The cubic Lagrange interpolation basis functions $ \big\{ \phi^{x}_{i,s} (x) \big\}^{3}_{s=0} $ are bounded, i.e., for $ 1 \leq i \leq N^{H}_{x}-2 $, we have
\begin{equation*}
	\begin{array}{l}
		\left\vert \phi^{x}_{i,0} (x) \right\vert \leq 
		\left\{
		\begin{split}
			& 1, \qquad   \qquad \  \   \   x \in ( x^{H}_{i-1},x^{H}_{i} ),\\
			& \frac{\sqrt{3}}{27}, \qquad  \quad  \  \  \   x \in ( x^{H}_{i},x^{H}_{i+1} ), \\
			& \frac{\sqrt{3}}{27}, \qquad \quad \   x \in ( x^{H}_{i+1},x^{H}_{i+2} ),
		\end{split}
		\right.
	\end{array}
\qquad 
	\begin{array}{l}
	\left\vert \phi^{x}_{i,1}(x)  \right\vert \leq 
	\left\{
	\begin{split}
		& \frac{7\sqrt{7} + 10}{27}, \qquad  \quad  x \in ( x^{H}_{i-1},x^{H}_{i} ),\\
		& 1, \qquad  \qquad \qquad   \  \  \    x \in ( x^{H}_{i},x^{H}_{i+1} ), \\
		& \frac{7\sqrt{7} - 10}{27}, \qquad  \  \  x \in ( x^{H}_{i+1},x^{H}_{i+2} ),
	\end{split}
	\right.
\end{array}
\end{equation*}
\begin{equation*}
	\begin{array}{l}
		\left\vert \phi^{x}_{i,2}(x)  \right\vert \leq 
		\left\{
		\begin{split}
			& \frac{7\sqrt{7} - 10}{27}, \qquad  \quad   x \in ( x^{H}_{i-1},x^{H}_{i} ),\\
			& 1, \qquad  \qquad \qquad   \  \  \   x \in ( x^{H}_{i},x^{H}_{i+1} ), \\
			& \frac{7\sqrt{7} + 10}{27}, \qquad \  \  x \in ( x^{H}_{i+1},x^{H}_{i+2} ),
		\end{split}
		\right.
	\end{array}
	\qquad 
	\begin{array}{l}
		\left\vert \phi^{x}_{i,3}(x)  \right\vert \leq 
		\left\{
		\begin{split}
			& \frac{\sqrt{3}}{27}, \qquad  \   \quad \   x \in ( x^{H}_{i-1},x^{H}_{i} ),\\
			& \frac{\sqrt{3}}{27}, \qquad  \quad  \  \     x \in ( x^{H}_{i},x^{H}_{i+1} ), \\
			& 1, \qquad   \qquad \  x \in ( x^{H}_{i+1},x^{H}_{i+2} ).
		\end{split}
		\right.
	\end{array}
\end{equation*}
\end{lem}
{\bf Proof.}
Consider  auxiliary function 
$$ Z (x)=\left(x - x^{H}_{i}\right) \left(x - x^{H}_{i+1}\right) \left(x - x^{H}_{i+2}\right)= (x - iH_{x})  (x - (i+1)H_{x}) (x - (i+2)H_{x}).$$ 
Note that its first derivative $ Z^{\prime}(x) = 3 x^{2} - 6(i+1) H_{x}x + (3 i^{2} + 6i + 2) H_{x}^{2}$ has two zero-points $ x_{-} = \left( i + 1 - \frac{\sqrt{3}}{3} \right)H_{x} $ and $ x_{+} = \left( i + 1 + \frac{\sqrt{3}}{3} \right)H_{x} $. Thus, we have  
\begin{equation*}
	\begin{aligned}
		& \max_{ x \in ( x^{H}_{i-1},x^{H}_{i} ) } \left\vert Z(x) \right\vert = \max\left\{ Z\left( \left(i-1\right)H_{x} \right), Z\left( iH_{x} \right) \right\} = 6H_{x}^{3}, \\
		& \max_{ x \in ( x^{H}_{i},x^{H}_{i+1} ) } \left\vert Z(x) \right\vert = \max\left\{ Z\left( iH_{x} \right), Z\left( \left(i+1-\frac{\sqrt{3}}{3}\right)H_{x} \right) , Z\left( \left(i+1\right)H_{x} \right)\right\} = \frac{2\sqrt{3}}{9} H_{x}^{3}, \\
		& \max_{ x \in ( x^{H}_{i+1},x^{H}_{i+2} ) } \left\vert Z(x) \right\vert = \max\left\{ Z\left( \left(i+1\right)H_{x} \right), Z\left( \left(i+1+\frac{\sqrt{3}}{3}\right)H_{x} \right) , Z\left( \left(i+2\right)H_{x} \right)\right\} = \frac{2\sqrt{3}}{9} H_{x}^{3},
	\end{aligned}
\end{equation*}
which, together with the definition of $ \phi^{x}_{i,0}(x) $, leads to the first conclusion. The remaining conclusions can be similarly proved. 
\qed

Next, we apply Lemma \ref{lem:basis_bounds} to establish the boundedness results in the discrete $ L^{2} $ and $ L^{\infty} $ norms for the piecewise bi-cubic Lagrange interpolation operator $\Pi_{H}$. 
\begin{lem}\label{lem:Lagrange_bounds_L2} For any $w\in\mathcal{V}_H^{0}$, the following estimate holds
\begin{equation*}
	\begin{aligned}
		\| \Pi_{H} w \|_{h} \leq C_{3} \| w \|_{H},
	\end{aligned}
\end{equation*}
where $ C_{3} = 4 \left(\frac{3 + 27^{2} + ( 10 + 7\sqrt{7} )^{2}}{27^{2}}\right) $. 
\end{lem}
{\bf Proof.}
Denote $ \xi_{i,j} =   \Pi_{H,x} w_{i,j} $, and then 
\begin{equation}\label{L2_Bound_e1}
	\begin{aligned}
		\| \Pi_{H} w \|^{2}_{h} = h_{x} h_{y} \sum^{N^{h}_{x}-1}_{i=1} \sum^{N^{h}_{y}-1}_{j=1} \left( \Pi_{H} w_{i,j} \right)^{2} = h_{x} \sum^{N^{h}_{x}-1}_{i=1} \left( h_{y} \sum^{N^{h}_{y}-1}_{j=1} \left( \Pi_{H,y} \xi_{i,j} \right) ^{2} \right).
	\end{aligned}
\end{equation}
For fixed $ 1 \leq i \leq N^{h}_{x}-1 $, we know
\begin{equation}\label{L2_Bound_e2}
	\begin{aligned}
		h_{y} \sum^{N^{h}_{y}-1}_{j=1} \left( \Pi_{H,y} \xi_{i,j} \right)^{2} & = h_{y} \sum^{N^{H}_{y}}_{k=1} \sum^{kM_{y}}_{j = (k-1)M_{y} + 1} \left( \Pi_{H,y} \xi_{i,j} \right)^{2} \\
		& = h_{y} \sum^{M_{y}}_{j=1} \left( \Pi_{H,y} \xi_{i,j} \right)^{2} + h_{y}\sum^{N^{H}_{y}-1}_{k=2} \sum^{kM_{y}}_{j = (k-1)M_{y} + 1} \left( \Pi_{H,y} \xi_{i,j} \right)^{2} + h_{y} \sum^{ N^{h}_{y}-1 }_{j = (N^{H}_{y}-1) M_{y}+ 1 } \left( \Pi_{H,y} \xi_{i,j} \right)^{2} \\
		& =: I_{1} + I_{2} + I_{3}.
	\end{aligned}
\end{equation}

We use the identity $ \left(\sum^{4}_{i=1} a_{i}\right)^{2} \leq 4 \sum^{4}_{i=1} a^{2}_{i} $ and the definition \eqref{Def:cubic_L} of cubic Lagrange interpolation operator $ \Pi_{H,y} $ to obtain 
\begin{equation}\label{L2_Bound_e3}
	\begin{aligned}
		I_{1} = h_{y} \sum^{M_{y}}_{j=1} \left( \sum^{3}_{s=0} \phi^{y}_{0,s}(y^{h}_{j}) \xi_{i,s} \right)^{2} \leq 4 h_{y} \sum^{M_{y}}_{j=1} \sum^{3}_{s=0} \left(\phi^{y}_{0,s}(y^{h}_{j})\right)^{2} \left(\xi_{i,s}\right)^{2} ,
	\end{aligned}
\end{equation}
which, together with Lemma \ref{lem:basis_bounds}, yields 
\begin{equation}\label{L2_Bound_e4}
	\begin{aligned}
		I_{1} & \leq 4 h_{y} \sum^{M_{y}}_{j=1} \left( \xi_{i,0}^{2} + \frac{(7\sqrt{7} + 10)^{2}}{27^{2}} \xi_{i,1}^{2} + \frac{(7\sqrt{7} - 10)^{2}}{27^{2}} \xi_{i,2}^{2} + \frac{3}{27^{2}} \xi_{i,3}^{2} \right) \\
		& = 4 H_{y} \left( \xi_{i,0}^{2} + \frac{(7\sqrt{7} + 10)^{2}}{27^{2}} \xi_{i,1}^{2} + \frac{(7\sqrt{7} - 10)^{2}}{27^{2}} \xi_{i,2}^{2} + \frac{3}{27^{2}} \xi_{i,3}^{2} \right).
	\end{aligned}
\end{equation}
Analogous to the process \eqref{L2_Bound_e3}--\eqref{L2_Bound_e4}, we can easily derive
\begin{equation}\label{L2_Bound_e5}
	\begin{aligned}
			I_{2} \leq 4 H_{y} \sum^{N^{H}_{y}-1}_{k=2} \left( \frac{3}{27^{2}} \xi_{i,k-2}^{2} + \xi_{i,k-1}^{2} + \xi_{i,k}^{2} + \frac{3}{27^{2}} \xi_{i,k+1}^{2} \right).
		\end{aligned}
\end{equation}
\begin{equation}\label{L2_Bound_e6}
	\begin{aligned}
		I_{3} \leq 4 H_{y} \left( \frac{3}{27^{2}} \xi_{i,N^{H}_{y}-3}^{2} + \frac{(7\sqrt{7}-10)^{2}}{27^{2}} \xi_{i,N^{H}_{y}-2}^{2} + \frac{(7\sqrt{7}+10)^{2}}{27^{2}} \xi_{i,N^{H}_{y}-1}^{2} + \xi_{i,N^{H}_{y}}^{2} \right).
	\end{aligned}
\end{equation}

Inserting \eqref{L2_Bound_e4}--\eqref{L2_Bound_e6} into \eqref{L2_Bound_e2} gives us 
\begin{equation}\label{L2_Bound_e7}
	\begin{aligned}
		h_{y} \sum^{N^{h}_{y}-1}_{j=1} \left( \Pi_{H,y} \xi_{i,j} \right)^{2} & \leq 4 H_{y} \left( \frac{3 + 27^{2} + (7\sqrt{7} + 10)^{2}}{27^{2}} \xi_{i,1}^{2} + \frac{3 + 2 \times 27^{2} + (7\sqrt{7} - 10)^{2}}{27^{2}} \xi_{i,2}^{2} \right.		\\
		& \qquad \quad 
		+ \frac{9 + 2 \times 27^{2} }{27^{2}} \xi_{i,3}^{2} + \frac{6 + 2 \times 27^{2} }{27^{2}} \sum^{N^{H}_{y}-4}_{j=4} \xi_{i,j}^{2} + \frac{9 + 2 \times 27^{2} }{27^{2}} \xi_{i,N^{H}_{y}-3}^{2}\\
		& \qquad \quad \left.
		+ \frac{3 + 2 \times 27^{2} + (7\sqrt{7} - 10)^{2}}{27^{2}} \xi_{i,N^{H}_{y}-2}^{2} + \frac{3 + 27^{2} + (7\sqrt{7} + 10)^{2}}{27^{2}} \xi_{i,N^{H}_{y}-1}^{2}
		\right) \\
		& \leq C_{3} H_{y} \sum^{N^{H}_{y}-1}_{j=1} \xi_{i,j}^{2},
	\end{aligned}
\end{equation}
which further implies
\begin{equation}\label{L2_Bound_e8}
	\begin{aligned}
		\| \Pi_{H} w \|^{2}_{h} = h_{x} \sum^{N^{h}_{x}-1}_{i=1} \left( h_{y} \sum^{N^{h}_{y}-1}_{j=1} \left( \Pi_{H,y} \xi_{i,j} \right)^{2} \right) \leq C_{3} H_{y} \sum^{N^{H}_{y}-1}_{j=1} \left( h_{x} \sum^{N^{h}_{x}-1}_{i=1} \left( \Pi_{H,x} w_{i,j}\right)^{2} \right).
	\end{aligned}
\end{equation}

For fixed $ 1 \leq j \leq N^{H}_{y} -1 $, we replace $ \left\{ \xi, y, j \right\} $ with $ \left\{ w, x, i \right\} $ in \eqref{L2_Bound_e2}--\eqref{L2_Bound_e8} to similarly obtain 
\begin{equation}\label{L2_Bound_e9}
	\begin{aligned}
		h_{x} \sum^{N^{h}_{x}-1}_{i=1} \left( \Pi_{H,x} w_{i,j} \right)^{2} & \leq  C_{3} H_{x} \sum^{N^{H}_{x}-1}_{i=1} w_{i,j}^{2}.
	\end{aligned}
\end{equation}
Consequently, we have
\begin{equation}\label{L2_Bound_e10}
	\begin{aligned}
		\| \Pi_{H} w \|^{2}_{h}  \leq C_{3}^{2} H_{x} H_{y} \sum^{N^{H}_{x}-1}_{i=1} \sum^{N^{H}_{y}-1}_{j=1} w_{i,j}^{2} = C_{3}^{2} \| w \|^{2}_{H}.
	\end{aligned}
\end{equation}
The proof is completed.
\qed

\begin{lem}\label{lem:Lagrange_bounds_Inf} For any $w\in\mathcal{V}_H^{0}$, the following estimate holds
	\begin{equation*}
		\begin{aligned}
			\| \Pi_{H} w \|_{h,\infty} \leq C_{4} \| w \|_{H,\infty},
		\end{aligned}
	\end{equation*}
	where $ C_{4} = \left( \frac{ 54 + 2\sqrt{3} }{27} \right)^{2} $. 
\end{lem}
{\bf Proof.}
Suppose $ \| \Pi_{H} w \|_{h,\infty} = \left\vert \Pi_{H} w_{i^{*},j^{*}} \right\vert $ where $ i M_{x} \leq i^{*} \leq (i+1) M_{x} $ and $ j M_{y} \leq j^{*} \leq (j+1) M_{y} $. 
Denote $ \xi_{i^*,j} = \Pi_{H,x} w_{i^*,j} $ for $ 0 \leq j \leq N^{H}_{y}$, then the triangle inequality and Lemma \ref{lem:basis_bounds} give us 
\begin{equation}\label{Inf_eq1}
	\begin{aligned}
		\| \Pi_{H} w \|_{h,\infty} & = \left\vert \Pi_{H,y} \xi_{i^{*}, j^{*} } \right\vert = \left\vert \sum^{3}_{s=0} \phi^{y}_{j,s}(y^{h}_{j^{*}}) \xi_{i^{*},j-1+s} \right\vert \\
		& \leq \frac{ \sqrt{3} }{ 27 } \left\vert \xi_{i^{*},j-1} \right\vert + \left\vert \xi_{i^{*},j} \right\vert + \left\vert \xi_{i^{*},j+1} \right\vert + \frac{ \sqrt{3} }{ 27 } \left\vert \xi_{i^{*},j+2} \right\vert \\
		& \leq \frac{ 54 + 2\sqrt{3} }{27} \max_{ 1 \leq k \leq N^{H}_{y} - 1 } \left\vert \xi_{i^{*},k} \right\vert, \qquad \text{ if } \  M_{y}  \leq j^{*} \leq (N^{H}_{y}-1) M_{y},
	\end{aligned}
\end{equation}
and
\begin{equation}\label{Inf_eq3}
	\begin{aligned}
		\| \Pi_{H} w \|_{h,\infty} & \leq \frac{ 7 \sqrt{7} + 10 }{27} \left\vert \xi_{i^{*},1} \right\vert + \frac{ 7 \sqrt{7} - 10 }{27} \left\vert \xi_{i^{*},2} \right\vert + \frac{ \sqrt{3} }{27} \left\vert \xi_{i^{*},3} \right\vert \\
		& \leq \frac{ 14\sqrt{7} + \sqrt{3} }{27} \max_{ 1 \leq k \leq N^{H}_{y} - 1 } \left\vert \xi_{i^{*},k} \right\vert, \qquad \text{ if } \  1  \leq j^{*} < M_{y},
	\end{aligned}
\end{equation}
\begin{equation}\label{Inf_eq4}
	\begin{aligned}
		\| \Pi_{H} w \|_{h,\infty} \leq \frac{ 14\sqrt{7} + \sqrt{3} }{27} \max_{ 1 \leq k \leq N^{H}_{y} - 1 } \left\vert \xi_{i^{*},k} \right\vert, \qquad \text{ if } \  (N^{H}_{y} - 1)M_{y} < j^{*} \leq N^{h}_{y} - 1.
	\end{aligned}
\end{equation}

Next, it is necessary to estimate  $\max_{ 1 \leq k \leq N^{H}_{y} - 1 } \left\vert \xi_{i^{*},k} \right\vert$.  For fixed index $ 1 \leq k \leq N^{H}_{y}-1 $, 
let $\left\vert  \xi_{i^{**},k} \right\vert= \max_{ 1 \leq \ell \leq N^{h}_{x}-1 } \left\vert \xi_{\ell,k} \right\vert$ where $ i M_{x} \leq i^{**} \leq (i+1) M_{x} $. Thus, with a similar treatment to the above estimates, we derive
\begin{equation}\label{Inf_eq5}
	\begin{aligned}
		\left\vert \xi_{i^{*},k} \right\vert =\left\vert \Pi_{H,x} w_{i^{*},k} \right\vert \leq \left\vert \Pi_{H,x} w_{i^{**},k} \right\vert \leq \frac{ 54 + 2\sqrt{3} }{27} \max_{ 1 \leq \ell \leq N^{H}_{x} - 1 } \left\vert w_{\ell,k} \right\vert,
	\end{aligned}
\end{equation}
which, together with \eqref{Inf_eq1}--\eqref{Inf_eq4}, completes the proof. \qed

\section{A nonlinear variable-step compact difference scheme }\label{sec:Non}
In this section, combined with the variable-step BDF2 method, we are committed to establishing a nonlinear compact difference scheme for the semilinear parabolic equation \eqref{Model:Bihar}--\eqref{Bound:Bihar} enclosed with  Dirichlet boundary condition. Meanwhile, we shall develop the corresponding error estimates by imposing the following regularity assumption
\begin{equation}\label{Regu:Bihar}
	\begin{aligned}
		u \in C(0,T;H^{6}(\Omega)) \cap C^{3}(0,T;L^2(\Omega)).
	\end{aligned}
\end{equation}
By the Sobolev embedding theorem, the above assumption implies that the solution of problem \eqref{Model:Bihar}--\eqref{Bound:Bihar} is bounded, i.e., there exist two constants $ m $ and $ M $ such that $ u \in B := [m,M] $. Furthermore, for a fixed small $ \delta > 0 $, we denote $ B_{\delta} := [m-\delta, M+ \delta] $. In this paper, we only assume that the nonlinear term $ f \in C^{2} ( B_{\delta}) $ holds  locally  on $ B_{\delta} $, that is, 
\begin{equation}\label{Regu:Bihar:e1}
\vert  f(w) \vert+ \vert f'(w) \vert + \vert f''(w) \vert \leq K_{f} \qquad \text{for } \  w \in B_{\delta}.
\end{equation}

\subsection{Nonlinear compact difference scheme}
To construct variable temporal stepsize schemes, we consider a nonuniform temporal mesh partition $ 0 = t_{0} < t_{1} < \cdots < t_{N} = T$ with  temporal stepsize $ \tau_{k} := t_{k} - t_{k-1} $ for $ 1 \leq k \leq N $. Denote the maximum temporal stepsize $ \tau := \max_{1 \leq k \leq N} \tau_{k}$ and adjacent temporal stepsize ratio $ r_{k} := \tau_{k}/\tau_{k-1}$ $ (k \geq 2) $. Set $ w^{k} = w(t_{k}) $ and $ \nabla_{\tau} w^{k} = w^{k} - w^{k-1}  $, then the variable-step BDF2 formula is defined by 
$$
\mathcal{D}_{2} w^{n} := \frac{ 1 + 2 r_{n} }{\tau_{n} (1 + r_{n}) } \nabla_{\tau} w^{n} - \frac{ r_{n}^{2} }{\tau_{n} (1 + r_{n}) } \nabla_{\tau} w^{n-1}\quad \text{for } \   n \ge 2.
$$
In particular, when $n=1$,  we use the BDF1 (i.e., backward Euler) formula $\mathcal{D}_{1} w^{1} = \frac{1}{\tau_{1}} \nabla_{\tau} w^{1}$ for the first time level discretization.
Let $ r_{1} = 0 $, we then rewrite the above variable-step BDF formula as a unified discrete convolution summation
\begin{equation}\label{Def:BDF2}
	\begin{aligned}
		\mathcal{D}_{2} w^{n} = \sum^{n}_{k=1} b^{(n)}_{n-k} \nabla_{\tau} w^{k}\quad \text{for } \  n \geq 1,
	\end{aligned}
\end{equation}
where the discrete convolution kernels $ b^{(n)}_{n-k} $ are defined by $ b^{(1)}_{0} := 1/\tau_{1} $ and 
\begin{equation}\label{Def:DCK}
	\begin{aligned}
		b^{(n)}_{0} := \frac{ 1 + 2 r_{n} }{\tau_{n} (1 + r_{n}) }, \quad b^{(n)}_{1} := - \frac{ r_{n}^{2} }{\tau_{n} (1 + r_{n}) }, \quad \text{and} \quad b^{(n)}_{j} := 0 \   \   \text{for} \  2 \leq j \leq n.
	\end{aligned}
\end{equation}

Next, we introduce the discrete orthogonal convolution (DOC) kernels of $\left\lbrace b^{(n)}_{n-k} \right\rbrace $ by \cite{2021_MOC_Liao,JSC_Liao_2022,JSC_Zhang_2022} 
\begin{equation}\label{Def:DOC}
	\begin{aligned}
		\sum^{n}_{m = k} \theta^{(n)}_{n-m} b^{(m)}_{m-k} = \delta_{nk}, \quad  1 \leq k \leq n, \  1 \leq n \leq N,
	\end{aligned}
\end{equation}
where $ \delta_{nk} = 1$ if $ n = k $ and $ \delta_{nk} = 0$ if $ n \neq k $. By exchanging the summation order and using definition \eqref{Def:DOC}, it is easy to check that
\begin{equation}\label{DOC_P}
	\begin{aligned}
	    \sum_{m=1}^n \theta_{n-m}^{(n)} \mathcal{D}_2 w^m = \sum_{k=1}^n \nabla_\tau w^k \sum_{m=k}^n \theta_{n-m}^{(n)} b_{m-k}^{(m)} = \nabla_{\tau} w^n, \quad  1 \leq n \leq N.
    \end{aligned}
\end{equation}

Let $U^{n}_{i,j}:=u(x^{h}_{i}, y^{h}_{i}, t_{n})$ be the exact nodal solutions. Then, we apply the compact difference operator $\mathcal{A}_{h} $ in space and variable-step BDF method in time to get the following equation
\begin{equation}\label{Non_eq_1}
	\begin{aligned}
		\mathcal{D}_2 \mathcal{A}_{h} U^{n}_{i,j}  - c \Lambda_{h} U^{n}_{i,j} = \mathcal{A}_{h} f ( U^{n}_{i,j} ) + \mathcal{A}_{h} g^{n}_{i,j} + R^{n}_{i,j}, \quad  1 \leq n \leq N,
	\end{aligned}
\end{equation}
for $ 1 \leq i \leq N^{h}_{x}-1 $, $ 1 \leq j \leq N^{h}_{y}-1 $. Here $R^{n}_{i,j} := (R_{t}^{n})_{i,j} + (R_{s}^{n})_{i,j} $ is  the local truncation error at point $(x^{h}_{i}, y^{h}_{i}, t_{n})$, where
$$
(R_{t}^{n})_{i,j} = \mathcal{A}_{h} \left( u_{t}(x^{h}_{i}, y^{h}_{j}, t_{n}) - \mathcal{D}_2 U^{n}_{i,j} \right) \quad \text{ and } \quad (R_{s}^{n})_{i,j} = -c \mathcal{A}_{h} \Delta u(x^{h}_{i}, y^{h}_{i}, t_{n}) + c \Lambda_{h} U^{n}_{i,j}. 
$$
Under the regularity condition \eqref{Regu:Bihar}, by using the Taylor expansion and well-known Bramble-Hilbert Lemma, it is easy to see
\begin{equation}\label{Trun_err_t}
	\begin{aligned}
\| R_{t}^{1} \|_{h} \leq C_{5} \tau_{1}, \quad \| R_{t}^{n} \|_{h} \leq C_{6} \tau_{n} \tau \quad \text{for} \  2 \leq n \leq N,
	\end{aligned}
\end{equation}
\begin{equation}\label{Trun_err_s}
	\begin{aligned}
\| R_{s}^{n} \|_{h} \leq C_{7} h^{4} \quad \text{for} \  1\leq n \leq N.
	\end{aligned}
\end{equation}
Let $u^{n}_{i,j}$ be the finite difference approximations to $U^{n}_{i,j}$, then we drop the local truncation errors in \eqref{Non_eq_1} to obtain the nonlinear compact difference scheme
\begin{equation}\label{sch:Non_CDM}
	\begin{aligned}
		\mathcal{D}_2 \mathcal{A}_{h} u^{n}_{i,j}  - c \Lambda_{h} u^{n}_{i,j} = \mathcal{A}_{h} f ( u^{n}_{i,j} ) + \mathcal{A}_{h} g^{n}_{i,j}, \quad (i,j) \in  \omega_{h}, ~1 \leq n \leq N.
	\end{aligned}
\end{equation}

In the following, we aim to prove the unique solvability and error analysis for the nonlinear variable temporal stepsize compact difference scheme \eqref{sch:Non_CDM}. For this purpose, we introduce an auxiliary solution $ \bar{u}^{n} =\{\bar{u}^{n}_{i,j}\}$ satisfying the following scheme
\begin{equation}\label{sch:Non_AE}
	\begin{aligned}
		\mathcal{D}_2 \mathcal{A}_{h} \bar{u}^{n}_{i,j}  - c \Lambda_{h} \bar{u}^{n}_{i,j} = \mathcal{A}_{h} \bar{f} ( \bar{u}^{n}_{i,j} ) + \mathcal{A}_{h} g^{n}_{i,j}, \quad (i,j) \in  \omega_{h}, ~ 1 \leq n \leq N,
	\end{aligned}
\end{equation}
 subject to the same initial and boundary conditions \eqref{Initial:Bihar}--\eqref{Bound:Bihar}. Here $ \bar{f} (w) $ represents a cut-off function of $ f $ on $ B_{\delta} $ such that
\begin{equation}\label{f_cut}
	\begin{array}{l}
		\bar{f}(w) :=
		\left\{
		\begin{split}
			& f(m - \delta), \qquad \  w < m-\delta,\\
			& f(w), \qquad \qquad  w \in B_{\delta},\\
			& f(M+\delta), \qquad \  w > M+\delta.
		\end{split}
		\right.
	\end{array}
\end{equation}
It is obvious that $ \bar{f}(w) $ is globally Lipschitz continuous on $ \mathbb{R} $ with  Lipschitz constant $ K_{f} $.

Below we shall show the unique solvability and error analysis for the nonlinear auxiliary scheme \eqref{sch:Non_AE} instead of  \eqref{sch:Non_CDM}.
Several useful lemmas are presented below which will be used for our theoretical analysis. 
\begin{lem}[\cite{JSC_Zhang_2022,JSC_Liao_2022}]\label{lem:DOC_PD} Assume that the adjacent temporal stepsize ratios $ r_{k} $ satisfy $ 0 < r_{k} < 4.8645 $. The DOC kernels $ \{\theta^{(n)}_{n-k} \}$ defined in \eqref{Def:DOC} are positive semi-definite, i.e., for any real sequence $ \{ w_{k} \}^{n}_{k=1} $, it holds that 
	\begin{equation*}
		\begin{aligned}
			\sum_{k=1}^n w_k \sum_{m=1}^k \theta_{k-m}^{(k)} w_m \geq 0 \quad {\rm for} \   n \geq 1.
		\end{aligned}
	\end{equation*}
\end{lem}

\begin{lem}[\cite{2021_MOC_Liao,JSC_Zhang_2022}]\label{lem:DOC_estimate} The DOC kernels $ \{\theta^{(n)}_{n-k} \}$  defined in \eqref{Def:DOC} have the following properties
	\begin{equation*}
		\begin{aligned}
			\theta_{n-m}^{(n)}>0 \quad   {\rm for} \  1 \leq m \leq n, \quad \sum_{m=1}^n \theta_{n-m}^{(n)}=\tau_n \quad  {\rm for} \  n \geq 1, \\
			\sum^{n}_{m=1} \sum^{n}_{k = m} \theta^{(k)}_{k-m} = t_{n}, \quad {\rm and} \quad \sum^{n}_{k = m} \theta^{(k)}_{k-m} \leq 2 \tau \quad {\rm for} \  1 \leq  m \leq n.
		\end{aligned}
	\end{equation*}
\end{lem}

\begin{lem}[\cite{2021_MOC_Liao}]\label{lem:Gronwall} Let $ \lambda \geq 0 $, and the time sequences $ \{ \xi^{k} \}_{k=0}^{N} $ and $ \{ w^{k} \}_{k=1}^{N} $ be nonnegative. If
	\begin{equation}\label{Gron:1}
		\begin{aligned}
			w^{n} \leq \lambda \sum^{n-1}_{k=1} \tau_{k} w^{k} + \sum^{n}_{k=0} \xi^{k} \quad {\rm for} \  1 \leq n \leq N,
		\end{aligned}
	\end{equation}
	then it holds that
	\begin{equation}\label{Gron:2}
		\begin{aligned}
			w^{n} \leq \exp(\lambda t_{n-1}) \sum^{n}_{k=0} \xi^{k} \quad {\rm for} \  1 \leq n \leq N.
		\end{aligned}
	\end{equation}
\end{lem}

\subsection{ Unique solvability }
It is noticed that in \cite{SINUM_Liao_2020}, the authors proposed a nonlinear variable-step BDF2 scheme  for the Allen--Cahn equation with a polynomial type double-well nonlinear potential, and they proved its unique solvability by showing that the solution of the nonuniform
BDF2 scheme is equivalent to a minimization problem with strictly convex energy functional. In this subsection, the unique solvability of the auxiliary BDF2 scheme  \eqref{sch:Non_AE} for  semilinear parabolic equations with general nonlinearity will be discussed by the Browder's fixed point theorem (see e.g. \cite{Browder_2_IMA}). By homogenization treatment, it suffices to  consider the corresponding homogeneous case. 
\begin{thm}\label{thm:Non_solva} 
	The auxiliary nonlinear compact difference scheme \eqref{sch:Non_AE} is solvable if the maximum temporal stepsize satisfies $ \tau \leq \frac{1}{ 4 K_{f} } $.
\end{thm}
{\bf Proof.}
Note that \eqref{sch:Non_AE} can be equivalently rewritten as
\begin{equation}\label{Solv_e1}
	\begin{aligned}
		b^{(n)}_{0} \mathcal{A}_{h} \bar{u}^{n}_{i,j} - G^{n}_{i,j} - c \Lambda_{h} \bar{u}^{n}_{i,j} - \mathcal{A}_{h} \bar{f} ( \bar{u}^{n}_{i,j} ) - \mathcal{A}_{h} g^{n}_{i,j} = 0, \quad (i,j) \in  \omega_{h}, 
	\end{aligned}
\end{equation}
where $ G^{n}_{i,j} $ is defined as
\begin{equation*}
	\begin{aligned}
		G^{n}_{i,j} := ( b^{(n)}_{0} - b^{(n)}_{1} ) \mathcal{A}_{h} \bar{u}^{n-1}_{i,j} + b^{(n)}_{1} \mathcal{A}_{h} \bar{u}^{n-2}_{i,j}.
	\end{aligned}
\end{equation*}
Denote the mapping
$$
[ T \bar{u} ]^{n}_{i,j} := b^{(n)}_{0} \mathcal{A}_{h} \bar{u}^{n}_{i,j} - G^{n}_{i,j} - c \Lambda_{h} \bar{u}^{n}_{i,j} - \mathcal{A}_{h} \bar{f} ( \bar{u}^{n}_{i,j} ) - \mathcal{A}_{h} g^{n}_{i,j}, \quad (i,j) \in  \omega_{h}. 
$$
Then, taking the inner product of $ T \bar{u}^{n} $ with $ \bar{u}^{n} $ in the sense of $(\cdot,\cdot)_h$ gives us
\begin{equation}\label{Solv_e2}
	\begin{aligned}
		\left( T \bar{u}^{n}, \bar{u}^{n} \right)_{h} & = b^{(n)}_{0} \left( \mathcal{A}_{h} \bar{u}^{n}, \bar{u}^{n} \right)_{h} - \left( G^{n}, \bar{u}^{n} \right)_{h} - c \left( \Lambda_{h} \bar{u}^{n}, \bar{u}^{n} \right)_{h} \\
		& \quad - \left( \mathcal{A}_{h} \bar{f} ( \bar{u}^{n} ) , \bar{u}^{n} \right)_{h} - \left( \mathcal{A}_{h} g^{n} , \bar{u}^{n} \right)_{h} := \sum^{5}_{k=1} S_{k}.
	\end{aligned}
\end{equation}

Next, we estimate \eqref{Solv_e2} term-by-term. For $ S_{2} $, utilizing Cauchy-Schwarz inequality yields 
\begin{equation*}
	\begin{aligned}
		S_{2} \geq - \left( ( b^{(n)}_{0} - b^{(n)}_{1} ) \| \bar{u}^{n-1} \|_{\mathcal{A}, h} - b^{(n)}_{1}\| \bar{u}^{n-2} \|_{\mathcal{A}, h} \right) \| \bar{u}^{n} \|_{\mathcal{A}, h}.
	\end{aligned}
\end{equation*}
Moreover, by summation by parts and homogeneous boundary conditions, and noting that $ \mathcal{A}_{h,x} $ and $ \mathcal{A}_{h,y} $ are self-adjoint and positive definite operators,  
we have
\begin{equation*}
	\begin{aligned}
		S_{3} = \left( \mathcal{A}_{h,y} d_{h,x} \bar{u}^{n}, d_{h,x} \bar{u}^{n} \right)_{h,x} + \left(\mathcal{A}_{h,x}  d_{h,y} \bar{u}^{n}, d_{h,y} \bar{u}^{n} \right)_{h,y} \geq 0.
	\end{aligned}
\end{equation*}
Besides, we use Cauchy-Schwarz inequality to get
\begin{equation*}
	\begin{aligned}
		S_{4} + S_{5} = - \left( \mathcal{A}_{h} \left( \bar{f}( \bar{u}^{n} ) + g^{n} \right), \bar{u}^{n} \right)_{h} \geq - \interleave \bar{f}( \bar{u}^{n} ) + g^{n} \interleave_{h} \| \bar{u}^{n} \|_{h},
	\end{aligned}
\end{equation*}
which further, together with the globally Lipschitz continuous of $ \bar{f} $ and Lemma \ref{lem:cpmpact1}, leads to
\begin{equation*}
	\begin{aligned}
		S_{4} + S_{5} & \geq - K_{f} \left( \| \bar{u}^{n} \|_{h} + \| u^{0} \|_{h} \right) \| \bar{u}^{n} \|_{h} - \interleave \bar{f}( u^{0} ) + g^{n} \interleave_{h} \| \bar{u}^{n} \|_{h} \\
		& \geq - 3 K_{f} \| \bar{u}^{n} \|^{2}_{\mathcal{A},h} - \sqrt{3} \left( K_{f} \| u^{0} \|_{h} + \interleave \bar{f}( u^{0} ) + g^{n} \interleave_{h} \right) \| \bar{u}^{n} \|_{\mathcal{A},h}.
	\end{aligned}
\end{equation*}

Denote 
\begin{equation*}
	\begin{aligned}
		Q := ( b^{(n)}_{0} - b^{(n)}_{1} ) \| \bar{u}^{n-1} \|_{\mathcal{A}, h} - b^{(n)}_{1}\| \bar{u}^{n-2} \|_{\mathcal{A}, h} + \sqrt{3} \left( K_{f} \| u^{0} \|_{h} + \interleave \bar{f}( u^{0} ) + g^{n} \interleave_{h} \right).
	\end{aligned}
\end{equation*}
Then inserting the above estimates into \eqref{Solv_e2} gives
\begin{equation*}
	\begin{aligned}
		 \left( T\bar{u}^{n}, \bar{u}^{n} \right)_{h} \geq \left( b^{(n)}_{0} - 3 K_{f} \right) \| \bar{u}^{n} \|^{2}_{\mathcal{A},h} - Q \| \bar{u}^{n} \|_{\mathcal{A},h}.
	\end{aligned}
\end{equation*}
Thus, if $ \tau_{n} \leq \frac{1}{ 4 K_{f} } $, we conclude that
\begin{equation*}
	\begin{aligned}
		\left( T\bar{u}^{n}, \bar{u}^{n} \right)_{h} \geq  \left( K_{f} \| \bar{u}^{n} \|_{\mathcal{A},h}  - Q \right) \| \bar{u}^{n} \|_{\mathcal{A},h} \geq 0,
	\end{aligned}
\end{equation*}
when $ \| \bar{u}^{n} \|_{\mathcal{A},h} = Q/K_{f} $. Consequently, Browder's fixed point theorem shows there exists a $ \bar{u}^{n} \in \mathcal{V}^{0}_{h} $ such that $ T \bar{u}^{n} = 0$ and $ \| \bar{u}^{n} \|_{h} \leq \sqrt{3} \| \bar{u}^{n} \|_{\mathcal{A},h} \leq \sqrt{3} Q/K_{f} $, which implies the solvability of scheme \eqref{sch:Non_AE}.
\qed

\begin{thm}\label{thm:Non_unique} 
	The solution of the auxiliary nonlinear compact difference scheme \eqref{sch:Non_AE} is unique if the adjacent temporal stepsize ratios $ r_{k} $ satisfy $ 0 < r_{k} < 4.8645 $ and the maximum temporal stepsize $ \tau \leq \frac{1}{ 4 \sqrt{3} K_{f} } $.
\end{thm}
{\bf Proof.}
The argument is by contradiction. Suppose that there are two solutions $ \bar{u}_{1}^{n} $ and $ \bar{u}_{2}^{n} $ satisfying \eqref{sch:Non_AE} and the same initial and boundary conditions \eqref{Initial:Bihar}--\eqref{Bound:Bihar}, i.e.,
\begin{equation*}
	\begin{aligned}
		\mathcal{D}_2 [\mathcal{A}_{h} \bar{u}_{1}]^{n}_{i,j}  - c [\Lambda_{h} \bar{u}_{1}]^{n}_{i,j} = \mathcal{A}_{h} \bar{f} ( [\bar{u}_{1}]^{n}_{i,j} ) + \mathcal{A}_{h} g^{n}_{i,j},
\\
		\mathcal{D}_2 [\mathcal{A}_{h} \bar{u}_{2}]^{n}_{i,j}  - c [\Lambda_{h} \bar{u}_{2}]^{n}_{i,j} = \mathcal{A}_{h} \bar{f} ( [\bar{u}_{2}]^{n}_{i,j} ) + \mathcal{A}_{h} g^{n}_{i,j},
	\end{aligned}
\end{equation*}
for $ 1 \leq i \leq N^{h}_{x} - 1$, $ 1 \leq  j \leq N^{h}_{y}-1 $ and $ 1 \leq n \leq N $. Denote $ v^{n} := \bar{u}_{1}^{n} - \bar{u}_{2}^{n}$ and subtract these two  equations to get
\begin{equation}\label{Unique_e1}
	\begin{aligned}
		\mathcal{D}_2 \mathcal{A}_{h} v^{n}_{i,j}  - c \Lambda_{h} v^{n}_{i,j} = \mathcal{A}_{h} \bar{f} ( [\bar{u}_{1}]^{n}_{i,j} ) - \mathcal{A}_{h} \bar{f} ( [\bar{u}_{2}]^{n}_{i,j} ). 
	\end{aligned}
\end{equation}

In \eqref{Unique_e1}, we let $ n = m $, and then, multiplying it by the DOC kernels $ \theta^{(k)}_{k-m} $ and summing $ m $ from $ 1 $ to $ k $ yields 
\begin{equation}\label{Unique_e2}
	\begin{aligned}
		\nabla_{\tau} \mathcal{A}_{h} v^{k}_{i,j} - c \sum^{k}_{m=1} \theta^{(k)}_{k-m} \Lambda_{h} v^{m}_{i,j} = \sum^{k}_{m=1} \theta^{(k)}_{k-m} \mathcal{A}_{h} \left( \bar{f} ( [\bar{u}_{1}]^{m}_{i,j} ) - \bar{f} ( [\bar{u}_{2}]^{m}_{i,j} ) \right),
	\end{aligned}
\end{equation}
where we have used the orthogonal identity \eqref{Def:DOC}.

Furthermore, taking the inner product of the above equation with $ 2 v^{k} $, and summing the resulting
equation from $ k = 1 $ to $ n $ gives  
\begin{equation}\label{Unique_e3}
	\begin{aligned}
		\| v^{n} \|_{\mathcal{A},h}^{2} - 2 c \sum^{n}_{k=1} \sum^{k}_{m=1} \theta^{(k)}_{k-m} ( \Lambda_{h} v^{m}, v^{k} )_{h} & \leq 2 \sum^{n}_{k=1} \sum^{k}_{m=1} \theta^{(k)}_{k-m} \left(  \mathcal{A}_{h} \bar{f} ( \bar{u}_{1}^{n} )  - \mathcal{A}_{h} \bar{f} ( \bar{u}_{2}^{n} ) , v^{k} \right)_{h},
	\end{aligned}
\end{equation}
where $v^{0}=0 $ has been used.

Due to $ \mathcal{A}_{h,x} $ and $ \mathcal{A}_{h,y} $ are self-adjoint and positive definite operators, there exist $ \eta_{x} $ and $ \eta_{y} $ such that $ \mathcal{A}_{h,x} = \eta^{2}_{x} $ and $ \mathcal{A}_{h,y} = \eta^{2}_{y} $. 
Thus, by summation by parts and homogeneous boundary conditions, the second left-hand side term of \eqref{Unique_e3} could be rewritten as
\begin{equation}\label{Unique_e4}
	\begin{aligned}
	 & \quad	- 2 c \sum^{n}_{k=1} \sum^{k}_{m=1} \theta^{(k)}_{k-m} ( \Lambda_{h} v^{m}, v^{k} )_{h}\\
	 & = 2 c \sum^{n}_{k=1} \sum^{k}_{m=1} \theta^{(k)}_{k-m} \left( ( \mathcal{A}_{h,x} d_{h,y} v^{m}, d_{h,y} v^{k} )_{h} + ( \mathcal{A}_{h,y} d_{h,x} v^{m}, d_{h,x} v^{k} )_{h} \right) \\
	& = 2 c \sum^{n}_{k=1} \sum^{k}_{m=1} \theta^{(k)}_{k-m} \left( ( \eta_{x} d_{h,y} v^{m}, \eta_{x}d_{h,y} v^{k} )_{h} + ( \eta_{y} d_{h,x} v^{m}, \eta_{y} d_{h,x} v^{k} )_{h} \right) \geq 0,
	\end{aligned}
\end{equation}
where the positive semi-definiteness of the DOC kernels $ \theta^{(k)}_{k-m} $ (cf. Lemma \ref{lem:DOC_PD}) has been used in the last inequality.  

Using Cauchy-Schwarz inequality and Lemma \ref{lem:cpmpact1} we see
\begin{equation*}
	\begin{aligned}
		2 \sum^{n}_{k=1} \sum^{k}_{m=1} \theta^{(k)}_{k-m} \left(  \mathcal{A}_{h} \bar{f} ( \bar{u}_{1}^{n} )  - \mathcal{A}_{h} \bar{f} ( \bar{u}_{2}^{n} ), v^{k} \right)_{h} & \leq 2 \sum^{n}_{k=1} \sum^{k}_{m=1} \theta^{(k)}_{k-m} \| \bar{f} ( \bar{u}_{1}^{n} ) - \bar{f} ( \bar{u}_{2}^{n} ) \|_{\mathcal{A},h} \| v^{k} \|_{\mathcal{A},h} \\
		& \leq 2 \sum^{n}_{k=1} \sum^{k}_{m=1} \theta^{(k)}_{k-m} \| \bar{f} ( \bar{u}_{1}^{n} ) - \bar{f} ( \bar{u}_{2}^{n} ) \|_{h} \| v^{k} \|_{\mathcal{A},h},
	\end{aligned}
\end{equation*}
which, together with the global Lipschitz continuous property of $\bar{f}$, gives
\begin{equation}\label{Unique_e5}
	\begin{aligned}
		2 \sum^{n}_{k=1} \sum^{k}_{m=1} \theta^{(k)}_{k-m} \left(  \mathcal{A}_{h} \bar{f} ( \bar{u}_{1}^{n} )  - \mathcal{A}_{h} \bar{f} ( \bar{u}_{2}^{n} ), v^{k} \right) & \leq 2 K_{f} \sum^{n}_{k=1} \sum^{k}_{m=1} \theta^{(k)}_{k-m} \| v^{m} \|_{h} \| v^{k} \|_{\mathcal{A},h} \\
		& \leq 2 \sqrt{3} K_{f} \sum^{n}_{k=1} \| v^{k} \|_{\mathcal{A},h} \sum^{k}_{m=1} \theta^{(k)}_{k-m} \| v^{m} \|_{\mathcal{A}, h}.
	\end{aligned}
\end{equation}

Now, inserting \eqref{Unique_e4}--\eqref{Unique_e5} into \eqref{Unique_e3}, we obtain 
\begin{equation*}
	\begin{aligned}
		\| v^{n} \|_{\mathcal{A},h}^{2} \leq 2 \sqrt{3} K_{f} \sum^{n}_{k=1} \| v^{k} \|_{\mathcal{A},h} \sum^{k}_{m=1} \theta^{(k)}_{k-m} \| v^{m} \|_{\mathcal{A}, h}. 
	\end{aligned}
\end{equation*}
Choosing $ n^{*}$ $ ( 1 \leq n^{*} \leq n )$ such that $ \| v^{n^{*}} \|_{\mathcal{A},h} = \max_{ 1 \leq k \leq n } \| v^{n} \|_{\mathcal{A},h}$. Then, the above inequality yields 
\begin{equation*}
	\begin{aligned}
		\| v^{n^{*}} \|_{\mathcal{A},h}^{2} \leq 2 \sqrt{3} K_{f} \sum^{n^{*}}_{k=1} \| v^{k} \|_{\mathcal{A},h} \sum^{k}_{m=1} \theta^{(k)}_{k-m} \| v^{n^{*}} \|_{\mathcal{A}, h}, 
	\end{aligned}
\end{equation*}
which, eliminating  $ \| v^{n^{*}} \|_{\mathcal{A},h} $ from both sides and using Lemma \ref{lem:DOC_estimate}, implies
\begin{equation}\label{Unique_e6}
	\begin{aligned}
		\| v^{n} \|_{\mathcal{A},h} \leq \| v^{n^{*}} \|_{\mathcal{A},h} \leq 2 \sqrt{3} K_{f} \sum^{n^{*}}_{k=1} \| v^{k} \|_{\mathcal{A},h} \sum^{k}_{m=1} \theta^{(k)}_{k-m} \leq 2 \sqrt{3} K_{f} \sum^{n}_{k=1} \tau_{k} \| v^{k} \|_{\mathcal{A},h} .
	\end{aligned}
\end{equation}
Then, for $ \tau_{n} \leq 1/(4 \sqrt{3} K_{f}) $,  an application of the discrete Gr${\rm \ddot{o}}$nwall inequality in Lemma \ref{lem:Gronwall} gives  
\begin{equation}\label{Unique_e7}
	\begin{aligned}
		\| v^{n} \|_{\mathcal{A},h} \leq 4 \sqrt{3} K_{f} \sum^{n-1}_{k=1} \tau_{k} \| v^{k} \|_{\mathcal{A},h} \leq 0,
	\end{aligned}
\end{equation}
which proves the uniqueness of the solutions to scheme \eqref{sch:Non_AE}.
\qed

\subsection{ Error analysis }
Denote $ \bar{e}^{n}:= U^{n} - \bar{u}^{n} $ and $ e^{n} := U^{n}- u^{n} $ for  $ 0 \leq n \leq N $. It is easy to see that $ \bar{e}^{0} = e^{0}=0 $. Next, we shall first give an error estimates in the discrete $L^{2}$ norm between the exact solution and the numerical solution yielded by the auxiliary BDF2 scheme   \eqref{sch:Non_AE}.
\begin{thm}\label{thm:Conver_Non} 
	Assume that the solution of \eqref{Model:Bihar}--\eqref{Bound:Bihar} satisfy the regularity assumption condition \eqref{Regu:Bihar} and the adjacent temporal stepsize ratios $ r_{k} $ satisfy $ 0 < r_{k} < 4.8645 $. If the maximum temporal stepsize $
	\tau \leq \frac{1}{4 \sqrt{3} K_{f} } $,
	the following estimate holds for the auxiliary nonlinear compact difference scheme \eqref{sch:Non_AE} that
	\begin{equation*}\label{non:err_L2}
 			\|U^{n} - \bar{u}^{n}\|_{h} \leq C_{8} \left(\tau^{2} + h^{4} \right) \quad {\rm for } \  1 \leq n \leq N,
 	\end{equation*}
	where $ C_{8} := 12 \exp \big(4\sqrt{3}K_{f} T\big)  \max\{ 2C_{5} + C_{6} T, 3 C_{7} T \} $.
\end{thm}
{\bf Proof.}
Subtracting \eqref{sch:Non_AE} from \eqref{Non_eq_1}, we can get the following error equation 
\begin{equation}\label{EE:non_e1}
	\begin{aligned}
		\mathcal{D}_2 \mathcal{A}_{h} \bar{e}^{n}_{i,j}  - c \Lambda_{h} \bar{e}^{n}_{i,j} = \mathcal{A}_{h} \bar{f}  ( U^{n}_{i,j} ) -  \mathcal{A}_{h} \bar{f} ( \bar{u}^{n}_{i,j} ) +  R^{n}_{i,j}, \quad 1 \leq n \leq N.
	\end{aligned}
\end{equation}
for $ 1 \leq i \leq N^{h}_{x}-1 $ and $ 1 \leq j \leq N^{h}_{y}-1 $. Analogous to the proof of  \eqref{Unique_e1}--\eqref{Unique_e2}, it follows from the orthogonality \eqref{Def:DOC} of the DOC kernels $\{\theta^{(k)}_{k-m} \} $ that
\begin{equation}\label{EE:non_e2}
	\begin{aligned}
		\nabla_{\tau} \mathcal{A}_{h} \bar{e}^{k}_{i,j} - c \sum^{k}_{m=1} \theta^{(k)}_{k-m} \Lambda_{h} \bar{e}^{m}_{i,j} = \sum^{k}_{m=1} \theta^{(k)}_{k-m} \mathcal{A}_{h} \left( \bar{f} (U^{m}_{i,j}) - \bar{f} (\bar{u}^{m}_{i,j}) \right) + \sum^{k}_{m = 1} \theta^{(k)}_{k-m} R^{m}_{i,j},
	\end{aligned}
\end{equation}
which, by taking the inner product with $ 2 \bar{e}^{k} $ and summing the resulting
equality from $ k = 1 $ to $ n $, gives 
\begin{equation}\label{EE:non_e3}
	\begin{aligned}
		\| \bar{e}^{n} \|_{\mathcal{A},h}^{2} - 2 c \sum^{n}_{k=1} \sum^{k}_{m=1} \theta^{(k)}_{k-m} ( \Lambda_{h} \bar{e}^{m}, \bar{e}^{k} )_{h} & \leq 2 \sum^{n}_{k=1} \sum^{k}_{m=1} \theta^{(k)}_{k-m} \left(  \mathcal{A}_{h} f(U^{m})  - \mathcal{A}_{h} f(\bar{u}^{m}) , \bar{e}^{k} \right)_{h} \\
		& \qquad  + 2 \sum^{n}_{k=1} \sum^{k}_{m=1} \theta^{(k)}_{k-m} \left( R^{m}, \bar{e}^{k} \right)_{h}.
	\end{aligned}
\end{equation}

With a similar treatment to \eqref{Unique_e4}--\eqref{Unique_e5}, the above inequality yields
\begin{equation}\label{EE:non_e4}
	\begin{aligned}
		\| \bar{e}^{n} \|_{\mathcal{A},h}^{2} \leq 2 \sqrt{3} K_{f} \sum^{n}_{k=1} \| \bar{e}^{k} \|_{\mathcal{A},h} \sum^{k}_{m=1} \theta^{(k)}_{k-m} \| \bar{e}^{m} \|_{\mathcal{A}, h} + 2 \sqrt{3} \sum^{n}_{k=1} \| \bar{e}^{k} \|_{\mathcal{A},h} \sum^{k}_{m=1} \theta^{(k)}_{k-m} \| R^{m} \|_{h}. 
	\end{aligned}
\end{equation}
Similarly, if $ \| \bar{e}^{n^{*}} \|_{\mathcal{A},h} = \max_{ 1 \leq k \leq n } \| \bar{e}^{n} \|_{\mathcal{A},h}$, using Lemma \ref{lem:DOC_estimate} we have
\begin{equation}\label{EE:non_e7}
	\begin{aligned}
		\| \bar{e}^{n} \|_{\mathcal{A},h} \leq \| \bar{e}^{n^{*}} \|_{\mathcal{A},h} & \leq 2 \sqrt{3} K_{f} \sum^{n^{*}}_{k=1} \| \bar{e}^{k} \|_{\mathcal{A},h} \sum^{k}_{m=1} \theta^{(k)}_{k-m} + 2 \sqrt{3} \sum^{n^{*}}_{k=1} \sum^{k}_{m=1} \theta^{(k)}_{k-m} \| R^{m} \|_{h} \\
		& \leq 2 \sqrt{3} K_{f} \sum^{n}_{k=1} \tau_{k} \| \bar{e}^{k} \|_{\mathcal{A},h} + 2 \sqrt{3} \sum^{n}_{k=1} \sum^{k}_{m=1} \theta^{(k)}_{k-m} \| R^{m} \|_{h}. 
	\end{aligned}
\end{equation}
For the last term of \eqref{EE:non_e7}, we exchange the order of summation and utilize Lemma \ref{lem:DOC_estimate} to obtain 
\begin{equation*}
	\begin{aligned}
	     \sum^{n}_{k=1} \sum^{k}_{m=1} \theta^{(k)}_{k-m} \| R^{m} \|_{h} & = \sum^{n}_{m=1} \| R^{m} \|_{h} \sum^{n}_{k=m} \theta^{(k)}_{k-m} = \| R^{1} \|_{h} \sum^{n}_{k=1} \theta^{(k)}_{k-1} + \sum^{n}_{m=2} \| R^{m} \|_{h} \sum^{n}_{k=m} \theta^{(k)}_{k-m} \\
	     & \leq 2 \tau \| R^{1} \|_{h} + \max_{2 \leq m \leq n} \| R^{m} \|_{h} \sum^{n}_{m=2} \sum^{n}_{k=m} \theta^{(k)}_{k-m} \leq 2 \tau \| R^{1} \|_{h} + t_{n} \max_{2 \leq m \leq n} \| R^{m} \|_{h}, 
	\end{aligned}
\end{equation*}
which, together with the triangle inequality  and \eqref{Trun_err_t}--\eqref{Trun_err_s}, gives 
\begin{equation}\label{EE:non_e7b}
	\begin{aligned}
		\sum^{n}_{k=1} \sum^{k}_{m=1} \theta^{(k)}_{k-m} \| R^{m} \|_{h} \leq \big( 2C_{5} + C_{6} t_{n} \big)\ \tau^{2} + C_{7}\ \big( 2 \tau + t_{n} \big) h^{4}. 
	\end{aligned}
\end{equation}

Now, insert the above inequality into \eqref{EE:non_e7}, for $ \tau \leq 1 / (4\sqrt{3} K_{f}) $, we have 
\begin{equation}\label{EE:non_e8}
	\begin{aligned}
		\| \bar{e}^{n} \|_{\mathcal{A},h} \leq 4 \sqrt{3} K_{f} \sum^{n-1}_{k=1} \tau_{k} \| \bar{e}^{k} \|_{\mathcal{A},h} + ( 8\sqrt{3}C_{5} + 4\sqrt{3}C_{6} T ) \tau^{2} + 12\sqrt{3} C_{7} T h^{4}. 
	\end{aligned}
\end{equation}
Then, applying the discrete Gr\"{o}nwall inequality to \eqref{EE:non_e8} and using Lemma \ref{lem:cpmpact1} yields
\begin{equation*}
	\begin{aligned}
	\| \bar{e}^{n} \|_{h} \le \sqrt{3}	\| \bar{e}^{n} \|_{\mathcal{A},h} \leq 12\exp(4\sqrt{3}K_{f} T) \left( ( 2C_{5} + C_{6} T ) \tau^{2} + 3 C_{7} T h^{4} \right) \leq C_{8} \left( \tau^{2} + h^{4} \right),
	\end{aligned}
\end{equation*}
which implies the theorem.
\qed

Finally, we would like to give the $L^{2}$ norm error estimates  between the exact solution and the numerical solution yielded by the variable-step BDF2 compact scheme   \eqref{sch:Non_CDM}. Note that by Theorem \ref{thm:Conver_Non} and Lemma \ref{lem:inverse}, we see
\begin{equation*}
	\begin{aligned}
		\| \bar{e}^{n} \|_{h,\infty} \leq C_{0}  h^{-1}\|\bar{e}^{n}\|_{h} \leq C_{0} C_{8} h^{-1}(\tau^{2} + h^{4} ).
	\end{aligned}
\end{equation*}
Consequently, if $ \tau = o(h^{\frac{1}{2}}) $,  and for $ \tau $, $ h $ sufficiently small,  it holds that  $ \bar{u}^{n}_{i,j} \in B_{\delta} $, which further implies  $ \bar{f}( \bar{u}^{n}_{i,j} ) = f( \bar{u}^{n}_{i,j} ) $ and thus  in this context  $ \bar{u}^{n}_{i,j} \equiv u^{n}_{i,j} $. We summarize these conclusions in the following theorem.
\begin{thm}\label{thm:Conver_Non_New} 
	Under the conditions in Theorem \ref{thm:Conver_Non} and if the maximum temporal stepsize satisfies $ \tau = o( h^{\frac{1}{2}} ) $, the nonlinear compact difference scheme \eqref{sch:Non_CDM} admits a unique solution satisfying
	\begin{equation*}\label{non:err_L2_new}
		\begin{aligned}
			\|U^{n} - u^{n}\|_{h} \leq C_{8} \left(\tau^{2} + h^{4} \right) \quad {\rm for } \  1 \leq n \leq N.
		\end{aligned}
	\end{equation*}
\end{thm}

\section{ An efficient variable-step two-grid compact difference scheme  }\label{sec:TG}
In order to solve the nonlinear system \eqref{Model:Bihar}--\eqref{Bound:Bihar} efficiently, we shall propose a two-grid compact finite difference scheme based on the variable-step BDF2 method and the piecewise bi-cubic Lagrange interpolation developed in Section \ref{sec:Map_lag} as follows. 

\textbf{Step 1.} On the coarse grid,   solve a small-scale nonlinear compact  difference scheme to find a rough solution $u^{n}_{H}=\{u^{n}_{H,i,j}\}$ by
\begin{equation}\label{sch:coar_equ1}
	\mathcal{D}_{2} \mathcal{A}_{H} u^{n}_{H,i,j} - c \Lambda_{H} u^{n}_{H,i,j} = \mathcal{A}_{H} f(u^{n}_{H,i,j}) + \mathcal{A}_{H} g^{n}_{i,j}, \quad (i,j) \in  \omega_{H}, 
\end{equation}
 subject to the initial and boundary conditions \eqref{Initial:Bihar}--\eqref{Bound:Bihar} defined on coarse grid.

\textbf{Step 2.} On the fine grid,  solve a large-scale linearized compact difference scheme to produce a corrected solution $ u^{n}_{h}=\{u^{n}_{h,i,j}\} $ based on the rough solution $u^{n}_{H}$ in Step 1 by
\begin{equation}\label{sch:fine_equ1}
	\begin{aligned}
		\mathcal{D}_{2} \mathcal{A}_{h} u^{n}_{h,i,j} - c \Lambda_{h} u^{n}_{h,i,j} = \mathcal{A}_{h} F^{n}_{i,j} + \mathcal{A}_{h} g^{n}_{i,j},  \quad (i,j) \in  \omega_{h}, 
	\end{aligned}
\end{equation}
subject to the initial and boundary conditions \eqref{Initial:Bihar}--\eqref{Bound:Bihar},
where $F^{n}_{i,j}$ represents a Newton linearization from coarse grid to fine grid defined as
\begin{equation}\label{def:F}
	\begin{aligned}
		F^{n}_{i,j} := f\left( \left[\Pi_{H}u_{H}^{n}\right]_{i,j} \right) + f'\left( \left[\Pi_{H}u_{H}^{n}\right]_{i,j} \right)\left( u^{n}_{h,i,j} - \left[\Pi_{H}u_{H}^{n}\right]_{i,j} \right).
	\end{aligned}
\end{equation}

On the coarse grid, denote $e_{H}^{n}= U^{n}- u_{H}^{n}$ for  $ 0 \leq n \leq N $. Analogous to Theorem \ref{thm:Conver_Non_New}, we can immediately reach the following conclusions.
\begin{thm}\label{thm:Conver_Coarse} 
Under the conditions in Theorem \ref{thm:Conver_Non} and if the maximum temporal stepsize satisfies $ \tau = o(H^{\frac{1}{2}}) $, the nonlinear compact difference scheme \eqref{sch:coar_equ1} defined on the coarse grid admits a unique solution satisfying
\begin{equation*}\label{corase:err_L2}
	\begin{aligned}
		\|U^{n} -u_{H}^{n}\|_{H} \leq C_{8}\left( \tau^{2} + H^{4} \right) \quad {\rm for } \  1 \leq n \leq N.
	\end{aligned}
\end{equation*}
\end{thm}

Based on this theorem and Lemmas \ref{lem:Lagrange_bounds_L2}--\ref{lem:Lagrange_bounds_Inf}, the following corollary can be derived. 
\begin{cor}\label{cor:Inter_err} Under the conditions in Theorem \ref{thm:Conver_Coarse}, the numerical solution $u^{n}_{H}$ of the nonlinear compact difference scheme \eqref{sch:coar_equ1} satisfies
	\begin{align}\label{cubic:err_L2}
		&\|U^{n} - \Pi_{H}u^{n}_{H}\|_{h} \leq C_{9}  \left( \tau^{2} + H^{4} \right)  \quad  {\rm for } \  1 \leq n \leq N,\\\label{cubic:err_Linfty}
		&\|U^{n} - \Pi_{H}u^{n}_{H}\|_{h,\infty} \leq C_{10} \left(  H^{-1}\tau^{2} + H^{3}\right)  \quad  {\rm for } \  1 \leq n \leq N,
	\end{align}
	and consequently the interpolation solution $\Pi_{H} u_{H}^{n}$ is bounded on the fine grid and satisfies  
	\begin{equation}\label{cubic:bound}
		\Pi_{H} u_{H}^{n} \in B_\delta,
	\end{equation}
	where $ C_{9} = C_{2} + C_{3}C_{8} $, $ C_{10} := C_{1} + C_{0}C_{4}C_{8} $. 
\end{cor}
{\bf Proof.}
We perform the splitting 
$$
\|U^{n} - \Pi_{H}u^{n}_{H}\|_{h} \leq \| U^{n} - \Pi_{H}U^{n} \|_{h} + \| \Pi_{H} U^{n} - \Pi_{H}u^{n}_{H} \|_{h},
$$
where the first right-hand side term could be bounded as $\| U^{n} - \Pi_{H}U^{n} \|_{h} \leq C_{2} H^{4} $ via Lemma \ref{lem:cubic_error}, and due to the linear property of the piecewise bi-cubic Lagrange interpolation operator  with respect ro the interpolated function, Lemma \ref{lem:Lagrange_bounds_L2} and Theorem \ref{thm:Conver_Coarse} give us
$$
\| \Pi_{H} U^{n} - \Pi_{H}u^{n}_{H} \|_{h} = \| \Pi_{H}e^{n}_{H} \|_{h} \leq C_{3} \| e^{n}_{H} \|_{H} \leq  C_{3} C_{8} ( \tau^{2} + H^{4}) ,
$$
which completes the proof of estimate \eqref{cubic:err_L2}. 

Analogous to the proof of \eqref{cubic:err_L2}, by Lemmas \ref{lem:inverse}, \ref{lem:cubic_error}, \ref{lem:Lagrange_bounds_Inf} and Theorem \ref{thm:Conver_Coarse}, we can conclude that
\begin{equation*}
	\begin{aligned}
        \|U^{n} - \Pi_{H}u^{n}_{H}\|_{h,\infty} & \leq \| U^{n} - \Pi_{H}U^{n} \|_{h,\infty} + \| \Pi_{H} U^{n} - \Pi_{H}u^{n}_{H} \|_{h,\infty} \\
        & \leq C_{1} H^{4} + C_{4} \| e^{n}_{H} \|_{H,\infty} \leq C_{1} H^{4} + C_{0}C_{4} H^{-1} \| e^{n}_{H} \|_{H} \\
        & \leq ( C_{1} + C_{0}C_{4}C_{8} ) ( H^{-1}\tau^{2} + H^{3} ) .
		\end{aligned}
\end{equation*}
Furthermore, an application of the triangle inequality yields 
$$
\|\Pi_{H} u_{H}^{n}\|_{h,\infty}  \leq \| U^{n} \|_{h,\infty} + \| U^{n} - \Pi_{H} U^{n}_H \|_{h,\infty} \leq \| U^{n}\|_{h,\infty} + C_{10} ( H^{-1} \tau^{2} + H^{3} ),
$$
which completes the proof of \eqref{cubic:bound}.
\qed

\begin{rem}\label{rem:cubic_solution} \rm Corollary \ref{cor:Inter_err} declares that the interpolation solution $ \Pi_{H} u^{n}_{H} $ is bounded under the discrete $ L^{\infty} $ norm on the fine grid, if  the maximum temporal stepsize satisfies $ \tau = o(H^{\frac{1}{2}}) $. This  restriction  is compatible with  previous work, see Refs. \cite[Theorem 4.2]{SINUM_Dawson_1998} and \cite[Theorem 4.2]{2022_ANM_Fu}. Notice that the boundedness of $ \Pi_{H} $ in the sense of $ L^{\infty} $ norm in Lemma \ref{lem:Lagrange_bounds_Inf} plays an important role in the proof of \eqref{cubic:err_Linfty}--\eqref{cubic:bound}. Actually, if Lemma \ref{lem:Lagrange_bounds_Inf} does not hold, we have
\begin{equation*}
	\begin{aligned}
		\| \Pi_{H} U^{n} - \Pi_{H}u^{n}_{H} \|_{h,\infty} \leq C_{0} h^{-1} \| \Pi_{H}e^{n}_{H} \|_{h} \leq C_{0}C_{3}C_{8} h^{-1} ( \tau^{2} + H^{4} ),
	\end{aligned}
\end{equation*}
which leads to a much worse restriction $ \tau = o(h^{\frac{1}{2}}) $. On the other hand, in practical computation the condition $ \tau =o(H^{\frac{1}{2}}) $ is
not too severe, as the coarse grid size $H$ is large enough compared with the fine grid size $h$, and a mild choice $\tau=\mathcal{O} ( H^{2})$ suggested by Theorem \ref{thm:Conver_Coarse} can naturally satisfy the restrictive condition.
\end{rem}

At last, we shall give an error estimate for $e_{h}^{n}= U^{n} - u_{h}^{n}$ on the fine grid for the linearized compact difference scheme \eqref{sch:fine_equ1}.
\begin{thm}\label{thm:Conver_Fine}
	Under the conditions in Theorem \ref{thm:Conver_Coarse}, the following error estimate holds 
for	the solution of the two-grid compact difference scheme \eqref{sch:coar_equ1}--\eqref{sch:fine_equ1}  
	$$
	\|U^{n} - u_{h}^{n}\|_{h} \leq C_{11} \left( \tau^{2} + h^{4} + H^{7} \right)  \quad {\rm for}  \  1 \leq n \leq N, 
	$$ 
    where $ C_{11} := 2 \sqrt{3} \exp \big(4\sqrt{3}K_{f} T\big) \max\{ 4\sqrt{3}C_{5} + 2\sqrt{3}C_{6} T +  C_{9} C_{10} K_{f} T, 6\sqrt{3} C_{7} T \}$.
\end{thm}
{\bf Proof.}
For the linearized scheme \eqref{sch:fine_equ1} on the fine grid, we can get a very similar error equation 
\begin{equation}\label{EE:TG_e1}
	\begin{aligned}
		\mathcal{D}_{2} \mathcal{A}_{h} e^{n}_{h,i,j}- c \Lambda_{h} e^{n}_{h,i,j}= \mathcal{A}_{h} f( U^{n}_{i,j} ) - \mathcal{A}_{h} F^{n}_{i,j} + R^{n}_{i,j}, \quad (i,j)\in \omega_h,
	\end{aligned}
\end{equation}
which, together with a similar treatment to \eqref{EE:non_e1}--\eqref{EE:non_e3}, leads to
\begin{equation}\label{EE:TG_e2}
	\begin{aligned}
		\| e_{h}^{n} \|_{\mathcal{A},h}^{2} - 2 c \sum^{n}_{k=1} \sum^{k}_{m=1} \theta^{(k)}_{k-m} ( \Lambda_{h} e_{h}^{m}, e_{h}^{k} ) & \leq 2 \sum^{n}_{k=1} \sum^{k}_{m=1} \theta^{(k)}_{k-m} \left(  \mathcal{A}_{h} f( U^{m} ) - \mathcal{A}_{h} F^{m}, e_{h}^{k} \right) \\
		& \quad  + 2 \sum^{n}_{k=1} \sum^{k}_{m=1} \theta^{(k)}_{k-m} \left( R^{m}, e_{h}^{k} \right).
	\end{aligned}
\end{equation}

Compared with the estimate \eqref{EE:non_e3} in Theorem \ref{thm:Conver_Non}, the only difference lies in the treatment of the nonlinear term. We apply Taylor expansion of $f( U^{m} ) $ at $ \Pi_{H} u_{H}^m $ to obtain
\begin{equation*}
	\begin{aligned}
		f(U_{i,j}^{m}) & = f\left( \left[\Pi_{H} u_{H}^{m} \right]_{i,j} \right) + f'\left( \left[\Pi_{H} u_{H}^{m} \right]_{i,j} \right) \left(U^{m}_{i,j} - \left[\Pi_{H} u_{H}^{m} \right]_{i,j} \right)  + \frac{1}{2} f''(\mu_{i,j}^{m}) \left(U^{m}_{i,j} - \left[\Pi_{H} u_{H} \right]^{m}_{i,j} \right)^{2},
	\end{aligned}
\end{equation*}
for some constant $ \mu^{m}_{i,j} $ between $ U^{m}_{i,j} $ and $[\Pi_{H} u_{H}^{m} ]_{i,j} $. Then, subtract $F^{m}_{i,j}$ in \eqref{def:F} from this equation, we have
\begin{equation*}
	\begin{aligned}
		f(U_{i,j}^{m}) - F^{m}_{i,j} = f'\left(  [\Pi_{H} u_{H}^{m} ]_{i,j} \right) \left(  U^{m}_{i,j} - u_{h,i,j}^{m}\right) + \frac{1}{2} f''(\mu_{i,j}^{m}) \left(  U^{m}_{i,j} - \left[\Pi_{H} u_{H}^{m} \right]_{i,j} \right)^{2}.
	\end{aligned}
\end{equation*}

By Corollary \ref{cor:Inter_err}, if $ \tau = o(H^{\frac{1}{2}}) $, we have $ \Pi_{H} u_{H}^{m}\in B_\delta$ and $\mu^{m} \in B_\delta$. Thus, assumption \eqref{Regu:Bihar:e1}
implies that $\big\vert  f'( [\Pi_{H} u_{H}^{m} ]_{i,j} )\big\vert  + \big\vert  f''(\mu_{i,j}^{m}) \big\vert \leq K_{f}$.
Furthermore, we apply Lemma \ref{lem:cpmpact1}, Cauchy-Schwarz inequality and  Corollary \ref{cor:Inter_err}  to derive
\begin{equation}\label{EE:TG_e3}
	\begin{aligned}
		& \quad \left(  \mathcal{A}_{h} f( U^{m} ) - \mathcal{A}_{h} F^{m}, e_{h}^{k} \right) \leq \| f( U^{m} ) - F^{m} \|_{h} \| e_{h}^{k} \|_{\mathcal{A},h} \\
		& \leq K_{f}  \| e_{h}^{m} \|_{h} \| e_{h}^{k} \|_{\mathcal{A},h} + \frac{1 }{2} K_{ f }\|  U^{m} - \Pi_{H} u_{H}^{m}  \|_{h,\infty} \| U^{m} - \Pi_{H} u_{H}^{m} \|_{h} \| e_{h}^{k} \|_{\mathcal{A},h} \\
%
       & \leq \sqrt{3} K_{f} \| e_{h}^{m} \|_{\mathcal{A},h} \| e_{h}^{k} \|_{\mathcal{A},h} + \frac{1}{2} C_{9} C_{10} K_{f} \big( \tau^2 + H^{7} \big) \| e^{k}_{h} \|_{\mathcal{A},h}.
	\end{aligned}
\end{equation} 
The other terms in \eqref{EE:TG_e2} could be estimated similarly as in the proof of Theorem \ref{thm:Conver_Non}. Then, substituting \eqref{EE:TG_e3} into \eqref{EE:TG_e2} gives us 
\begin{equation*}
	\begin{aligned}
		\| e_{h}^{n} \|_{\mathcal{A},h}^{2} & \leq 2 \sqrt{3} K_{f} \sum^{n}_{k=1} \| e_{h}^{k} \|_{\mathcal{A},h} \sum^{k}_{m=1} \theta^{(k)}_{k-m} \| e_{h}^{m} \|_{\mathcal{A}, h} +  C_{9} C_{10} K_{f} \big( \tau^2 + H^{7} \big) \sum^{n}_{k=1} \| e^{k}_{h} \|_{\mathcal{A},h} \sum^{k}_{m=1} \theta^{(k)}_{k-m}   \\
		& \qquad + 2 \sqrt{3} \sum^{n}_{k=1} \| e_{h}^{k} \|_{\mathcal{A},h} \sum^{k}_{m=1} \theta^{(k)}_{k-m} \| R^{m} \|_{h}. 
	\end{aligned}
\end{equation*}

As was done before, choosing an integer $ n^{*}$ $ ( 1 \leq n^{*} \leq n )$ such that $ \| e_{h}^{n^{*}} \|_{\mathcal{A},h} = \max_{ 1 \leq k \leq n } \| e_{h}^{n} \|_{\mathcal{A},h}$, we have 
\begin{equation*}
	\begin{aligned}
		\| e_{h}^{n} \|_{\mathcal{A},h} & \leq 2 \sqrt{3} K_{f} \sum^{n}_{k=1} \| e_{h}^{k} \|_{\mathcal{A},h} \sum^{k}_{m=1} \theta^{(k)}_{k-m} 
		         +  C_{9} C_{10} K_{f} \big( \tau^2 + H^{7} \big) \sum^{n}_{k=1} \sum^{k}_{m=1} \theta^{(k)}_{k-m}    
		          + 2 \sqrt{3} \sum^{n}_{k=1} \sum^{k}_{m=1} \theta^{(k)}_{k-m} \| R^{m} \|_{h}. 
	\end{aligned}
\end{equation*}
Therefore, for $ \tau \leq 1 / (4\sqrt{3} K_{f})$, we use Lemma \ref{lem:DOC_estimate} and estimate \eqref{EE:non_e7b} to obtain
\begin{equation}\label{EE:non_e9}
	\begin{aligned}
		 \| e_{h}^{n} \|_{\mathcal{A},h} 
		 & \leq 4 \sqrt{3} K_{f} \sum^{n-1}_{k=1} \tau_{k} \| e^{k} \|_{\mathcal{A},h} + ( 8\sqrt{3}C_{5} + 4\sqrt{3}C_{6} T + 2 C_{9} C_{10} K_{f} T ) \tau^{2} \\
		 & \qquad
		 + 12 \sqrt{3} C_{7} T h^{4}+ 2 C_{9} C_{10} K_{f} T H^{7},  
	\end{aligned}
\end{equation}
from which an application of the discrete Gr\"{o}nwall inequality in Lemma \ref{lem:Gronwall} and Lemma \ref{lem:cpmpact1} yield the following estimate
\begin{equation*}
	\begin{aligned}
		\| e_{h}^{n} \|_{h} & \leq \sqrt{3} \| e_{h}^{n} \|_{\mathcal{A},h} \\
		&\leq 2\sqrt{3} \exp(4\sqrt{3}K_{f} T) \left( ( 4\sqrt{3}C_{5} + 2\sqrt{3} C_{6} T +  C_{9} C_{10} K_{f} T ) \tau^{2} 
		+ 6\sqrt{3} C_{7} T h^{4}  +  C_{9} C_{10} K_{f} T H^{7}\right) \\
		& \leq  C_{11} ( \tau^{2} + h^{4} + H^{7} ).
	\end{aligned}
\end{equation*}
The proof is completed. 
\qed

\section{Extension to periodic boundary condition}\label{sec:PBC}
In this section, we extend the ideas and derivations in previous sections to the semilinear parabolic equation \eqref{Model:Bihar}--\eqref{Bound:Bihar}  with periodic boundary condition. Firstly, we denote the following spaces of grid functions on grids $ \bar{\omega}_{\kappa} $
$$
\mathcal{V}^{p}_{\kappa} = \left\{v | v \in \mathcal{V}_{\kappa} \  \  \text {and} \  \   v \  \text{is periodic} \right\}.
$$
Furthermore, for any grid functions $ w, q \in \mathcal{V}^{p}_{\kappa} $, the discrete inner product and corresponding norms are redefined as

\[
	(w,q)_{\kappa} = \kappa_{x} \kappa_{y} \sum^{N^{\kappa}_{x}}_{i=1} \sum^{N^{\kappa}_{y}}_{j=1} w_{i,j} q_{i,j}, ~~
	\| w \|_{\kappa} = \sqrt{(w, w)_{\kappa}}, ~~ \|w\|_{\mathcal{A},\kappa} = \sqrt{ ( \mathcal{A}_{\kappa} w, w )_{\kappa}},
\]
and a useful lemma is listed below.
\begin{lem}[\cite{2010_SINUM_Liao}]\label{lemp:cpmpact1} For any $ w \in \mathcal{V}^{p}_{\kappa} $, we have
	$
	\frac{2}{3} \| w \|_{\kappa} \leq \| w \|_{\mathcal{A},\kappa} \leq \| w \|_{\kappa}. 
	$
\end{lem}

In the context of periodic boundary case, we shall still adopt the piecewise bi-cubic Lagrange interpolation operator defined in Section \ref{sec:Map_lag} to construct 
high-order two-grid difference scheme. However, a small modification of the proof of Lemmas \ref{lem:Lagrange_bounds_L2}--\ref{lem:Lagrange_bounds_Inf} should be given to show the boundedness conclusions of the piecewise bi-cubic Lagrange interpolation operator under the discrete $L^{2}$ and $ L^{\infty} $ norms.

\begin{lem}\label{lemp:Lagrange_bounds_L2} For any $w\in\mathcal{V}_H^{p}$, the following estimate holds
	\begin{equation*}
		\begin{aligned}
			\| \Pi_{H} w \|_{h} \leq C_{3} \| w \|_{H},
		\end{aligned}
	\end{equation*}
	where $ C_{3} = 4 \left(\frac{3 + 27^{2} + ( 10 + 7\sqrt{7} )^{2}}{27^{2}}\right) $. 
\end{lem}
{\bf Proof.}
Denote $ \xi_{i,j} = \Pi_{H,x} w_{i,j} $, and then $ \| \Pi_{H} w \|^{2}_{h} $ can be rewritten as 
\begin{equation}\label{p_L2_Bound_e1}
	\begin{aligned}
		\| \Pi_{H} w \|^{2}_{h} = h_{x} h_{y} \sum^{N^{h}_{x}}_{i=1} \sum^{N^{h}_{y}}_{j=1} \left(  \Pi_{H} w \right) _{i,j}^{2} = h_{x} \sum^{N^{h}_{x}}_{i=1} \left( I_{1} + I_{2} + I_{4} \right),
	\end{aligned}
\end{equation}
where $ I_{1} $, $ I_{2} $ are defined in Lemma \ref{lem:Lagrange_bounds_L2} and 
\begin{equation*}
	\begin{aligned}
		I_{4} := h_{y} \sum^{ N^{h}_{y} }_{j = M_{y}(N^{H}_{y}-1) + 1 } \left( \Pi_{H,y} \xi_{i,j} \right)^{2}.
	\end{aligned}
\end{equation*}
While, the application of Lemma \ref{lem:basis_bounds} shows
\begin{equation*}
	\begin{aligned}
		I_{4} \leq 4 H_{y} \left( \frac{3}{27^{2}} \xi_{i,N^{H}_{y}-3}^{2} + \frac{(7\sqrt{7}-10)^{2}}{27^{2}} \xi_{i,N^{H}_{y}-2}^{2} + \frac{(7\sqrt{7}+10)^{2}}{27^{2}} \xi_{i,N^{H}_{y}-1}^{2} + \xi_{i,N^{H}_{y}}^{2} \right).
	\end{aligned}
\end{equation*}
This together with \eqref{L2_Bound_e4}--\eqref{L2_Bound_e5} implies
\begin{equation}\label{p_L2_Bound_e7}
	\begin{aligned}
		I_{1} + I_{2} + I_{4} \leq C_{3} H_{y} \sum^{N^{H}_{y}-1}_{j=1} \xi_{i,j}^{2} + H_{y} 
		\frac{24 + 8 \times 27^{2} }{27^{2}} \xi_{i,N^{H}_{y}}^{2} \leq C_{3} H_{y} \sum^{N^{H}_{y}}_{j=1} \xi_{i,j}^{2}.
	\end{aligned}
\end{equation}

Finally, the claimed result can be derived immediately by following a completely similar process as \eqref{L2_Bound_e8}--\eqref{L2_Bound_e10}.
\qed

\begin{lem}\label{lemp:Lagrange_bounds_Inf} For any $w\in\mathcal{V}_H^{p}$, the following estimate holds
	\begin{equation*}
		\begin{aligned}
			\| \Pi_{H} w \|_{h,\infty} \leq C_{4} \| w \|_{H,\infty},
		\end{aligned}
	\end{equation*}
	where $ C_{4} = \left( \frac{ 54 + 2\sqrt{3} }{27} \right)^{2} $. 
\end{lem}
{\bf Proof.} Similar as the proof of Lemma \ref{lem:Lagrange_bounds_Inf}, estimate \eqref{Inf_eq1} holds for  $  M_{y}  \leq j^{*} \leq (N^{H}_{y}-1) M_{y}$.  Besides, if $ 1  \leq j^{*} \leq M_{y} $, we have
\begin{equation*}
	\begin{aligned}
		\| \Pi_{H} w \|_{h,\infty} & \leq \frac{ 7 \sqrt{7} + 10 }{27} \left\vert \xi_{i^{*},1} \right\vert + \frac{ 7 \sqrt{7} - 10 }{27} \left\vert \xi_{i^{*},2} \right\vert + \frac{ \sqrt{3} }{27} \left\vert \xi_{i^{*},3} \right\vert + \frac{ \sqrt{3} }{27} \left\vert \xi_{i^{*},N_{y}^{H}} \right\vert\\
		& \leq \frac{ 54 + 2\sqrt{3} }{27} \max_{ 1 \leq k \leq N^{H}_{y} } \left\vert \xi_{i^{*},k} \right\vert,
	\end{aligned}
\end{equation*}
and similarly, if $ (N^{H}_{y} - 1)M_{y}  \leq j^{*} \leq N^{h}_{y} $, we have
\begin{equation*}
	\begin{aligned}
		\| \Pi_{H} w \|_{h,\infty} \leq \frac{ 54 + 2\sqrt{3} }{27} \max_{ 1 \leq k \leq N^{H}_{y} } \left\vert \xi_{i^{*},k} \right\vert.
	\end{aligned}
\end{equation*}

Finally, analogous to \eqref{Inf_eq5},  we see $ \left\vert \xi_{i^{*},k} \right\vert \leq \frac{ 54 + 2\sqrt{3} }{27} \max_{ 1 \leq i \leq N^{H}_{x} } \left\vert w_{i,k} \right\vert $ holds, which completes the proof.
\qed

Now, an efficient two-grid fourth-order compact difference scheme for model \eqref{Model:Bihar}--\eqref{Bound:Bihar} is proposed similarly as follows.

\textbf{Step 1.} On the coarse grid, solve a small-scale nonlinear compact finite difference scheme to find a rough solution $u^{n}_{H}$ by
\begin{equation}\label{schp:coar_equ1}
	\mathcal{D}_{2} \mathcal{A}_{H} u^{n}_{H,i,j} - c \Lambda_{H} u^{n}_{H,i,j} = \mathcal{A}_{H} f(u^{n}_{H,i,j}) + \mathcal{A}_{H} g^{n}_{i,j}, \quad (i,j) \in  \omega_{H}, 
\end{equation}
 subject to the initial condition \eqref{Initial:Bihar} and periodic boundary condition.

\textbf{Step 2.} On the fine grid, solve a large-scale linearized compact difference scheme to produce a corrected solution $ u^{n}_{h} $ based on the rough solution $u^{n}_{H}$ in Step 1 by
\begin{equation}\label{schp:fine_equ1}
	\begin{aligned}
		\mathcal{D}_{2} \mathcal{A}_{h} u^{n}_{h,i,j} - c \Lambda_{h} u^{n}_{h,i,j} = \mathcal{A}_{h} F^{n}_{i,j} + \mathcal{A}_{h} g^{n}_{i,j}, \quad (i,j) \in  \omega_{h}, 
	\end{aligned}
\end{equation}
 subject to the initial condition \eqref{Initial:Bihar} and periodic boundary condition.

Following the proofs of Corollary \ref{cor:Inter_err}, Theorems \ref{thm:Conver_Coarse} and \ref{thm:Conver_Fine}, together with Lemmas \ref{lemp:cpmpact1}--\ref{lemp:Lagrange_bounds_Inf}, the unique
solvability and error estimates for the two-grid algorithm \eqref{schp:coar_equ1}--\eqref{schp:fine_equ1} can be proved very similarly, and we skip the detailed proof here.
\begin{thm}\label{thmp:Conver_Non} 
	Assume that the solution of \eqref{Model:Bihar}--\eqref{Bound:Bihar} satisfy regularity assumption condition \eqref{Regu:Bihar} and the adjacent temporal stepsize ratios $ r_{k} $ satisfy $ 0 < r_{k} < 4.8645 $. If the maximum temporal stepsize $ \tau = o( H^{\frac{1}{2}} ) $ and $
	\tau \leq \frac{1}{6 K_{f} } $, the two-grid compact difference scheme \eqref{schp:coar_equ1}--\eqref{schp:fine_equ1} admits a unique solution satisfying
	$$
	\|U^{n} - u_{h}^{n}\|_{h}  \leq C_{12}\left ( \tau^{2} + h^{4} + H^{7}\right) \quad {\rm for}  \  1 \leq n \leq N, 
	$$ 
	where $ C_{12} := 3 \exp(6K_{f} T) \max\{ 6C_{5} + 3C_{6} T + C_{9} C_{10} K_{f} T, 9 C_{7} T \}$.
\end{thm}

\section{Numerical examples}\label{sec:exam}
In this section, we shall present several numerical experiments to test the effectiveness and efficiency of the variable-step two-grid compact difference scheme. 
In the computation, a Newton-type iterative procedure with tolerance error $1.0\times 10^{-13} $ is performed to solve the nonlinear algebra systems at each time level. 

\subsection{Accuracy and numerical stability tests on uniform temporal grids}
\begin{figure}[htbp]
	\centering
		\begin{minipage}[t]{0.65\linewidth}
			\centering
			\includegraphics[width=2.75in]{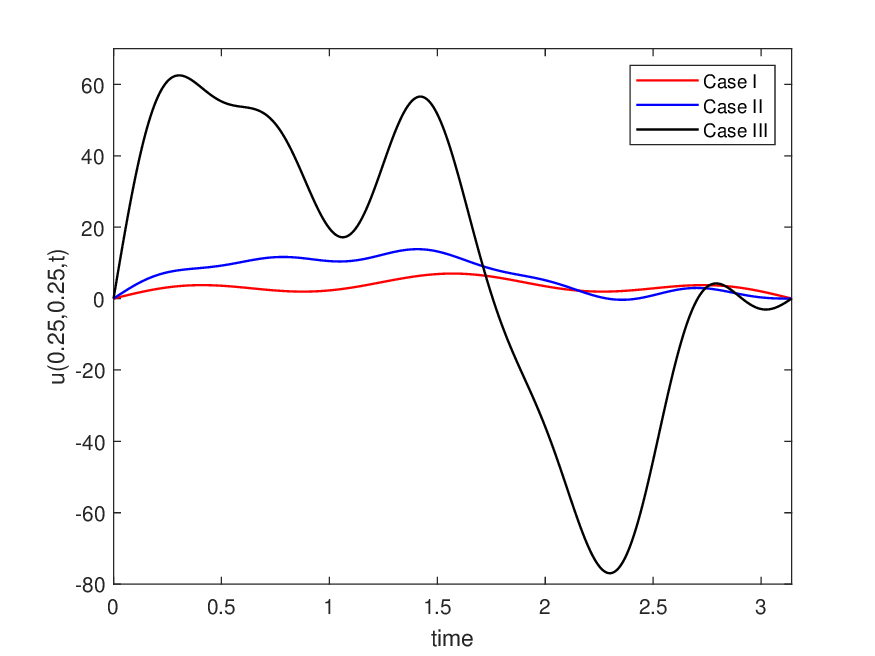}
				\caption{ Evolution of exact solution with time at point $(0.25,0.25) $}
			\label{figEX1_1}
		\end{minipage}
\end{figure}

In this subsection, we shall compare the numerical accuracy of the two-grid compact difference scheme \eqref{sch:coar_equ1}--\eqref{sch:fine_equ1} with the standard nonlinear scheme \eqref{sch:Non_CDM} as well as the following implicit-explicit scheme
\begin{equation}\label{sch:IM}
	\begin{aligned}
		\mathcal{D}_2 \mathcal{A}_{h} u^{n}_{i,j}  - c \Lambda_{h} u^{n}_{i,j} = \mathcal{A}_{h} f ( u^{n,*}_{i,j} ) + \mathcal{A}_{h} g^{n}_{i,j},
	\end{aligned}
\end{equation}
where
\begin{equation*}
	\begin{array}{l}
			u^{n,*} :=
			\left\{
			\begin{split}
					& 2 u^{n-1} - u^{n-2}, \qquad  n \geq 2,\\
					& u^{0}, \qquad \qquad \quad \quad \ n = 1.
				\end{split}
			\right.
		\end{array}
\end{equation*}
For this purpose, we consider the following semilinear parabolic equation
\begin{equation}\label{Ex:1}
	\begin{aligned}
		u_{t} - \Delta u = u - u^3 + g, \quad (x,y) \in (0,1)^2, \  t \in (0,\pi],
	\end{aligned}
\end{equation}
where the linear part $ g(x,y,t) $ is given such that the exact solution is one of the following three types, i.e.,
\begin{itemize}
	\item Case I:    $ u(x,y,t)=\left[ 5 \sin\left(t\right) + 2 \sin\left(5 t\right) \right] \sin\left(2\pi x\right)\sin\left(2\pi y\right) $;
	
	\item Case II:    $u(x,y,t)= \left[  10 \sin\left(t\right) + 5 \sin\left(2t\right) + 2 \sin\left(5 t\right) + \sin\left(10 t\right) \right]  \sin\left(2\pi x\right)\sin\left(2\pi y\right) $;
	
	\item Case III:    $u(x,y,t)= \left[  10 \sin\left(t\right) + 50 \sin\left(2t\right) + 30 \sin\left(5 t\right) + 10 \sin\left(10 t\right) \right]  \sin\left(2\pi x\right)\sin\left(2\pi y\right) $.
\end{itemize}
Figure \ref{figEX1_1} displays the evolution of these solutions with time, in which it can be clearly observed that the solution in Case I changes most smoothly while the solution in Case III changes most sharply with respect to time.

Firstly, we set $ N^{h}_{x} = N^{h}_{y} $ and $ M_{x} = M_{y} = 10 $ and adjust $N$ and $ N^{h}_{x}$ to investigate the spatial convergence of these three methods. Numerical results for Cases I--III are listed in Tables \ref{Ex:tabS1}--\ref{Ex:tabS3} respectively, which indicates the fourth-order spatial accuracy for both the nonlinear method and the two-grid method. But the implicit-explicit method is only successfully implemented for Case I and fails for the other two cases, which may be caused by improper treatment of the time dependence and nonlinearity.
\begin{table}[!tbp]
	\caption{Spatial convergence of nonlinear scheme \eqref{sch:Non_CDM}, two-grid scheme \eqref{sch:coar_equ1}--\eqref{sch:fine_equ1} and implicit-explicit scheme \eqref{sch:IM} for Case I on uniform temporal grids\label{Ex:tabS1}}%
	{\footnotesize\begin{tabular*}{\columnwidth}{@{\extracolsep\fill}ccccccc@{\extracolsep\fill}}
			\toprule
			& \multicolumn{2}{c}{Nonlinear scheme}  & \multicolumn{2}{c}{Two-grid scheme} & \multicolumn{2}{c}{Implicit-explicit scheme}\\
			\midrule
			$ (N, N^{h}_{x}) $ &  Error  &  Order  &  Error  &  Order  &  Error  &  Order \\
			\midrule 
			$ (80, 100) $  &  $8.16 \times 10^{-4}$    &   ---      &  $8.16 \times 10^{-4}$    &  ---    &  $9.40 \times 10^{-4}$  & ---       \\ 
			$ (180, 150) $  &  $1.64 \times 10^{-4}$   &   3.96      &  $1.64 \times 10^{-4}$    &  3.96   &  $1.77 \times 10^{-4}$  & 4.12      \\ 
			$ (320, 200) $  &  $5.21 \times 10^{-5}$    &   3.99      &  $5.21 \times 10^{-5}$    &  3.99   &  $5.49 \times 10^{-5}$   & 4.07       \\ 
			$ (500, 250) $  &  $2.14 \times 10^{-5}$    &   3.99      &  $2.14 \times 10^{-5}$    &  3.99   &  $2.23 \times 10^{-5}$  & 4.04      \\ 
			\bottomrule
	\end{tabular*}}
\end{table}
\begin{table}[!tbp]
	\caption{Spatial convergence of nonlinear scheme \eqref{sch:Non_CDM}, two-grid scheme \eqref{sch:coar_equ1}--\eqref{sch:fine_equ1} and implicit-explicit scheme \eqref{sch:IM} for Case II on uniform temporal grids\label{Ex:tabS2}}%
	{\footnotesize\begin{tabular*}{\columnwidth}{@{\extracolsep\fill}ccccccc@{\extracolsep\fill}}
			\toprule
			& \multicolumn{2}{c}{Nonlinear scheme}  & \multicolumn{2}{c}{Two-grid scheme} & \multicolumn{2}{c}{Implicit-explicit scheme}\\
			\midrule
			$ (N, N^{h}_{x}) $ &  Error  &  Order  &  Error  &  Order  &  Error  &  Order \\
			\midrule 
			$ (80, 100) $  &  $2.23 \times 10^{-3}$    &   ---      &  $2.23 \times 10^{-3}$    &  ---    &  Inf  & ---       \\ 
			$ (180, 150) $  &  $4.81 \times 10^{-4}$   &   3.78      &  $4.81 \times 10^{-4}$    &  3.78   &  Inf  & ---      \\ 
			$ (320, 200) $  &  $1.55 \times 10^{-4}$    &   3.94      &  $1.55 \times 10^{-4}$    &  3.94   &  Inf   & ---       \\ 
			$ (500, 250) $  &  $6.40 \times 10^{-5}$    &   3.97      &  $6.40 \times 10^{-5}$    &  3.97   &  $6.77 \times 10^{-5}$  & ---      \\ 
			\bottomrule
	\end{tabular*}}
\end{table}
\begin{table}[!tbp]
	\caption{Spatial convergence of nonlinear scheme \eqref{sch:Non_CDM}, two-grid scheme \eqref{sch:coar_equ1}--\eqref{sch:fine_equ1} and implicit-explicit scheme \eqref{sch:IM} for Case III on uniform temporal grids\label{Ex:tabS3}}%
	{\footnotesize\begin{tabular*}{\columnwidth}{@{\extracolsep\fill}ccccccc@{\extracolsep\fill}}
			\toprule
			& \multicolumn{2}{c}{Nonlinear scheme}  & \multicolumn{2}{c}{Two-grid scheme} & \multicolumn{2}{c}{Implicit-explicit scheme}\\
			\midrule
			$ (N, N^{h}_{x}) $ &  Error  &  Order  &  Error  &  Order  &  Error  &  Order \\
			\midrule 
			$ (80, 100) $  &  $1.84 \times 10^{-2}$    &   ---      &  $1.84 \times 10^{-2}$    &  ---    &  Inf  & ---       \\ 
			$ (180, 150) $  &  $4.03 \times 10^{-3}$   &   3.75      &  $4.03 \times 10^{-3}$    &  3.75   &  Inf  & ---      \\ 
			$ (320, 200) $  &  $1.30 \times 10^{-3}$    &   3.92      &  $1.30 \times 10^{-3}$    &  3.92   &  Inf   & ---       \\ 
			$ (500, 250) $  &  $5.38 \times 10^{-4}$    &   3.96      &  $5.38 \times 10^{-4}$    &  3.96   &  Inf  & ---      \\ 
			\bottomrule
	\end{tabular*}}
\end{table}

Secondly, to test the temporal convergence rates, we fix $ N^{h}_{x} = N^{h}_{y} = 10 N^{H}_{x} = 10 N^{H}_{y} = 300 $ and present the numerical results with respect to $N$ in Table \ref{Ex:tab2} for Case I. We can observe that these three methods all have second-order temporal accuracy as proved. However, when the solution changes dramatically over time (e.g. Case II or III), the implicit-explicit method becomes unstable ($ |u| \rightarrow \infty $), while both the nonlinear method and two-grid method can generate the desired numerical solutions with the same magnitude accuracy, as seen in Tables \ref{Ex:tab3}--\ref{Ex:tab4}. Furthermore, the numerical results show that if we further refine the temporal grids (e.g., changing $N$ from $256$ to $512$ in Table \ref{Ex:tab3}), the implicit-explicit method may produce a correct result. But as shown in Table \ref{Ex:tab4}, its stability requirement for the temporal grid is quite related to the smoothness of the solution with respect to time. It is seen that the implicit-explicit discretization has much more strict stability condition compared to the other two methods when the solution $ u $ changes sharply with respect to $t$, for which a convincing explanation is that approximating the nonlinear term via solutions at previous time levels may leads to inaccuracy in this context. 
\begin{table}[!hb]
	\caption{Temporal convergence of nonlinear scheme \eqref{sch:Non_CDM}, two-grid scheme \eqref{sch:coar_equ1}--\eqref{sch:fine_equ1} and implicit-explicit scheme \eqref{sch:IM} for Case I on uniform temporal grids\label{Ex:tab2}}%
	{\footnotesize\begin{tabular*}{\columnwidth}{@{\extracolsep\fill}ccccccc@{\extracolsep\fill}}
			\toprule
			& \multicolumn{2}{c}{Nonlinear scheme}  & \multicolumn{2}{c}{Two-grid scheme} & \multicolumn{2}{c}{Implicit-explicit scheme}\\
			\midrule
			$ N $ &  Error  &  Order  &  Error  &  Order  &  Error  &  Order \\
			\midrule 
			$ 32 $  &  $4.65 \times 10^{-3}$    &   ---      &  $4.65 \times 10^{-3}$    &  ---    &  $5.25 \times 10^{-3}$  & ---       \\ 
			$ 64 $  &  $1.26 \times 10^{-3}$   &   1.88      &  $1.26 \times 10^{-3}$    &  1.88   &  $1.48 \times 10^{-3}$  & 1.82      \\ 
			$ 128 $  &  $3.23 \times 10^{-4}$    &   1.97      &  $3.23 \times 10^{-4}$    &  1.97   &  $3.56 \times 10^{-4}$   & 2.06       \\ 
			$ 256 $  &  $8.13 \times 10^{-5}$    &   1.99      &  $8.13 \times 10^{-5}$    &  1.99   &  $8.63 \times 10^{-5}$  & 2.05      \\ 
			\bottomrule
	\end{tabular*}}
\end{table}
\begin{table}[!htbp]
	\caption{Temporal convergence of nonlinear scheme \eqref{sch:Non_CDM}, two-grid scheme \eqref{sch:coar_equ1}--\eqref{sch:fine_equ1} and implicit-explicit scheme \eqref{sch:IM} for Case II on uniform temporal grids\label{Ex:tab3}}%
	{\footnotesize\begin{tabular*}{\columnwidth}{@{\extracolsep\fill}ccccccc@{\extracolsep\fill}}
			\toprule
			& \multicolumn{2}{c}{Nonlinear scheme}  & \multicolumn{2}{c}{Two-grid scheme} & \multicolumn{2}{c}{Implicit-explicit scheme}\\
			\midrule
			$ N $ &  Error  &  Order  &  Error  &  Order  &  Error  &  Order \\
			\midrule 
			$ 128 $  &  $9.30 \times 10^{-4}$    &   ---      &  $9.30 \times 10^{-4}$    &  ---    &  Inf  & ---       \\ 
			$ 256 $  &  $2.41 \times 10^{-4}$    &   1.95      &  $2.41 \times 10^{-4}$    &  1.95    &  Inf  & ---       \\ 
			$ 512 $  &  $6.10 \times 10^{-5}$   &   1.98      &  $6.10 \times 10^{-5}$    &  1.98   &  $6.45 \times 10^{-5}$  & ---      \\ 
			$ 1024 $  &  $1.53 \times 10^{-5}$    &   1.99      &  $1.53 \times 10^{-5}$    &  1.99   &  $1.61 \times 10^{-5}$   & 2.01       \\ 
			\bottomrule
	\end{tabular*}}
\end{table}
\begin{table}[!htbp]
	\caption{Temporal convergence of nonlinear scheme \eqref{sch:Non_CDM}, two-grid scheme \eqref{sch:coar_equ1}--\eqref{sch:fine_equ1} and implicit-explicit scheme \eqref{sch:IM} for Case III on uniform temporal grids\label{Ex:tab4}}%
	{\footnotesize\begin{tabular*}{\columnwidth}{@{\extracolsep\fill}ccccccc@{\extracolsep\fill}}
			\toprule
			& \multicolumn{2}{c}{Nonlinear scheme}  & \multicolumn{2}{c}{Two-grid scheme} & \multicolumn{2}{c}{Implicit-explicit scheme}\\
			\midrule
			$ N $ &  Error  &  Order  &  Error  &  Order  &  Error  &  Order \\
			\midrule
			$ 128 $  &  $7.70 \times 10^{-3}$    &   ---      &  $7.70 \times 10^{-3}$    &  ---    &  Inf  & ---       \\  
			$ 256 $  &  $2.00 \times 10^{-3}$    &   1.94      &  $2.00 \times 10^{-3}$    &  1.94    &  Inf  & ---       \\ 
			$ 512 $  &  $5.13 \times 10^{-4}$   &   1.96      &  $5.13 \times 10^{-4}$    &  1.96   &  Inf  & ---      \\ 
			$ 1024 $  &  $1.29 \times 10^{-4}$    &   1.99      &  $1.29 \times 10^{-4}$    &  1.99   &  Inf   & ---       \\ 
			\bottomrule
	\end{tabular*}}
\end{table}

\subsection{Accuracy and efficiency tests on variable-step temporal grids}
To check  the accuracy and efficiency on variable-step temporal grids, we consider model \eqref{Model:Bihar}--\eqref{Bound:Bihar} on $ (0,1)^2 \times (0,1] $ with $ c = \frac{1}{8\pi^2} $ and $ f(u) = u - u^3 $. The linear part $g$ is determined such that the exact solution is $ u(x,y,t) = \sin(t) \sin(2 \pi x) \sin(2 \pi y) $. The variable-step temporal grids are generated randomly by 
$$
\tau_{k} := T \frac{\theta_{k}}{S}, \quad {\rm with } \  S = \sum^{N}_{k=1} \theta_{k}, 
$$
where $ \theta_{k} $ is randomly drawn from the uniform distribution on the interval $ ( 1/4.8645, 1 ) $ such that the adjacent temporal stepsize ratio $ r_{k} < 4.8645 $.

As the implicit-explicit scheme \eqref{sch:IM} may generate wrong results, we only test the nonlinear scheme \eqref{sch:Non_CDM} and two-grid scheme \eqref{sch:coar_equ1}--\eqref{sch:fine_equ1}. We firstly test the errors and convergence rates in spatial and temporal directions for both methods with $ N^{h}_{x} = N^{h}_{y} $ and $ M_{x} = M_{y} = 5 $. The corresponding numerical results are listed in Tables \ref{tab1}--\ref{tab2} respectively, which indicates the fourth-order accuracy in space and second-order accuracy in time as proved in Theorems \ref{thm:Conver_Non} and \ref{thm:Conver_Fine}. Moreover, we compare the CPU times consumed by the two methods in Table \ref{tab3} under $ N^{h}_{x} = N^{h}_{y} $ and $ M_{x} = M_{y} $, in which we can clearly observe that the proposed two-grid method has significantly improved the computational efficiency, for example, it takes about two and a half hours for the implementation of the nonlinear scheme when $ N = N^{h}_{x} = 480 $, while the two-grid scheme consumes only about one hour to desire the same error!
\begin{table}[!htbp]
	\caption{Spatial convergence of nonlinear scheme \eqref{sch:Non_CDM} and two-grid scheme \eqref{sch:coar_equ1}--\eqref{sch:fine_equ1} on variable-step temporal grids \label{tab1}}%
	{\footnotesize\begin{tabular*}{\columnwidth}{@{\extracolsep\fill}cccccc@{\extracolsep\fill}}
			\toprule
			&  & \multicolumn{2}{c}{Nonlinear scheme} & \multicolumn{2}{c}{Two-grid scheme}\\
			\midrule
			$ (N, N^{h}_{x}) $ & $ \max r_{k} $  & Error & Order  & Error & Order\\
			\midrule
			$ (80, 100) $   &   4.3125     & $1.76 \times 10^{-5}$    & ---      & $1.76 \times 10^{-5}$       &   ---\\ 
			$ (180, 150) $  &   4.4777     & $3.67 \times 10^{-6}$    & 3.87     & $3.67 \times 10^{-6}$       &   3.87\\
			$ (320, 200) $ &   4.5252     & $1.14 \times 10^{-6}$    & 4.06     & $1.14 \times 10^{-6}$       &   4.06\\
			$ (500, 250) $ &   4.5481     & $4.65 \times 10^{-7}$    & 4.02   & $4.65 \times 10^{-7}$     &   4.02\\
			\bottomrule
	\end{tabular*}}
\end{table}
\begin{table}[!htbp]
	\caption{Temporal convergence of nonlinear scheme \eqref{sch:Non_CDM} and two-grid scheme \eqref{sch:coar_equ1}--\eqref{sch:fine_equ1} on variable-step temporal grids \label{tab2}}%
	{\footnotesize\begin{tabular*}{\columnwidth}{@{\extracolsep\fill}cccccc@{\extracolsep\fill}}
			\toprule
			 & & \multicolumn{2}{c}{Nonlinear scheme} & \multicolumn{2}{c}{Two-grid scheme}\\
			\midrule
			$ (N, N^{h}_{x}) $ & $ \max r_{k} $  & Error & Order & Error & Order\\
			\midrule
			$ (40, 40) $   &   3.9128     & $7.36 \times 10^{-5}$    & ---      & $6.17 \times 10^{-5}$       &   ---\\ 
			$ (80, 80) $   &   3.9673     & $1.80 \times 10^{-5}$    & 2.03     & $1.82 \times 10^{-5}$       &   1.76\\
			$ (160, 160) $ &   4.4221     & $4.66 \times 10^{-6}$    & 1.95     & $4.53 \times 10^{-6}$       &   2.01\\
			$ (320, 320) $ &   4.5826     & $1.16 \times 10^{-6}$    & 2.01     & $1.18 \times 10^{-6}$       &   1.94\\
			\bottomrule
	\end{tabular*}}
\end{table}

\begin{table}[!htbp]
	\caption{Errors and CPU times of nonlinear scheme \eqref{sch:Non_CDM} and two-grid scheme \eqref{sch:coar_equ1}--\eqref{sch:fine_equ1} on variable-step temporal grids \label{tab3}}%
	{\footnotesize\begin{tabular*}{\columnwidth}{@{\extracolsep\fill}cccccc@{\extracolsep\fill}}
			\toprule
			&  &  \multicolumn{2}{c}{Nonlinear scheme} & \multicolumn{2}{c}{Two-grid scheme}\\
			\midrule
			$ (N, N^{h}_{x}) $ & $ N^{H}_{x} $  & Error  & CPU times  & Error  & CPU times  \\
			\midrule
			(240, 240)    &   60  &  $2.05 \times 10^{-6}$      & 33 m 42 s    & $2.05 \times 10^{-6}$   & 7 m 17 s       \\
			(320, 320)    &   80  &  $1.15 \times 10^{-6}$   & 1 h 7 m 20 s   & $1.15 \times 10^{-6}$   & 22 m 35 s       \\
			(400, 400)    &  100  &   $7.47 \times 10^{-7}$  & 1 h 52 m 28 s   & $7.47 \times 10^{-7}$   & 41 m 27 s       \\ 
			(480, 480)    &  120  &   $5.19 \times 10^{-7}$  & 2 h 27 m 7 s   & $5.19 \times 10^{-7}$   & 1 h 2 m 18 s       \\
			\bottomrule
	\end{tabular*}}
\end{table}

\subsection{Effectiveness of adaptive temporal stepsize strategy}
In this test, we consider model \eqref{Model:Bihar}--\eqref{Bound:Bihar} on $ (0,1)^2 \times (0,4] $ with $ c = 1 $ and $ f(u) = \sin u $, and the exact solution is chosen as $ u(x,y,t) = \big[ 1 + 20 e^{-40 (t-1)^2} + 30 e^{-60 (t-4)^2}\big] \sin(2 \pi x) \sin(2 \pi y) $. Figure \ref{fig1} (left) depicts the evolution of solution with respect to time at fixed point $(0.25,0.25)$, which consists of two peaks and admits multiple time scales. Therefore, the variable-step two-grid scheme based on the adaptive temporal stepsize strategy \cite{SISC_Huang_2020,JSC_Liao_2022} will be adopted to improve the temporal accuracy 
\begin{equation}\label{Algor:Adaptive}
	\begin{aligned}
		\tau_{n+1} = \min\left\{ \max\left\{ \tau_{\min} , \frac{\tau_{\max}}{ \sqrt{ 1 + \eta \| \partial_{\tau} u^{n} \|_{h}^{2} } } \right\}, r_{\max} \tau_{n} \right\},
	\end{aligned}
\end{equation}
where $ \partial_{\tau} u^{n} := \nabla_{\tau} u^{n}/\tau_{n} $ and $ r_{\max} = 4.8 $ which satisfies the restrictions in Theorems \ref{thm:Conver_Non} and \ref{thm:Conver_Fine}. Here $ \tau_{\min} $ and $ \tau_{\max} $ are the pre-determined minimum and maximum temporal stepsize and $ \eta $ is a pre-chosen parameter. In this example, we uniformly set $ \tau_{\max} = 0.2 $, $ \eta = 500 $ and gradually reduce $ \tau_{\min} $ to generate grids with distinct stepsizes. 
\begin{figure}[!htbp]
	\centering
	\subfigure{
		\begin{minipage}[t]{0.5\linewidth}
			\centering
			\includegraphics[width=2.75in]{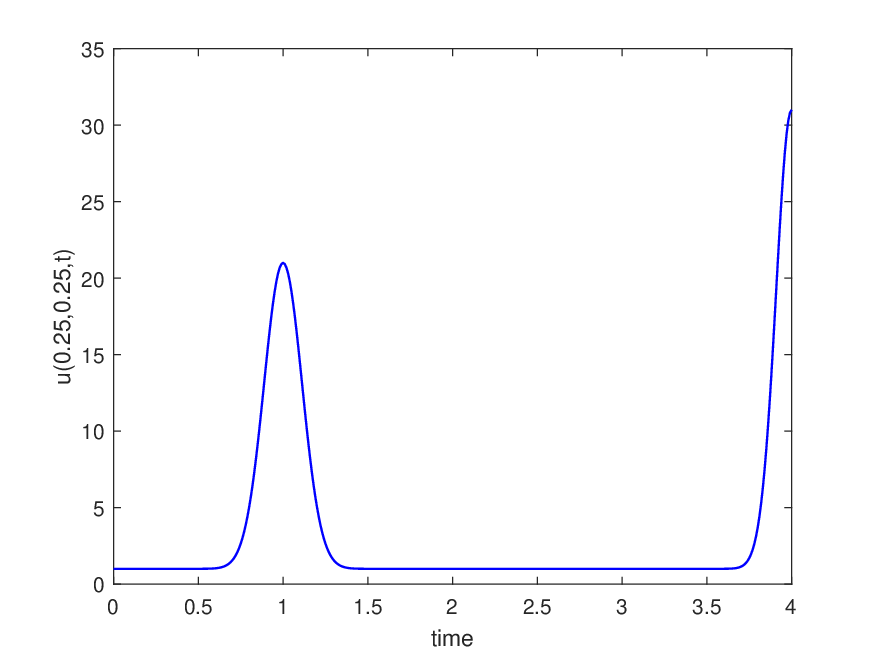}
		\end{minipage}%
	}%
	\subfigure{
		\begin{minipage}[t]{0.5\linewidth}
			\centering
			\includegraphics[width=2.75in]{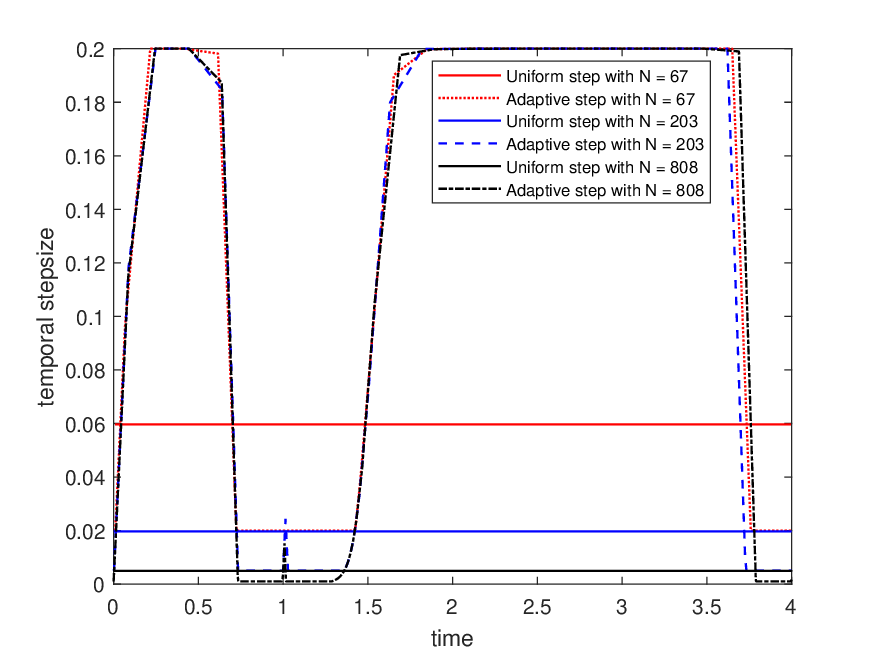}
		\end{minipage}%
	}%
	\centering
	\caption{Evolution of solution with time (left) at $(0.25,0.25)$ and temporal stepsize (right) of the two-grid scheme  and adaptive two-grid scheme}
	\label{fig1}
\end{figure}

\begin{table}[!htbp]
	\caption{Temporal convergence of two-grid scheme and adaptive two-grid scheme \label{tab4}}%
	{\footnotesize\begin{tabular*}{\columnwidth}{@{\extracolsep\fill}cccccccc@{\extracolsep\fill}}
			\toprule
			&   & \multicolumn{3}{c}{Two-grid scheme} & \multicolumn{3}{c}{Adaptive two-grid scheme}\\
			\midrule
			$ \tau_{\min} $ & N  & CPU times & Error  & Order & CPU times & Error & Order \\
			\midrule
			0.02    &   67   &  27.96 s   & $3.06 \times 10^{-1}$    & ---   &  28.73 s  & $2.18 \times 10^{-2}$   & ---     \\ 
			0.01    &   114  &  45.63 s   & $9.93 \times 10^{-2}$    & 2.12  &  48.43 s & $4.63 \times 10^{-3}$   & 2.91    \\
			0.005   &   203  & 1 m 21 s   & $2.53 \times 10^{-2}$    & 2.37  &   1 m 22 s  & $1.01 \times 10^{-3}$   & 2.64     \\
			0.002   &   489  & 3 m 13 s   & $3.12 \times 10^{-3}$    & 2.38   &   3 m 19 s  & $1.41 \times 10^{-4}$   & 2.24   \\
			0.001   &   808  & 5 m 30 s   & $1.00 \times 10^{-3}$    & 2.27  &   5 m 28 s  & $3.46 \times 10^{-5}$   & 2.80    \\
			\bottomrule
	\end{tabular*}}
\end{table}

We select $ N^{h}_{x} = N^{h}_{y} = 250 $ to test the temporal convergence rates yielded by the two-grid scheme on uniform and adaptive temporal grids with $ M_{x} = M_{y} = 10 $. The numerical results are presented in Table \ref{tab4}, which shows that the standard two-grid scheme is second-order accurate in time, while the variable-step two-grid scheme based on the adaptive temporal stepsize strategy \eqref{Algor:Adaptive} has better convergence and much smaller errors with the same number of grids. In other words, the standard two-grid scheme on uniform temporal grids requires more temporal steps, and of course more CPU times, to generate numerical solutions with the same magnitude accuracy as adaptive method which uses less temporal steps. The reason can be more intuitively observed from Figure \ref{fig1} (right), which shows that when the solution varies sharply, small temporal stepsizes are adaptively created to capture the fast evolution process, while otherwise large temporal stepsizes are generated to accelerate the time integration.

\subsection{ Application to phase-field Allen--Cahn equation }
In this subsection, we consider the following Allen--Cahn equation with a polynomial double-well potential, subject to periodic boundary condition
\begin{equation*}
	\begin{aligned}
		u_{t} - \varepsilon^2 \Delta u = u - u^3, \quad (x,y) \in \Omega, \  t \in (0,T],
	\end{aligned}
\end{equation*}
where $ \varepsilon $ is the interaction length that describes the thickness of the transition boundary between materials. It is well known that the energy dissipation law \cite{2019_SIAM_Shen,2023_JSC_Qiao}
\begin{equation*}
	\begin{aligned}
		E[u](t) \leq E[u](s), \quad \forall t > s
	\end{aligned}
\end{equation*}
holds for the Allen--Cahn equation, where $ E[u](t) $ represents the Lyapunov energy functional, namely
\begin{equation*}
	\begin{aligned}
		E[u](t) = \int_{\Omega} \frac{\varepsilon^{2}}{2} | \nabla u |^{2} + F(u) d\mathbf{x} \  \  {\rm with} \   \  F(u) = \frac{1}{4} ( 1 - u^2 )^2.
	\end{aligned}
\end{equation*}

Moreover, it has been observed that the evolution of the energy $ E[u](t) $ usually involves both fast and slow stages of change in the long time simulation. Thus, it is highly desirable for numerical methods to preserve the discrete energy dissipation law on the nonuniform temporal grids. Define the discrete energy functional $ \mathcal{E}[u^{n}] $ as
\begin{equation*}
	\begin{aligned}
			\mathcal{E}[u^{n}] := - \frac{\varepsilon^{2}}{2} h_x h_y\sum^{N^{h}_{x}}_{i=1} \sum^{N^{h}_{y}}_{j=1} u^{n}_{i,j}  \Delta_{h} u^{n}_{i,j} + \frac{1}{4} h_x h_y\sum^{N^{h}_{x}}_{i=1} \sum^{N^{h}_{y}}_{j=1} \left( 1 - ( u^{n}_{i,j} )^{2} \right)^{2}.
		\end{aligned}
\end{equation*}
In this test, we would also compare effectiveness and efficiency of the high-order nonlinear difference scheme \eqref{sch:Non_CDM}, the two-grid difference scheme \eqref{sch:coar_equ1}--\eqref{sch:fine_equ1} and the adaptive two-grid difference scheme using a similar adaptive temporal stepsize strategy \eqref{Algor:Adaptive} by replacing $ \partial_{\tau} u^{n} $ with $ \partial_{\tau} \mathcal{E}[u^{n}] $. 

\begin{example}
	In this example, we set $ \varepsilon = 0.02$ and apply these three methods to simulate the merging of four bubbles with an initial condition
	\begin{equation*}
		\begin{aligned}
			u_{0}(x,y) = & -\tanh \left( \left( \left(x-0.3\right)^2 + y^2 - 0.2^2 \right)/\varepsilon \right) \tanh \left( \left( \left(x+0.3\right)^2 + y^2 - 0.2^2 \right)/\varepsilon \right) \\
			& \times \tanh \left( \left( x^2 + \left(y-0.3\right)^2 - 0.2^2 \right)/\varepsilon \right) \tanh \left( \left( x^2 + \left(y+0.3\right)^2 - 0.2^2 \right)/\varepsilon \right).
		\end{aligned}
	\end{equation*}
\end{example}

In this simulation, a $ 384 \times 384 $ uniform mesh is taken to discretize the spatial domain $ \Omega = (-1,1)^2 $ and the ratio of coarse-fine grids are setted as $ M_{x} = M_{y} = 3 $. We start with the modeling of the solution by the nonlinear scheme and two-grid scheme with a constant temporal stepsize $ \tau = 0.1 $ until time $ T = 100 $. Parameters in the adaptive temporal stepsize strategy \eqref{Algor:Adaptive} are selected as $ \tau_{\min} = 0.1 $, $ \tau_{\max} = 1$ and $ \eta = 3200 $. In Figure \ref{fig2}, it displays a comparison on the evolution of solution snapshots  among the three methods, in which the gradually merging and shrinking process of the initial four-drops over time can be clearly observed. As can be seen in the figures, there seems no distinguishable differences among these methods. Next, we investigate the efficiency of the two-grid method on uniform grid and on adaptive grid for long time modeling. As seen in Figure \ref{fig3} (left), the evolution of the free energy with respect to time for these three methods coincide, which consists very well with the corresponding results in \cite{SINUM_Liao_2020}. Moreover, Table \ref{tab5} indicates that the adaptive two-grid scheme has significant  advantage in computational efficiency over the other two schemes. For example, it takes about 10 hours for the implementation of the nonlinear scheme up to $ T = 100 $, while the two-grid method using uniform temporal grid consumes about 5 hours. What is even more amazing is that the developed adaptive two-grid method using variable-step temporal grid takes only about one hour! In fact, the total number of adaptive temporal steps is only 156, while it takes 1000 steps for the uniform grid. Finally, the adaptive temporal stepsize curve of the adaptive two-grid method is plotted in Figure \ref{fig3} (right), which also demonstrates the superiority of the variable-step two-grid compact difference scheme.

\begin{figure}[htbp]
	\centering
	\subfigure{
		\begin{minipage}[t]{0.25\linewidth}
			\centering
			\includegraphics[width=1.5in]{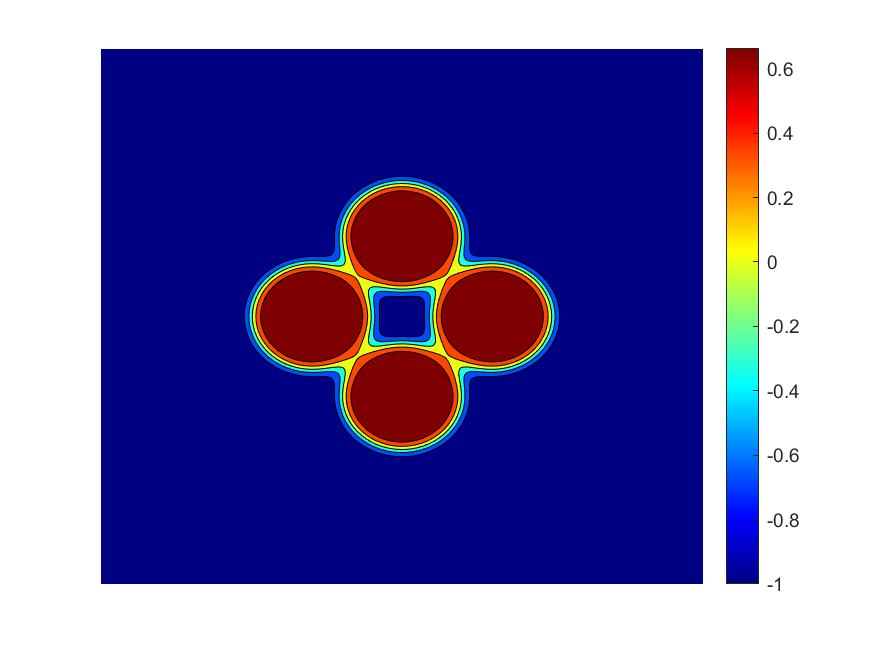}
		\end{minipage}%
	}%
	\subfigure{
		\begin{minipage}[t]{0.25\linewidth}
			\centering
			\includegraphics[width=1.5in]{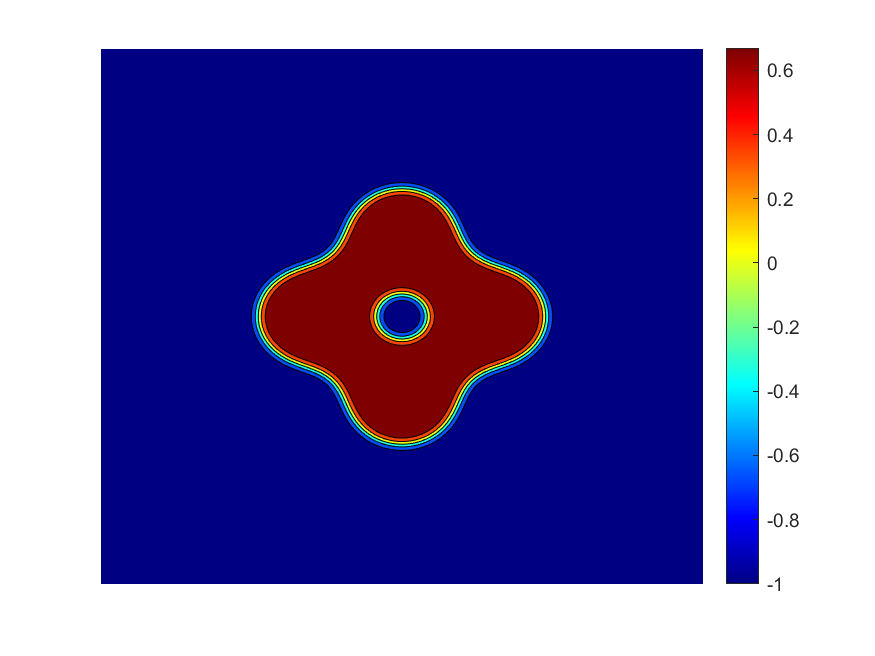}
		\end{minipage}%
	}%
	\subfigure{
		\begin{minipage}[t]{0.25\linewidth}
			\centering
			\includegraphics[width=1.5in]{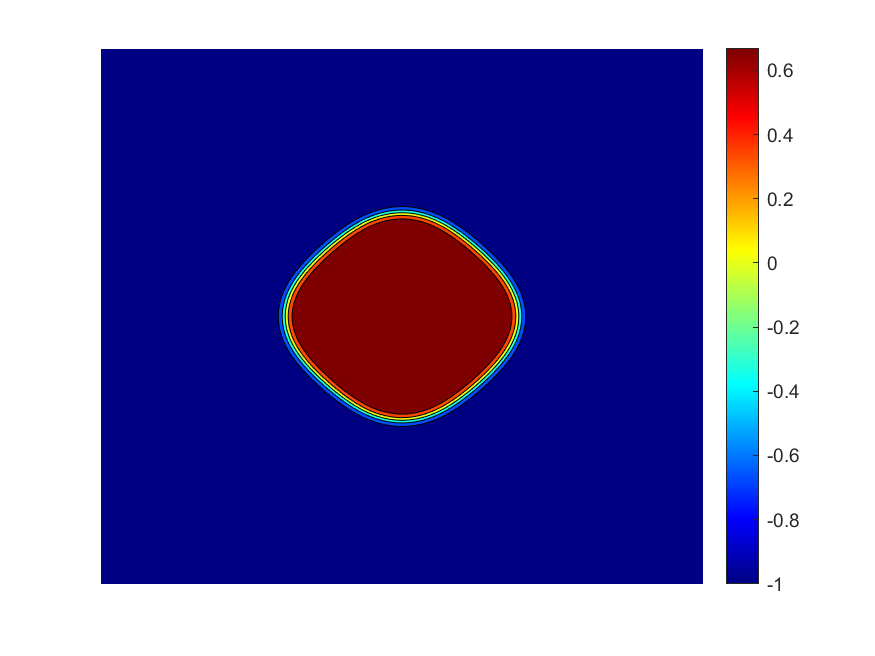}
		\end{minipage}%
	}%
	\subfigure{
		\begin{minipage}[t]{0.25\linewidth}
			\centering
			\includegraphics[width=1.5in]{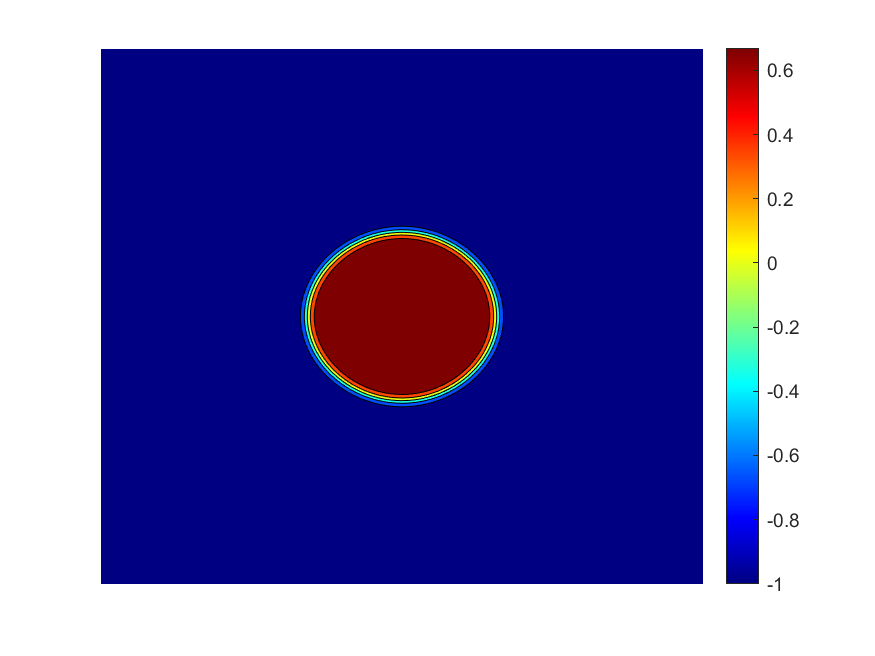}
		\end{minipage}%
	}%
	\vspace{-25pt}
	\setcounter{figure}{0}
	\renewcommand{\figurename}{}
	\renewcommand{\thefigure}{\Alph{figure}}
	\caption{Nonlinear scheme with fixed temporal stepsize $\tau = 0.1$}
	
	\subfigure{
		\begin{minipage}[t]{0.25\linewidth}
			\centering
			\includegraphics[width=1.5in]{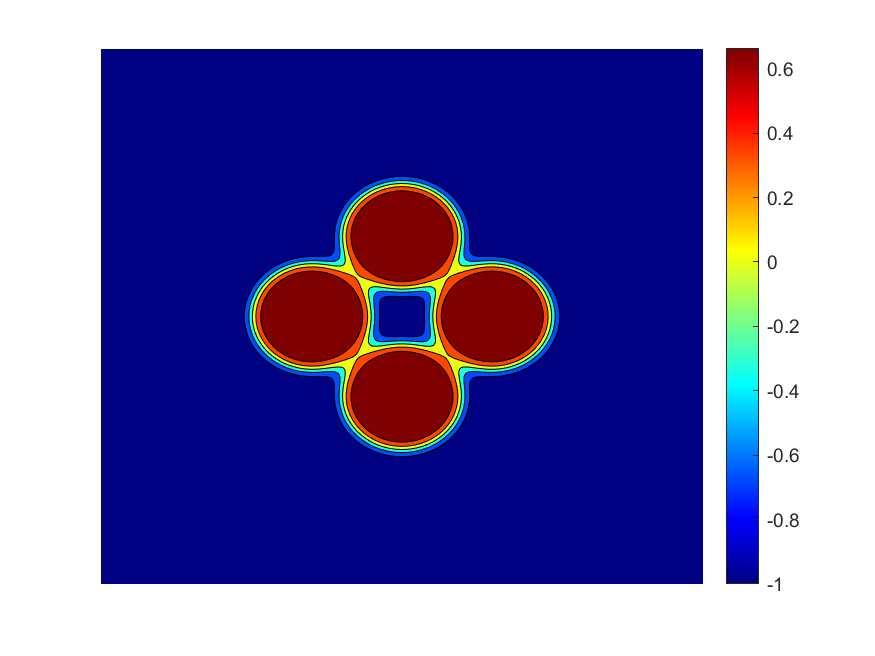}
		\end{minipage}%
	}%
	\subfigure{
		\begin{minipage}[t]{0.25\linewidth}
			\centering
			\includegraphics[width=1.5in]{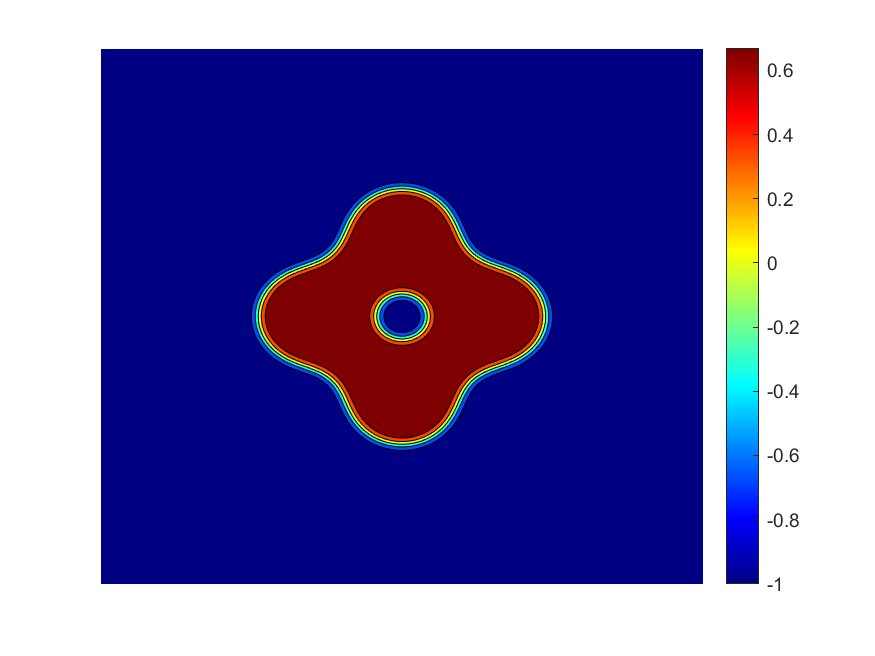}
		\end{minipage}%
	}%
	\subfigure{
		\begin{minipage}[t]{0.25\linewidth}
			\centering
			\includegraphics[width=1.5in]{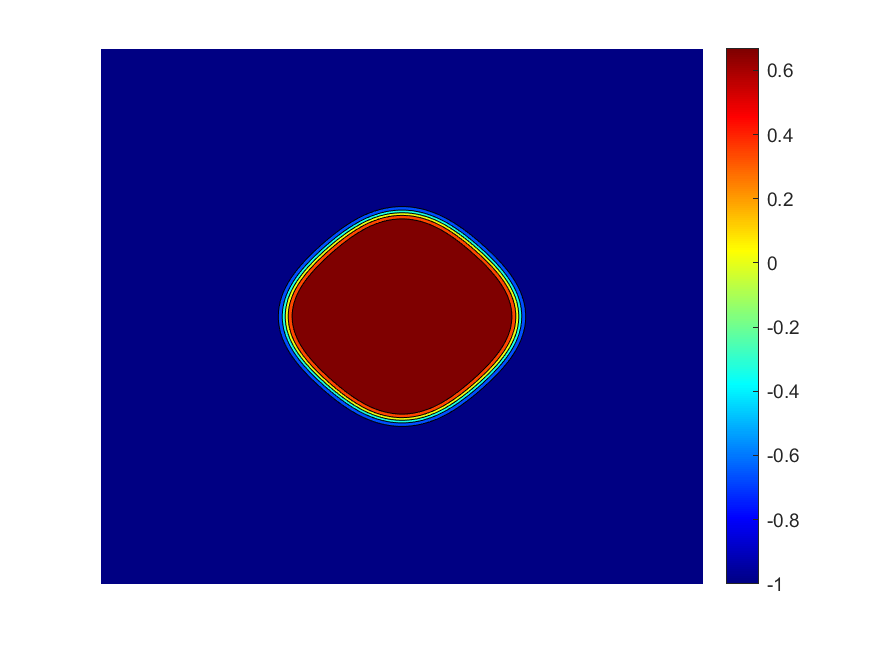}
		\end{minipage}%
	}%
	\subfigure{
		\begin{minipage}[t]{0.25\linewidth}
			\centering
			\includegraphics[width=1.5in]{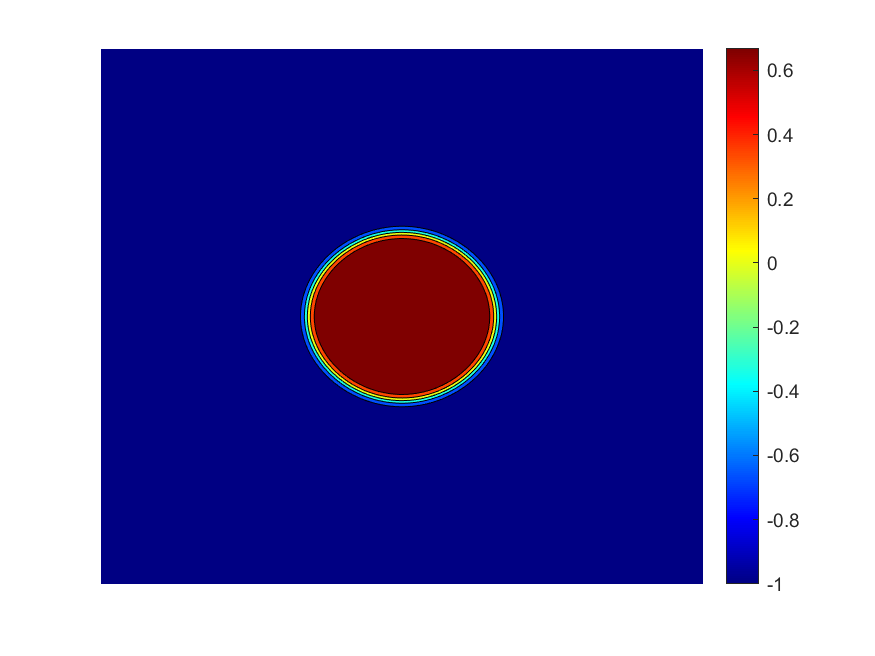}
		\end{minipage}%
	}%
	\vspace{-15pt}
	\caption{Two-grid scheme with fixed temporal stepsize $\tau = 0.1$}
	
	\subfigure{
		\begin{minipage}[t]{0.25\linewidth}
			\centering
			\includegraphics[width=1.5in]{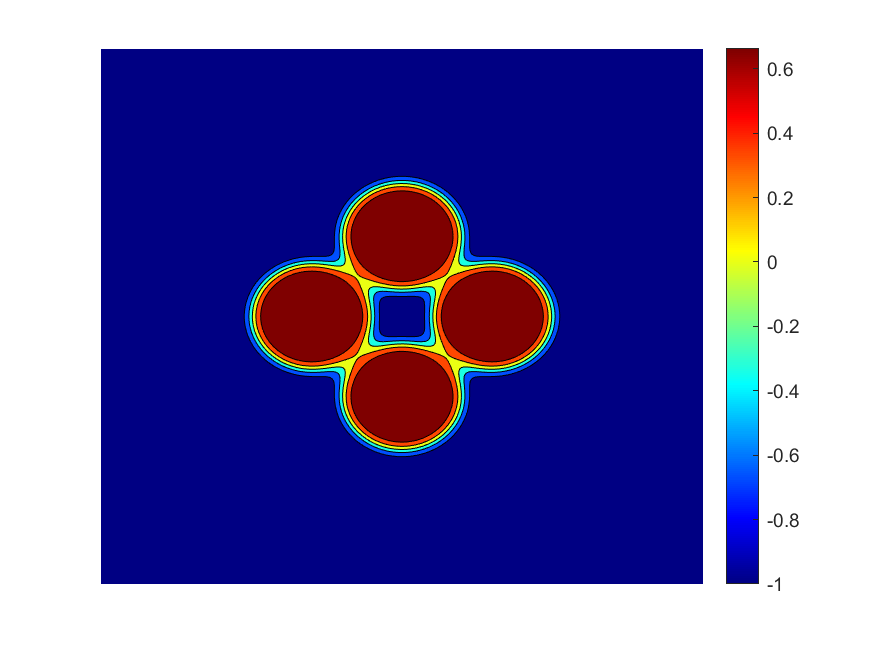}
		\end{minipage}%
	}%
	\subfigure{
		\begin{minipage}[t]{0.25\linewidth}
			\centering
			\includegraphics[width=1.5in]{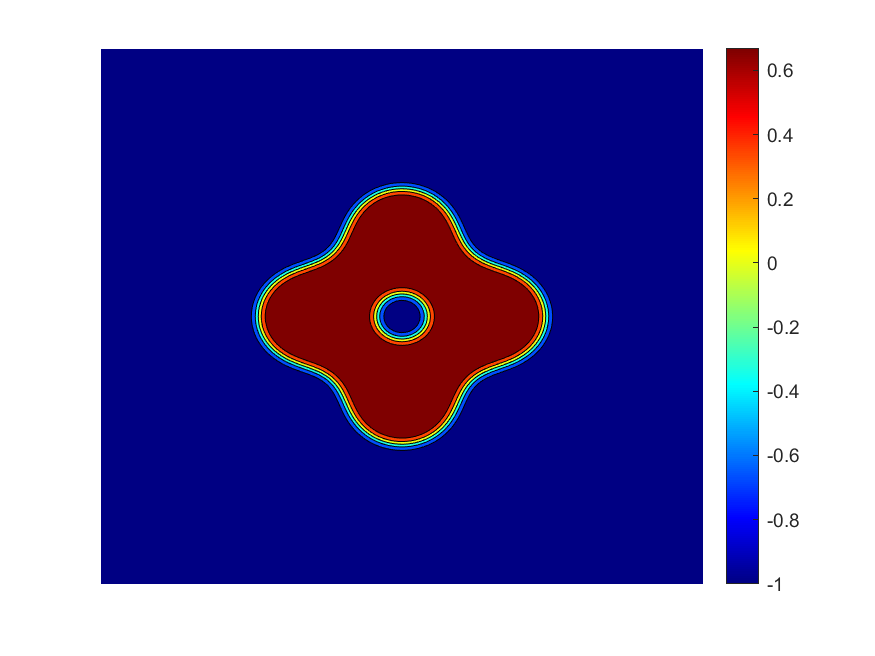}
		\end{minipage}%
	}%
	\subfigure{
		\begin{minipage}[t]{0.25\linewidth}
			\centering
			\includegraphics[width=1.5in]{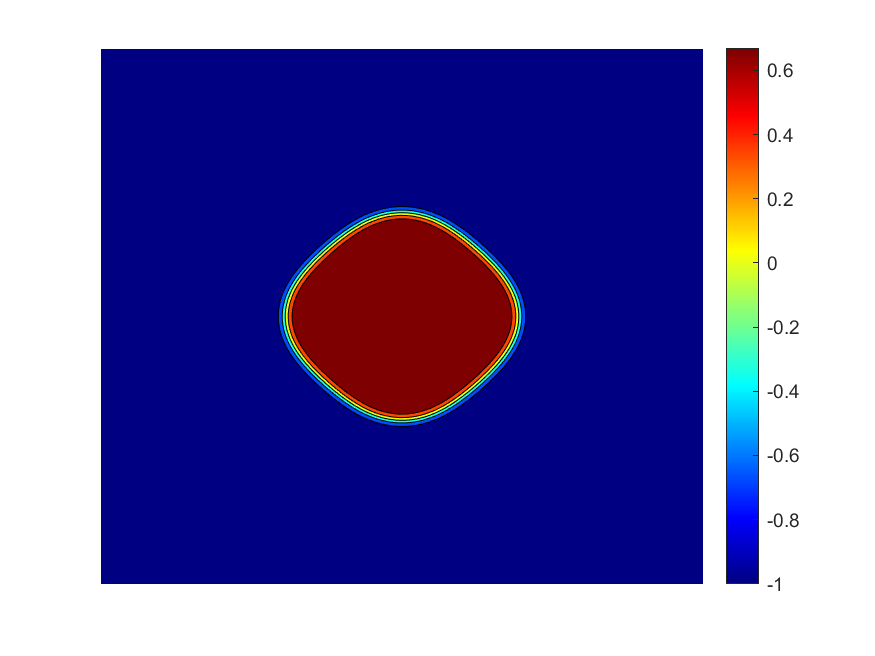}
		\end{minipage}%
	}%
	\subfigure{
		\begin{minipage}[t]{0.25\linewidth}
			\centering
			\includegraphics[width=1.5in]{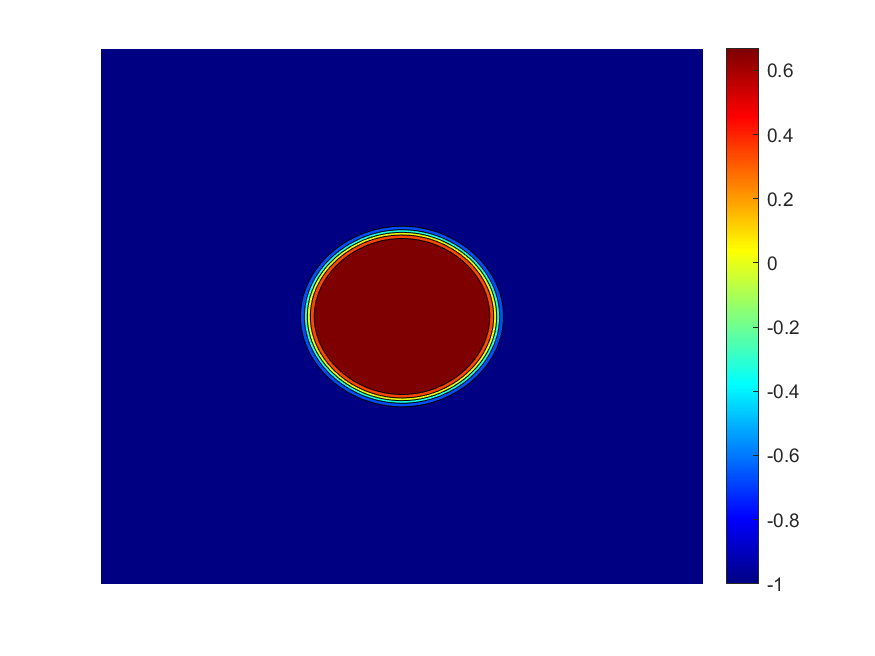}
		\end{minipage}%
	}%
	\vspace{-15pt}
	\caption{Adaptive two-grid scheme with $ \tau_{\min} = 0.1 $, $ \tau_{\max} = 1 $ and $ \eta = 3200 $}
	
	\centering
	\setcounter{figure}{2}
	\renewcommand{\figurename}{Figure}
	\renewcommand{\thefigure}{\arabic{figure}}
	\caption{ Solution snapshots of Allen--Cahn equation at $ t = 1, 10, 50, 100 $ (from left to right) yielded by nonlinear scheme, two-grid scheme and adaptive two-grid scheme }
	\label{fig2}
\end{figure}

\begin{figure}[htbp]
	\centering
	\subfigure{
		\begin{minipage}[t]{0.5\linewidth}
			\centering
			\includegraphics[width=2.75in]{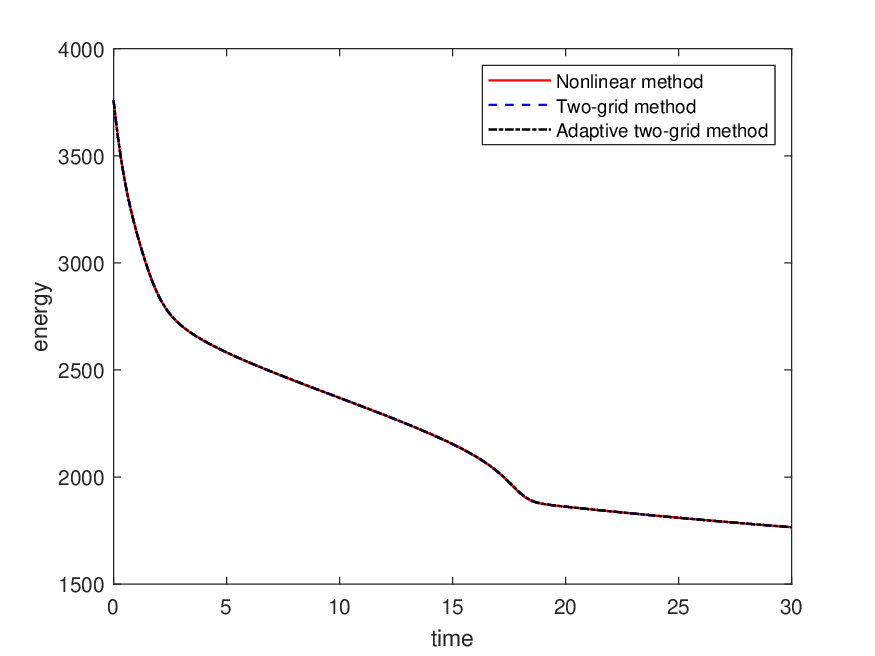}
		\end{minipage}%
	}%
	\subfigure{
		\begin{minipage}[t]{0.5\linewidth}
			\centering
			\includegraphics[width=2.75in]{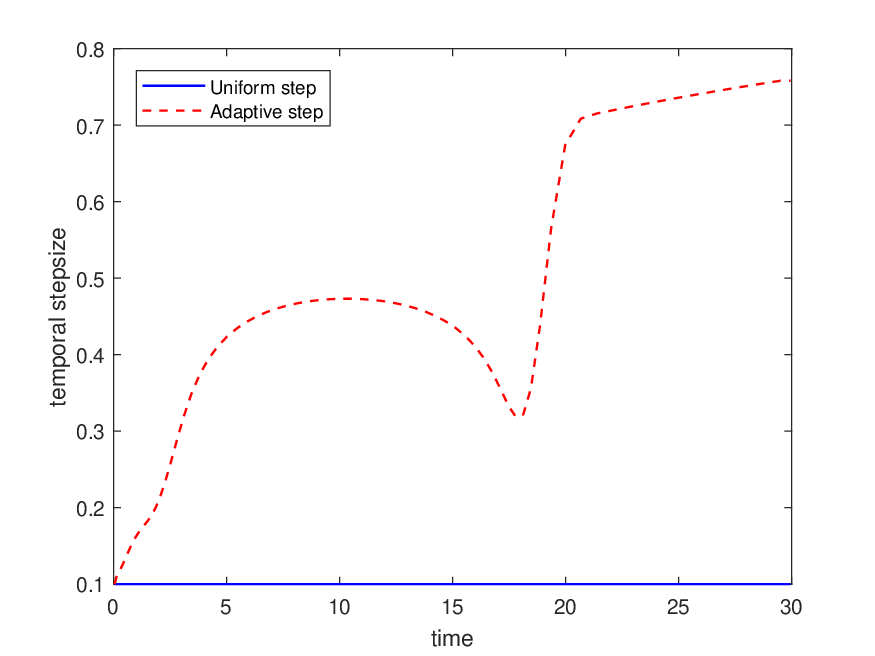}
		\end{minipage}%
	}%
	
	\centering
	\caption{ Evolutions of energy (left) and time steps (right) for the nonlinear scheme, two-grid scheme and adaptive two-grid scheme until time $ T = 30 $}
	\label{fig3}
\end{figure}

%

\begin{table}[!htbp]
	\caption{CPU times and the total number of temporal steps for three schemes\label{tab5}}%
	{\footnotesize\begin{tabular*}{\columnwidth}{@{\extracolsep\fill}ccccccc@{\extracolsep\fill}}
			\toprule
			& \multicolumn{2}{c}{Nonlinear scheme}  & \multicolumn{2}{c}{Two-grid scheme} & \multicolumn{2}{c}{Adaptive two-grid scheme}\\
			\midrule
			T & $ N $ & CPU times  & $ N $  & CPU times & $ N $ & CPU times \\
			\midrule 
			10  &  100    &   57 m 42 s      &  100    &  26 m 10 s    &  34  & 9 m 31 s       \\ 
			30  &  300    &   3 h 50 m 55 s      &  300    &  1 h 47 m 24 s    &  71  & 20 m 2 s       \\ 
			50  &  500    &   5 h 58 m 39 s      &  500    &  2 h 50 m 56 s    &  96   & 37 m 25 s       \\ 
			100  &  1000    &   9 h 43 m 51 s      &  1000    &  5 h 9 m 43 s    &  156  & 1 h 2 m 22 s      \\ 
			\bottomrule
	\end{tabular*}}
\end{table}

\begin{example}
In this example, we consider the coarsening process governed by the Allen--Cahn equation with the model parameter $ \varepsilon = 0.01 $ and computational domain $ \Omega = (0,1)^2 $. Here we choose a random initial condition $ u_{0}(x,y) = -0.05+0.1\times  rand(x,y) $.
\end{example}

\begin{figure}[htbp]
	\centering
	\subfigure{
		\begin{minipage}[t]{0.2\linewidth}
			\centering
			\includegraphics[width=1.1in]{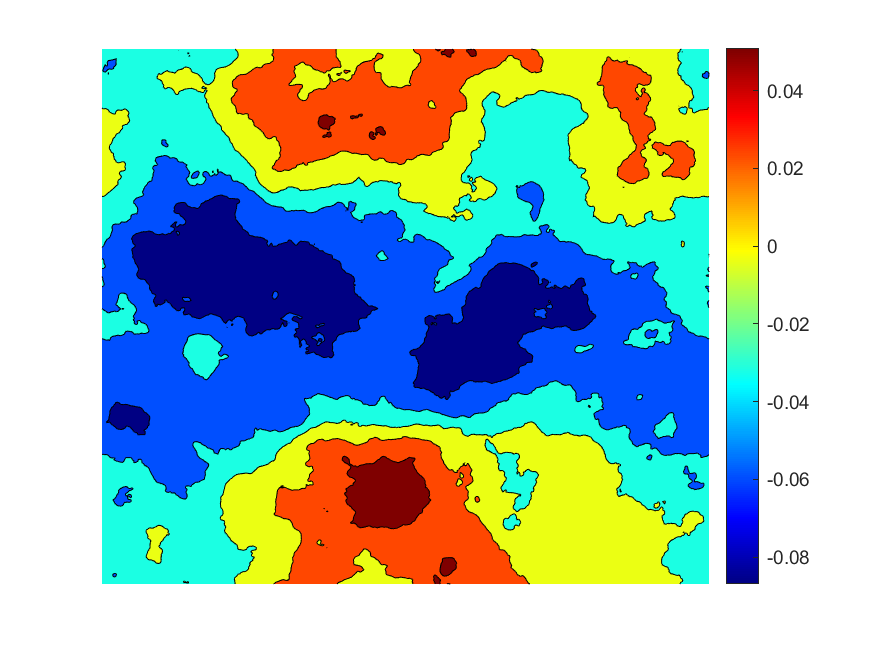}
		\end{minipage}%
	}%
	\subfigure{
		\begin{minipage}[t]{0.2\linewidth}
			\centering
			\includegraphics[width=1.1in]{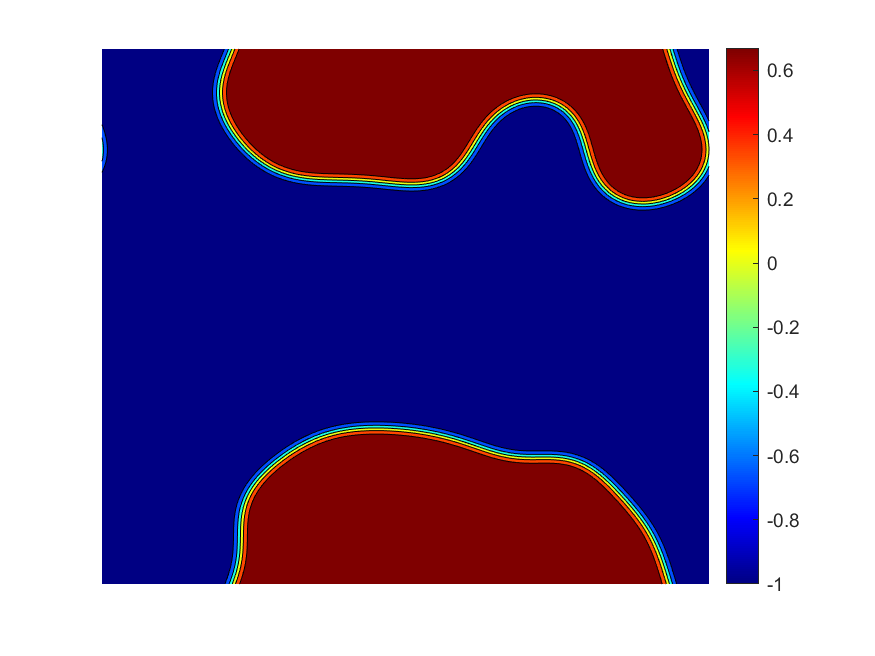}
		\end{minipage}%
	}%
	\subfigure{
		\begin{minipage}[t]{0.2\linewidth}
			\centering
			\includegraphics[width=1.1in]{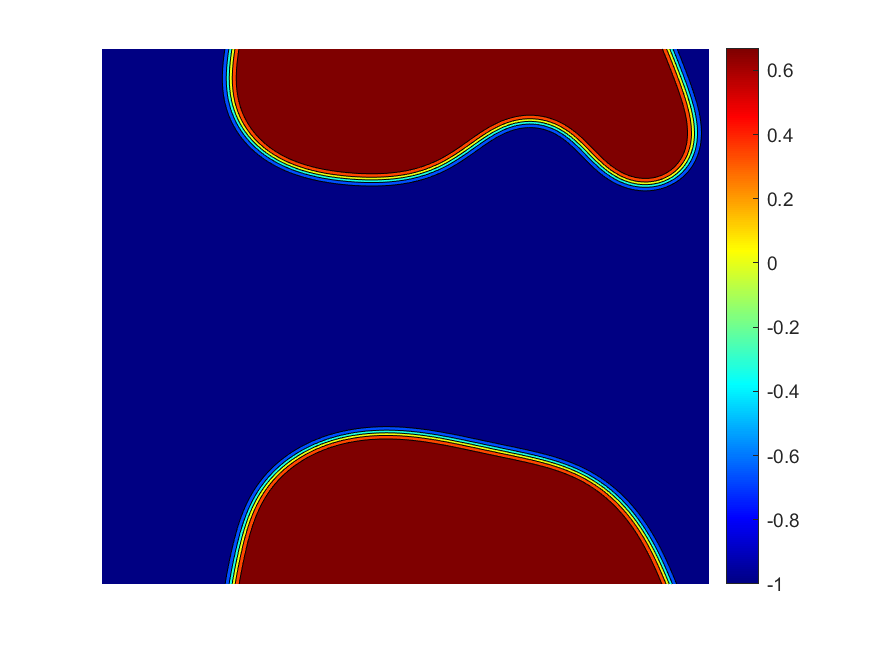}
		\end{minipage}%
	}%
	\subfigure{
		\begin{minipage}[t]{0.2\linewidth}
			\centering
			\includegraphics[width=1.1in]{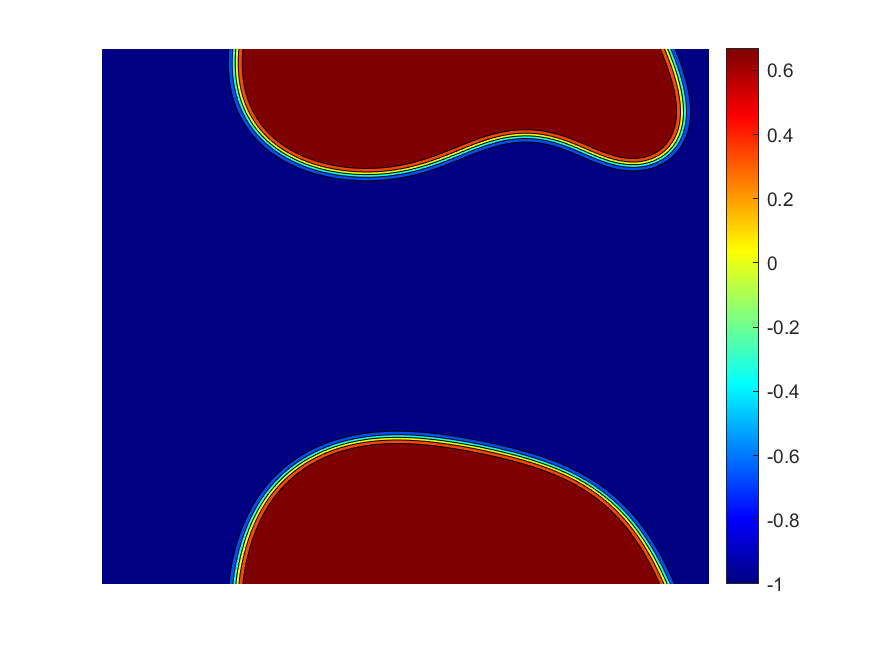}
		\end{minipage}%
	}%
	\subfigure{
		\begin{minipage}[t]{0.2\linewidth}
			\centering
			\includegraphics[width=1.1in]{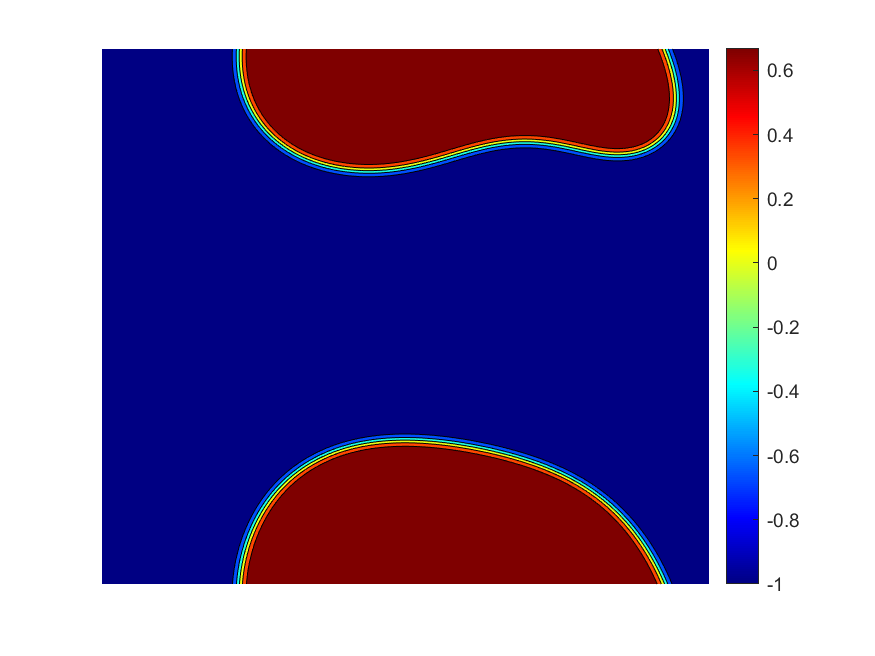}
		\end{minipage}%
	}%
	\vspace{-25pt}
	\setcounter{figure}{0}
	\renewcommand{\figurename}{}
	\renewcommand{\thefigure}{\Alph{figure}}
	\caption{Nonlinear scheme with fixed temporal stepsize $\tau = 1$}
	
	\subfigure{
		\begin{minipage}[t]{0.2\linewidth}
			\centering
			\includegraphics[width=1.1in]{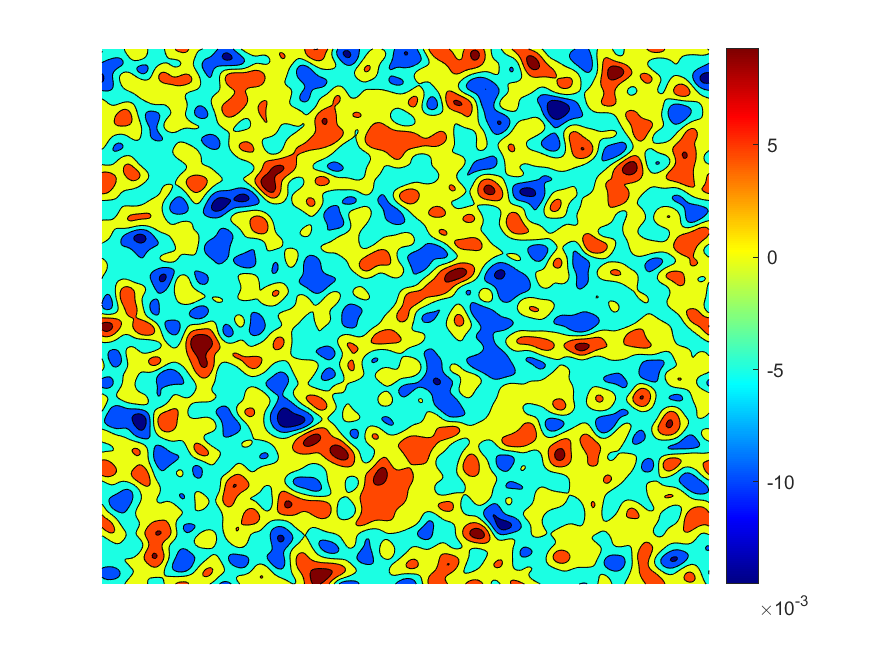}
		\end{minipage}%
	}%
	\subfigure{
		\begin{minipage}[t]{0.2\linewidth}
			\centering
			\includegraphics[width=1.1in]{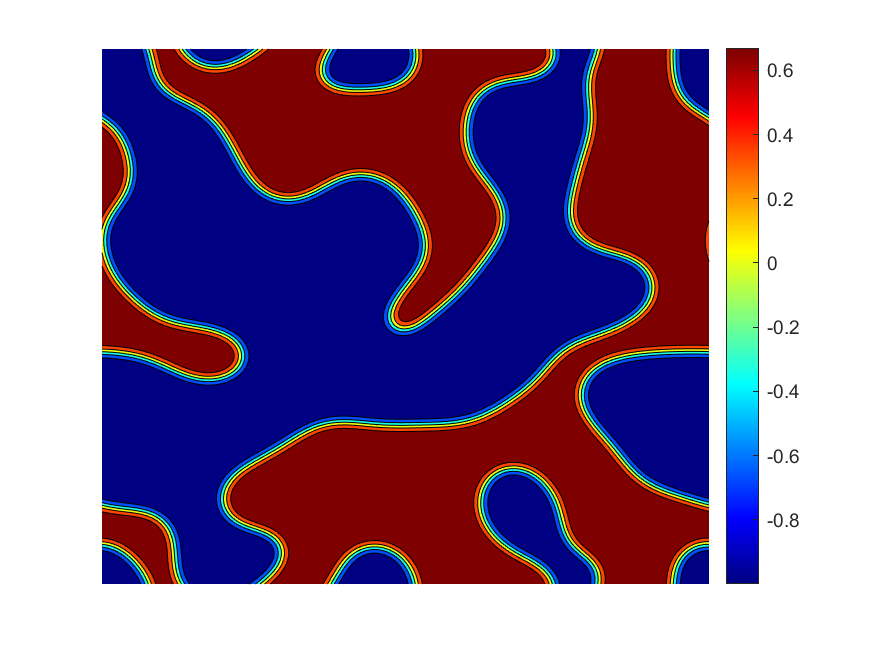}
		\end{minipage}%
	}%
	\subfigure{
		\begin{minipage}[t]{0.2\linewidth}
			\centering
			\includegraphics[width=1.1in]{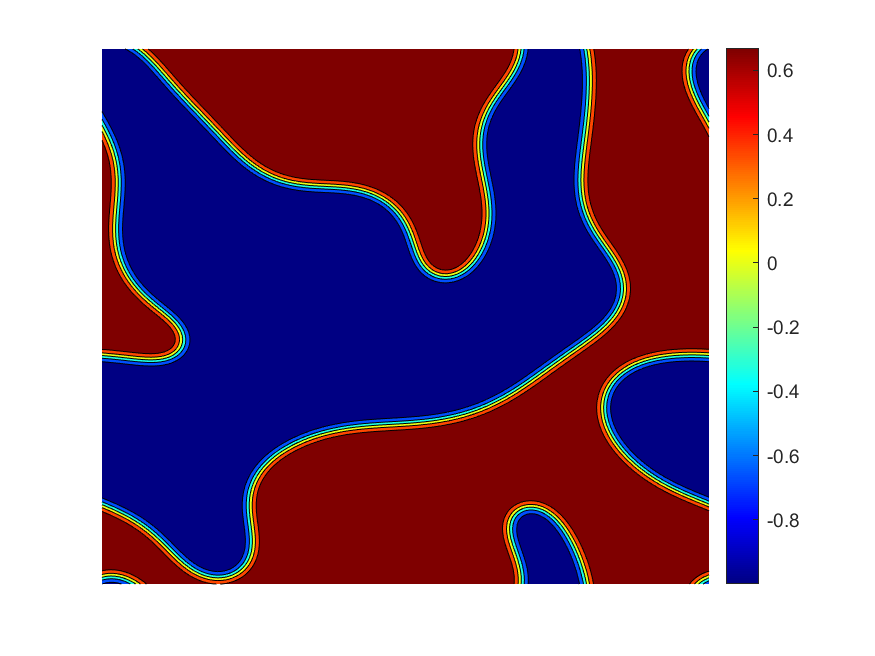}
		\end{minipage}%
	}%
	\subfigure{
		\begin{minipage}[t]{0.2\linewidth}
			\centering
			\includegraphics[width=1.1in]{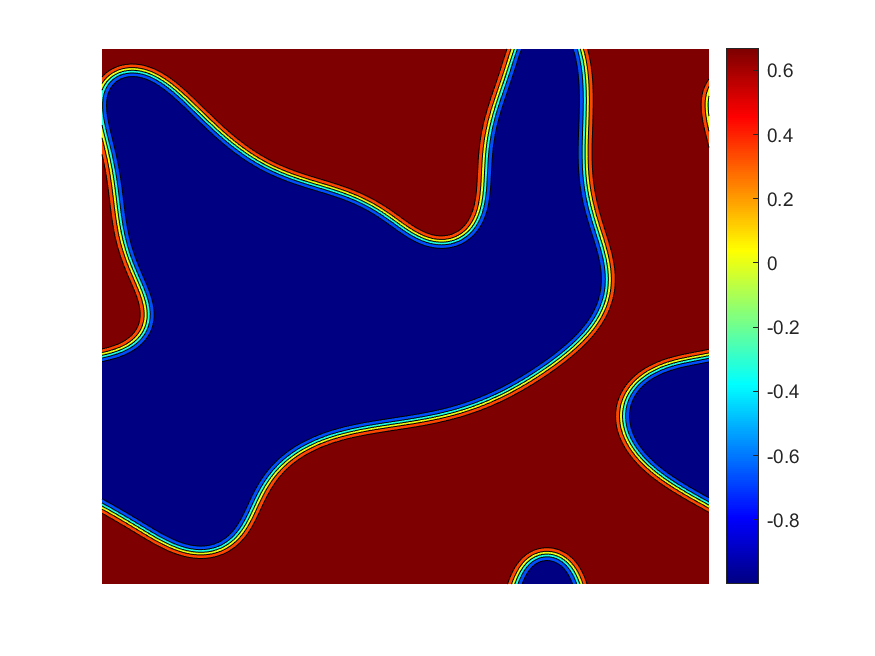}
		\end{minipage}%
	}%
	\subfigure{
		\begin{minipage}[t]{0.2\linewidth}
			\centering
			\includegraphics[width=1.1in]{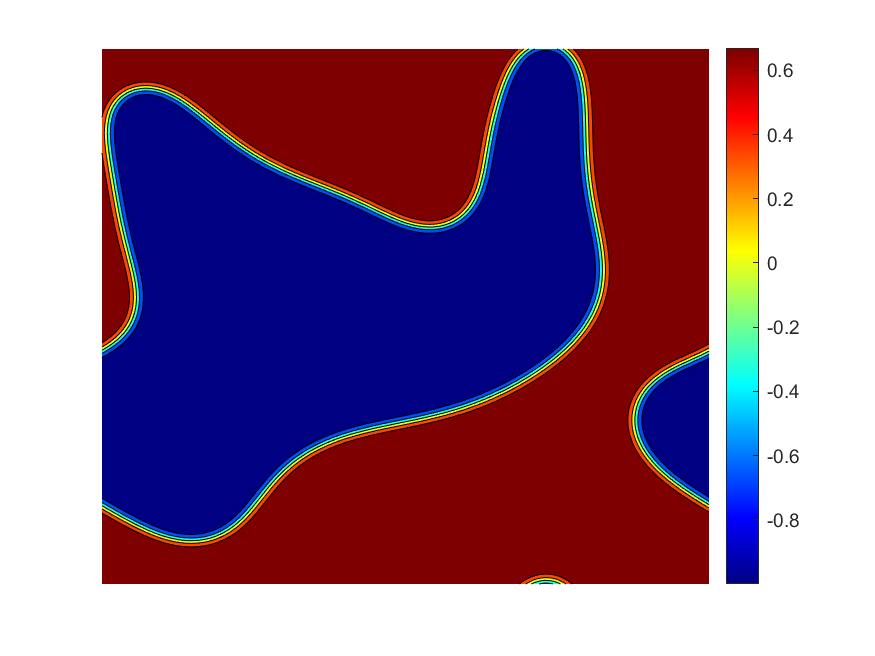}
		\end{minipage}%
	}%
	\vspace{-15pt}
	\caption{Nonlinear scheme with fixed temporal stepsize $\tau = 0.01$}
	
	\subfigure{
		\begin{minipage}[t]{0.2\linewidth}
			\centering
			\includegraphics[width=1.1in]{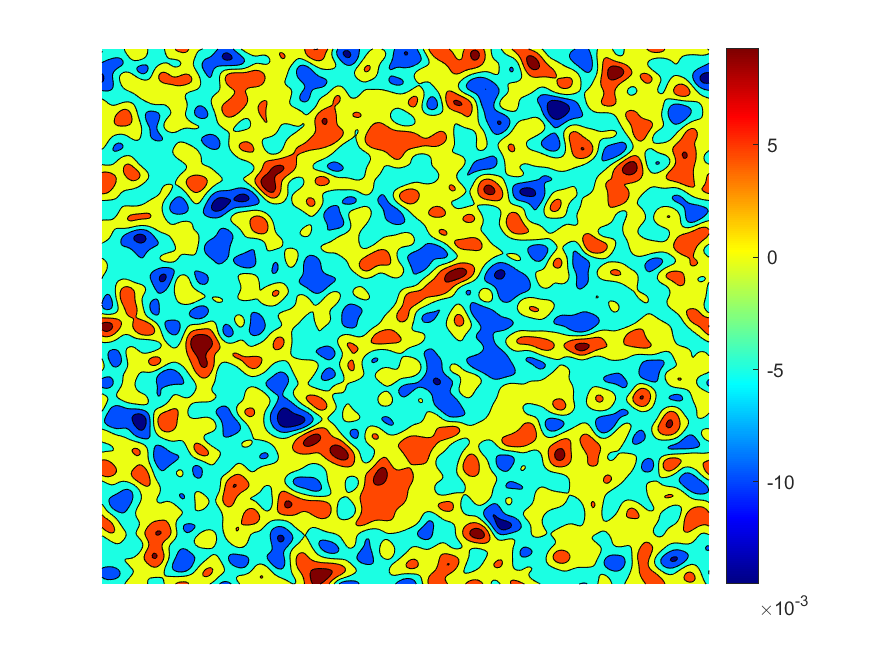}
		\end{minipage}%
	}%
	\subfigure{
		\begin{minipage}[t]{0.2\linewidth}
			\centering
			\includegraphics[width=1.1in]{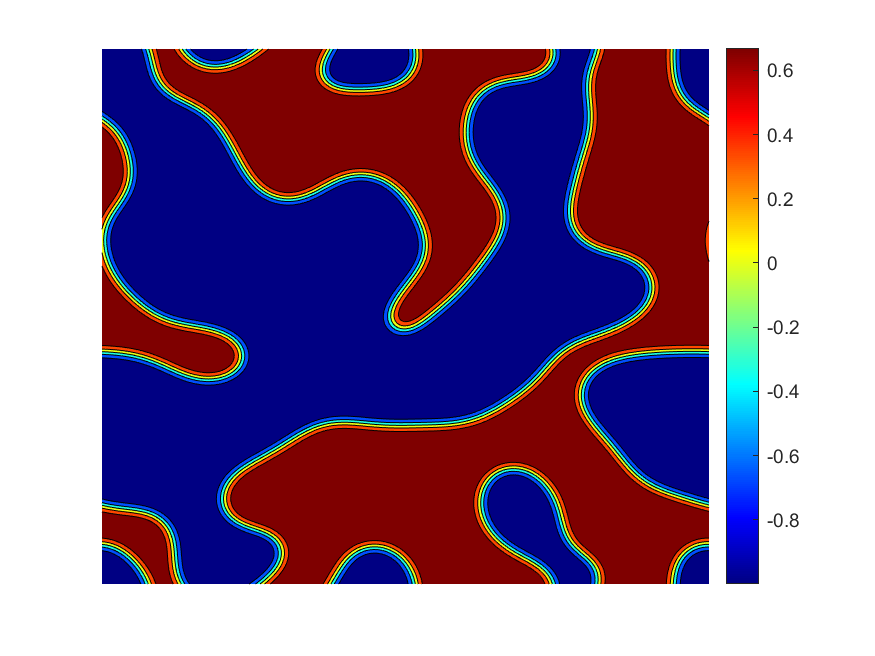}
		\end{minipage}%
	}%
	\subfigure{
		\begin{minipage}[t]{0.2\linewidth}
			\centering
			\includegraphics[width=1.1in]{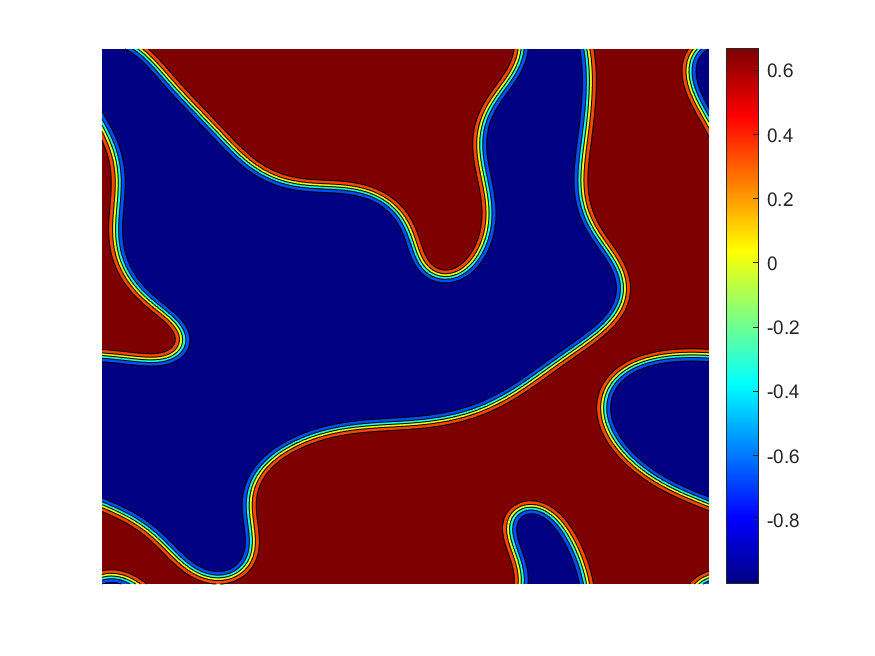}
		\end{minipage}%
	}%
	\subfigure{
		\begin{minipage}[t]{0.2\linewidth}
			\centering
			\includegraphics[width=1.1in]{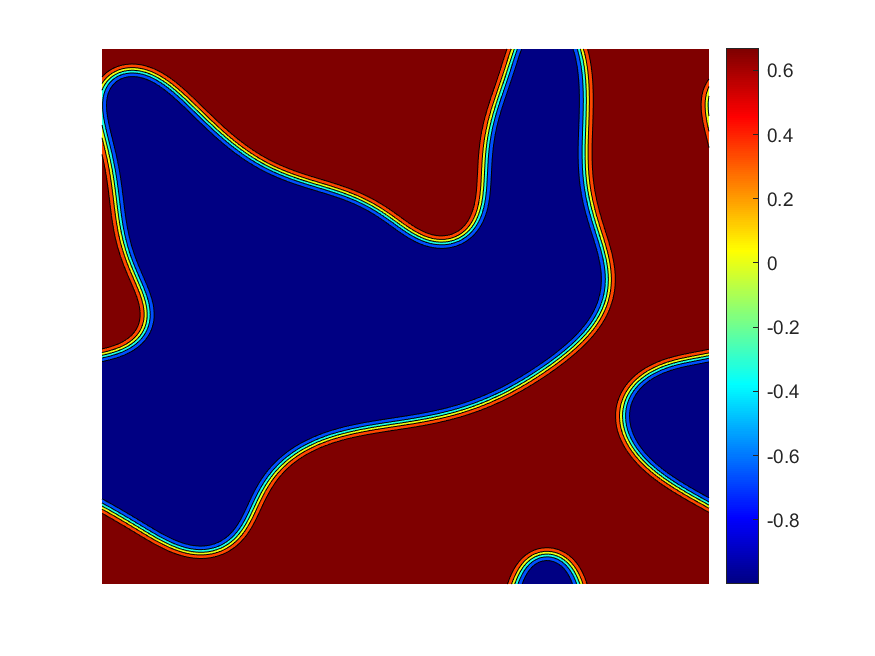}
		\end{minipage}%
	}%
	\subfigure{
		\begin{minipage}[t]{0.2\linewidth}
			\centering
			\includegraphics[width=1.1in]{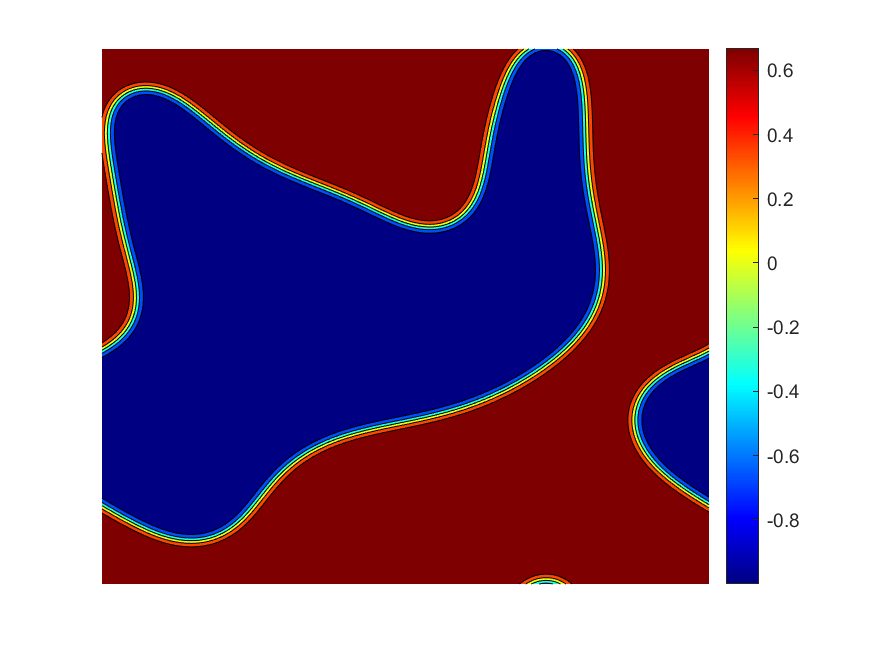}
		\end{minipage}%
	}%
	\vspace{-15pt}
	\caption{Two-grid scheme with fixed temporal stepsize $\tau = 0.01$}
	
	\subfigure{
		\begin{minipage}[t]{0.2\linewidth}
			\centering
			\includegraphics[width=1.1in]{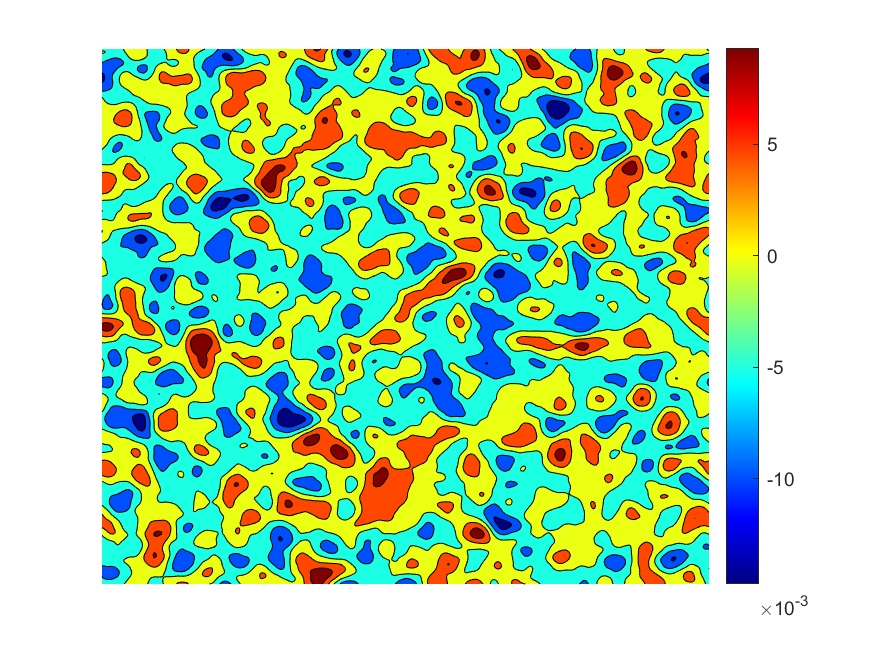}
		\end{minipage}%
	}%
	\subfigure{
		\begin{minipage}[t]{0.2\linewidth}
			\centering
			\includegraphics[width=1.1in]{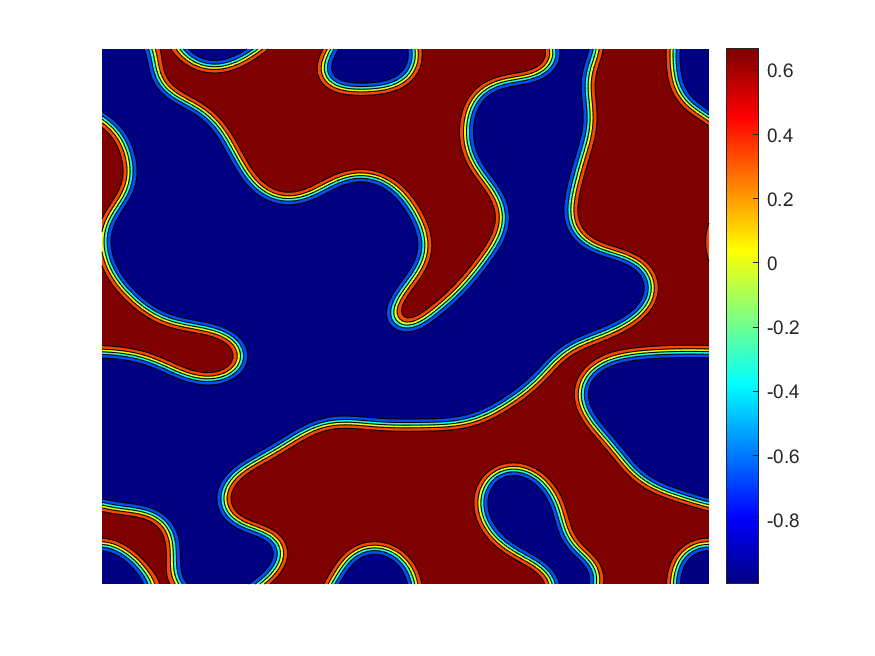}
		\end{minipage}%
	}%
	\subfigure{
		\begin{minipage}[t]{0.2\linewidth}
			\centering
			\includegraphics[width=1.1in]{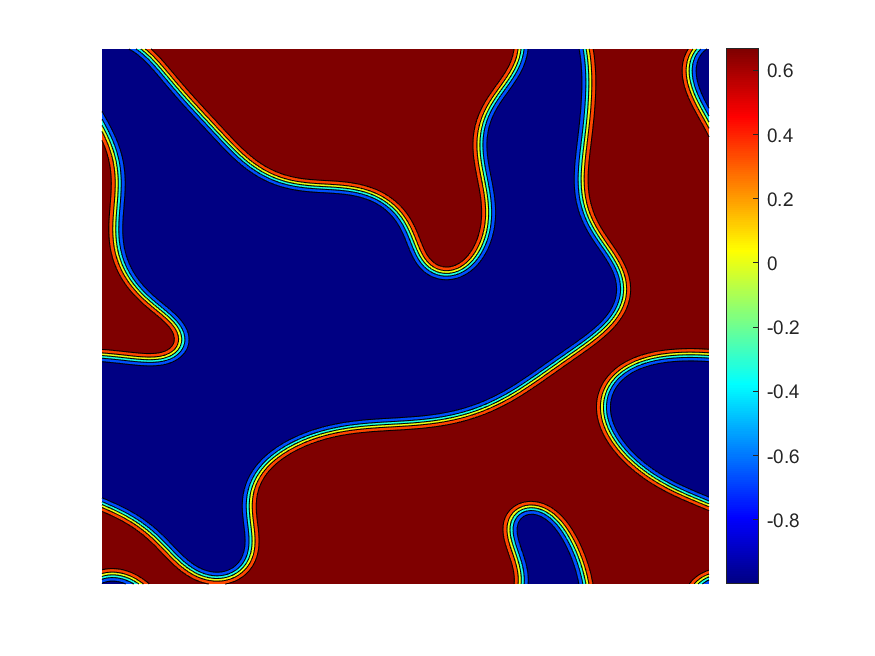}
		\end{minipage}%
	}%
	\subfigure{
		\begin{minipage}[t]{0.2\linewidth}
			\centering
			\includegraphics[width=1.1in]{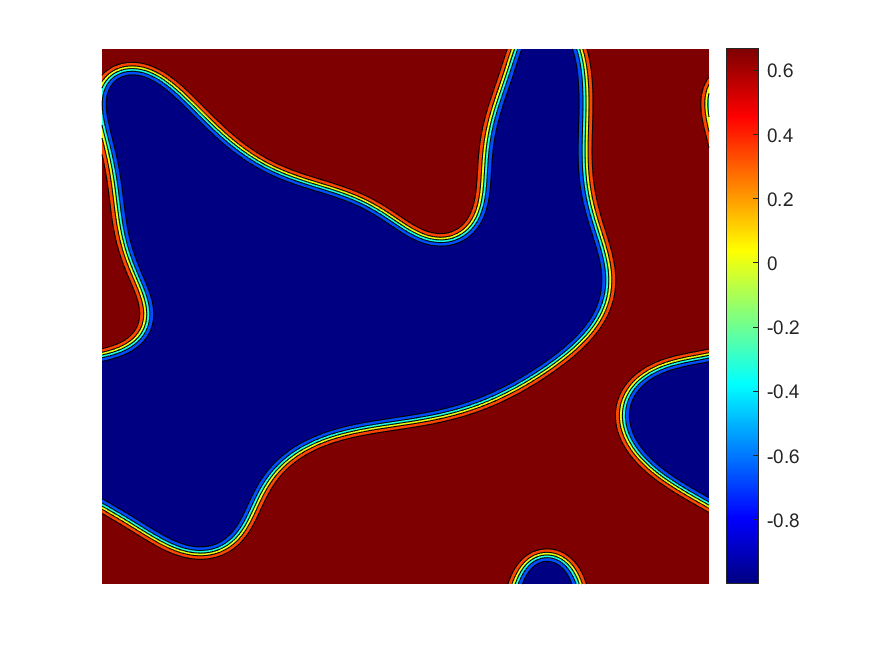}
		\end{minipage}%
	}%
	\subfigure{
		\begin{minipage}[t]{0.2\linewidth}
			\centering
			\includegraphics[width=1.1in]{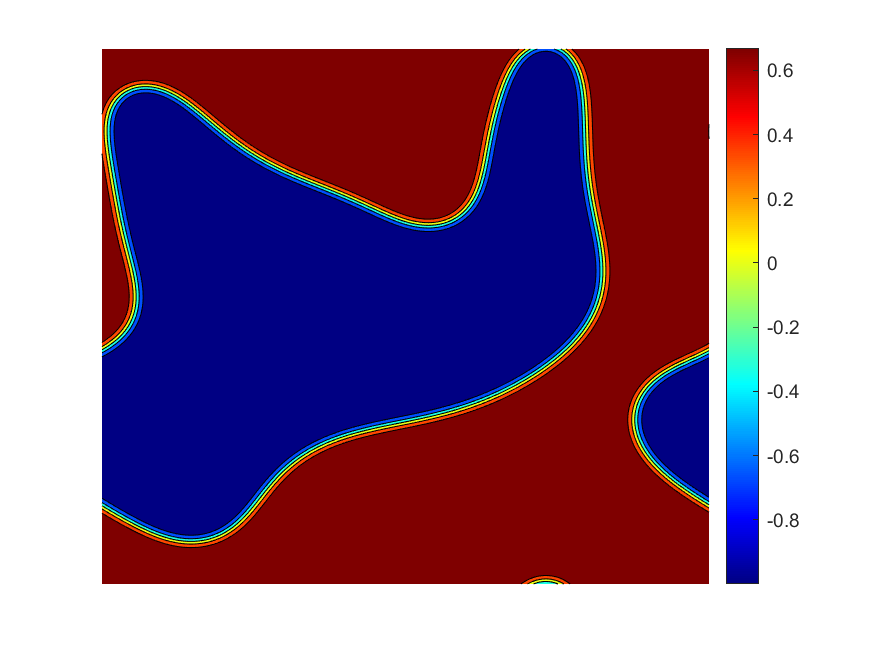}
		\end{minipage}%
	}%
	\vspace{-15pt}
	\caption{Adaptive two-grid scheme with $ \tau_{\min} = 0.01 $, $ \tau_{\max} = 1 $ and $ \eta = 8 \times 10^{4} $}
	
	\centering
	\setcounter{figure}{4}
	\renewcommand{\figurename}{Figure}
	\renewcommand{\thefigure}{\arabic{figure}}
	\caption{ Solution snapshots of coarsening dynamics for Allen--Cahn equation at $ t = 1, 20, 50, 80, 100 $ (from left to right) yielded by nonlinear scheme, two-grid scheme and adaptive two-grid scheme }
	\label{fig6}
\end{figure}

\begin{figure}[htbp]
	\centering
	\subfigure{
		\begin{minipage}[t]{0.5\linewidth}
			\centering
			\includegraphics[width=2.75in]{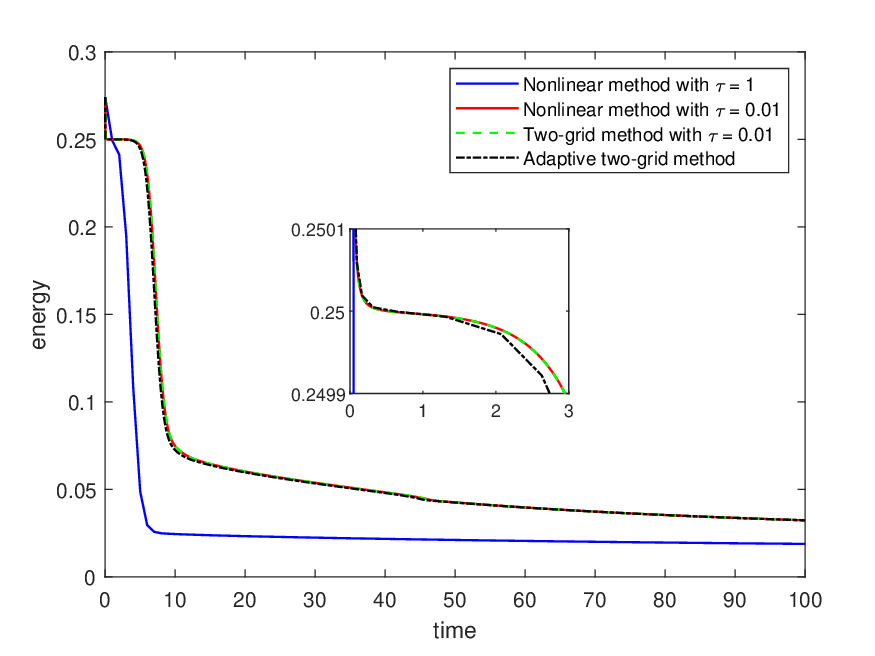}
		\end{minipage}%
	}%
	\subfigure{
		\begin{minipage}[t]{0.5\linewidth}
			\centering
			\includegraphics[width=2.75in]{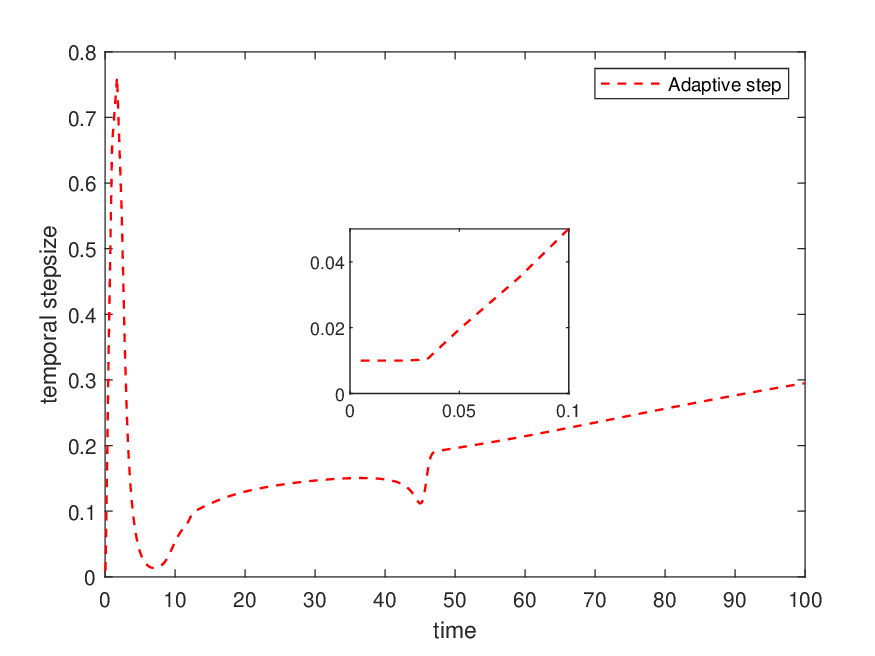}
		\end{minipage}%
	}%
	
	\centering
	\caption{ Evolutions of energy (left) and time steps (right) for the nonlinear scheme, two-grid scheme and adaptive two-grid scheme until time $ T = 100 $}
	\label{fig7}
\end{figure}

\begin{table}[!htbp]
	\caption{CPU times and the total number of temporal steps for three schemes\label{tab10}}%
	{\footnotesize\begin{tabular*}{\columnwidth}{@{\extracolsep\fill}ccccccc@{\extracolsep\fill}}
			\toprule
			& \multicolumn{2}{c}{Nonlinear scheme with $ \tau = 0.01 $}  & \multicolumn{2}{c}{Two-grid scheme with $ \tau = 0.01 $} & \multicolumn{2}{c}{Adaptive two-grid scheme}\\
			\midrule
			T & $ N $ & CPU times  & $ N $  & CPU times & $ N $ & CPU times \\
			\midrule 
			20  &  2000    &   4 h 29 m 59 s      &  2000    &  2 h 12 m 22 s    &  384  & 19 m 53 s       \\ 
			50  &  5000    &   10 h 33 m 53 s      &  5000    &  4 h 38 m 8 s    &  588  & 27 m 40 s       \\ 
			80  &  8000    &   15 h 41 m 54 s      &  8000    &  6 h 41 m 51 s    &  722   & 31 m 35 s       \\ 
			100  &  10000    &   18 h 43 m 29 s      &  10000    &  8 h 4 m 26 s    &  795  & 32 m 57 s      \\ 
			\bottomrule
	\end{tabular*}}
\end{table}

This simulation is performed under $ N_x^{h} = N_y^{h} = 384 $ and $ M_{x} = M_{y} = 3 $. Due to the fact that the initial values are randomly given, we provide the startup values on coarse grid for the two-grid scheme by implementing the nonlinear algorithm up to $ T = 0.5 $. In Figure \ref{fig6}, it displays the evolution of the coarsening dynamic for nonlinear and two-grid compact schemes with different time strategies. It is observed that the nonlinear scheme with large uniform temporal stepsize $ \tau = 1 $ yields inaccurate solution $ u $, while the adaptive two-grid scheme gives the correct coarsening pattern which is consistent with the results obtained by the nonlinear and two-grid methods with small uniform temporal stepsize $ \tau = 0.01 $. In Figure \ref{fig7}, we depict the evolution of discrete energies and temporal stepsizes with respect to time, which shows that the energy dissipation of the adaptive two-grid method agrees very well with the the nonlinear and two-grid methods using small uniform temporal stepsize. Moreover, the efficiency of the proposed two-grid scheme using the adaptive temporal stepsize strategy can also be seen from Figure \ref{fig7} (right) and Table \ref{tab10}. For example,  the adaptive two-grid method using variable-step temporal grid costs only 33 minutes for time marching to $T=100$, while the two-grid method using uniform temporal stepsize $\tau=0.01$ consumes more than 8 hours, even worse the implementation of the nonlinear scheme runs nearly 19 hours.

\section{Concluding remarks}\label{sec:Conclu}
High-order two-grid difference scheme for nonlinear PDEs are rarely studied in existing literature due to, e.g., the lack of the appropriate accuracy-preserving mapping operator. To address this issue, we introduce a piecewise bi-cubic Lagrange interpolation operator between two grids, and discuss its boundedness under $ L^{2} $ and $ L^{\infty} $ norms. Moreover, to effectively solve the nonlinear PDEs whose solutions may admit multiple time scales, variable-step temporal discretization methods, in particular, the variable-step multistep methods are natrually and valuable to improve accuracy for stiff problems. However, its numerical analysis is much more challenging than the single-step methods. As an illustration, combined with the variable-step BDF2 scheme, an efficient high-order two-grid difference method is developed for the semilinear parabolic equation with Dirichlet or periodic boundary conditions. The unique solvability of the nonlinear problem on coarse grid is shown by Browder's fixed point theorem. Moreover, with the help of DOC kernels and the boundedness of the high-order mapping operator, optimal-order error estimates for the two-grid method on both coarse and fine grids are rigourously proved under $ r_{k} := \tau_{k}/\tau_{k-1} < 4.8645 $ and the maximum temporal stepsize condition $ \tau = o(H^{\frac{1}{2}}) $, where the cut-off technique is used to reduce the regularity requirement on the nonlinear term $f$ to the local Lipschitz continuous condition. Several numerical examples are carried out to confirm the theoretical findings.

\section*{CRediT authorship contribution statement}
\textbf{Bingyin Zhang}: Methodology, Formal analysis, Software, Writing- Original draft.
\textbf{Hongfei Fu}: Conceptualization, Supervision, Writing- Reviewing and Editing, Funding acquisition.

\section*{Declaration of Competing Interest}
The authors declare that they have no competing interests.



 \bibliographystyle{spmpsci}      
\bibliography{Ref}

\end{document}